\documentclass{elsarticle}
\usepackage[a4paper]{geometry}
\usepackage{amsmath}
\usepackage{amssymb}
\usepackage{stmaryrd}

\newcommand{\re}{\mathbb{R}}
\newcommand{\jump}[1]{\left\llbracket #1 \right\rrbracket}
\newcommand{\avg}[1]{\left\{\!\!\left\{#1\right\}\!\!\right\}}

\newcommand{\ud}{\textrm{d}}
\newcommand{\cip}{C_\textrm{ip}}

\newcommand{\Gi}{\Gamma_{\mathcal{I}}}
\newcommand{\mom}[1]{ \left\langle #1 \right\rangle }
\newcommand{\tr}{\tau_R}
\newcommand{\fce}{f^{\textrm{ce}}}
\newcommand{\jmh}{{j-\half}}
\newcommand{\jph}{{j+\half}}
\newcommand{\Th}{\mathcal{T}_h}

\newcommand{\cfl}{\textrm{CFL}}
\newcommand{\half}{\frac{1}{2}}

\journal{}

\begin{document}

\title{Discontinuous Galerkin method for Navier-Stokes equations \\
using kinetic flux vector splitting}

\author{Praveen Chandrashekar}
\ead{praveen@math.tifrbng.res.in}
\address{TIFR Centre for Applicable Mathematics, Bangalore-560065, INDIA}

\begin{abstract}
Kinetic schemes for compressible flow of gases are constructed by exploiting the connection between Boltzmann equation and the Navier-Stokes equations. This connection allows us to construct a flux splitting for the Navier-Stokes equations based on the direction of molecular motion from which a numerical flux can be obtained. The naive use of such a numerical flux function in a discontinuous Galerkin (DG) discretization leads to an unstable scheme in the viscous dominated case. Stable schemes are constructed by adding additional terms either in a symmetric or non-symmetric manner which are motivated by the DG schemes for elliptic equations. The novelty of the present scheme is the use of kinetic fluxes to construct the stabilization terms. In the symmetric case, interior penalty terms have to be added for stability and the resulting schemes give optimal convergence rates in numerical experiments. The non-symmetric schemes lead to a cell energy/entropy inequality but exhibit sub-optimal convergence rates. These properties are studied by applying the schemes to a scalar convection-diffusion equation and the 1-D compressible Navier-Stokes equations. In the case of Navier-Stokes equations, entropy variables are used to construct stable schemes.
\end{abstract}

\maketitle
\section{Introduction}

Kinetic schemes for the compressible Navier-Stokes equations are constructed by exploiting the connection between the Boltzmann equation and the Navier-Stokes equations. The primary unknown in the Boltzmann equation is the velocity distribution function, which gives the distribution of molecules in velocity and physical space~\cite{cercignani1988boltzmann}. The macroscopic behaviour of the flow is obtained by averaging the Boltzmann equation over all the molecular velocities. In kinetic schemes based on finite volume or discontinuous finite elements, the intercell fluxes are approximated by exploiting this connection. Kinetic schemes were first developed for the solution of Euler equations governing inviscid compressible flows by making use of the Maxwell-Boltzmann distribution function, which is the equilibrium distribution for gases. While there are kinetic schemes based on a particle model which update the velocity distribution function, our interest in the present work is the use of kinetic model to construct numerical schemes at the continuum level. The kinetic model is utilized to derive a numerical flux function by using the upwinding principle with respect to the molecular velocities. The availability of such kinetic flux based schemes is useful for coupling the continuum solver with a particle-based method like direct simulation Monte Carlo in the case of rarefied flows where a pure continuum model would not be valid throughout the flow domain. The use of a DG scheme is expected to be helpful in coupling the continuum solver with the particle solver. Kinetic schemes at the continuum level based on finite volume discretization have been constructed in~\cite{DI1980231,smd-nasatp-2613}. In the approach of Mandal and Deshpande~\cite{Mandal1994447}, the two states on either side of a cell face are used to define two equilibrium distributions from which the inter-cell flux can be computed by taking the appropriate moments. This gives a numerical flux function which depends smoothly on the two states and can be used in a method of lines approach; the resulting set of ODE are solved using any time integration scheme like the Runge-Kutta method or an implicit scheme. Later, the Chapman-Enskog approximation to the BGK model~\cite{PhysRev.94.511} of the Boltzmann equation was used by Chou and Baganoff~\cite{Chou:1997:KFS:254115.254120} to derive a numerical flux function for the Navier-Stokes equations. The usual practice in finite volume methods is to use a Godunov scheme or approximate Riemann solver to construct the numerical flux function for the inviscid part of the flux and use a central difference-type scheme for the viscous part of the flux. But the use of a kinetic formulation gives a numerical flux function for the total flux, not just the inviscid flux. This numerical flux depends not only on the two states at a cell face but also on their derivatives, which are present due to the diffusion terms.

An alternate approach is followed in the gas kinetic finite volume schemes~\cite{Prendergast:1993:NHG:167236.167246,xu_martinelli_jameson_1995,Xu2001289} which are based on the Boltzmann equation with BGK collision model. The exact solution of the BGK model is used to compute the time integral of the inter-cell flux in a finite volume method. Thus in this approach, the time and space discretizations are coupled, similar to the Lax-Wendroff scheme. The use of the BGK model automatically gives a scheme for the Navier-Stokes equations. Ohwada~\cite{Ohwada2002156} has analyzed the consistency of the gas kinetic scheme, while a thorough analysis of the stability and accuracy of the BGK scheme for a scalar convection-diffusion equation has been made in~\cite{torr_xu}.

Discontinuous Galerkin finite element methods were first proposed in~\cite{reedhill} for the neutron transport equation. These methods use a discontinuous representation of the solution which allows them to compute solutions with steep gradients and shocks. They are useful for solving hyperbolic problems like the Euler equations of inviscid compressible flows and convection dominated problems like the compressible Navier-Stokes equations at high Reynolds numbers. The inter-cell discontinuity in the state is resolved by making use of any consistent numerical flux function for the convective flux. The vast developments in Godunov and approximate Riemann solvers can be exploited to construct these numerical flux functions. For the case of hyperbolic problems like Euler equations, the Runge-Kutta discontinuous Galerkin methods have been extensively developed~\cite{Cockburn:1989:TRL:69978.69982,Cockburn:1998:RDG:287244.287254} and their high order of accuracy is demonstrated in numerical experiments. DG methods provide a rational approach to construct high order accurate schemes for convection-dominated problems and can be considered as a generalization of finite volume methods, which have been very successful in computing discontinuous solutions and convection dominated, high Reynolds number flows. They also have many computational advantages like a compact stencil which aids parallelization, and the ability to easily develop hp-adaptive algorithms.

In the case of Navier-Stokes equations, the viscous fluxes which involve derivatives of the solution, must also be computed across the element boundaries. However, the derivatives are also possibly discontinuous across the elements due to the discontinuous numerical approximation. Averaging the derivatives to compute the inter-element flux leads to an unstable scheme, and is a generic problem for a DG scheme applied to any partial differential equation with elliptic terms. There are two approaches to construct stable schemes for the diffusion terms, which are sometimes refered to as the {\em flux formulation} and the {\em primal formulation}~\cite{arnold2002}. In the first approach, additional variables are introduced for the derivatives leading to a mixed finite element formulation. Some numerical flux function is required for the additional variables whose choice is dictated by the stability of the resulting scheme. Bassi and Rebay~\cite{Bassi1997267} developed such a scheme for compressible Navier-Stokes equations and applied it to laminar flow problems. Cockburn and Shu~\cite{Cockburn:1998:LDG:305653.305671} developed the local discontinuous Galerkin (LDG) method for time dependent nonlinear convection-diffusion problems by generalizing the method of~\cite{Bassi1997267}. In the second approach, the derivatives are averaged for computing the diffusive fluxes but additional terms are added to achieve stability of the scheme including interior penalty terms which help to enforce continuity of the solution across the elements. The resulting schemes, known as interior penalty scheme, may be symmetric (SIPG) or non-symmetric (NIPG) with respect to the diffusive terms. Arnold~\cite{arnold1982} proposed a symmetric DG scheme with interior penalty for the heat equation and showed optimal error  estimates in $L_2$ norm. Baumann and Oden~\cite{Oden1998491,Baumann1999311} introduced a non-symmetric scheme for diffusion and convection-diffusion problems, which has been used for compressible Navier-Stokes equations in~\cite{FLD:FLD338}. Symmetric schemes are preferable since they are the only ones which have adjoint consistency~\cite{arnold2002}, which is crucial for obtaining optimal error estimates in $L_2$ norm. A symmetric interior penalty scheme for compressible Navier-Stokes equations has been proposed in~\cite{Hartmann20089670} and its optimal accuracy has been shown through numerical experiments. A symmetric space-time DG scheme for Navier-Stokes equations using the appoach of~\cite{Bassi1997267} and diffusive Riemann solvers for the viscous terms has been developed in~\cite{Gassner:2008:DGS:1342061.1342078} and its optimal convergence rate has been demonstrated by numerical studies. The various DG formulations for an elliptic equation have been unified in~\cite{arnold2002} and their stability, consistency and adjoint consistency have been analyzed.

The use of a kinetic scheme has the attractive feature that it provides a numerical flux function for the total flux which can be used in a DG scheme. The numerical flux depends on the derivatives of the solution which are possibly discontinuous and a naive DG discretization can lead to an unstable scheme, especially in the viscous dominated situation. Xu~\cite{Xu:2004:DGB:996978.997050} has developed a DG discretization for the  Navier-Stokes equations using the BGK schemes, which however gives sub-optimal convergence rates similar to the NIPG scheme. Liu and Xu~\cite{springerlink:10.1007/BF03177419} proposed a modified version of discontinuous Galerkin BGK scheme for Navier-Stokes equations which shows improved convergence rates. These BGK-based DG schemes have not yet been analyzed from the point of view of stability and adjoint consistent as in~\cite{arnold2002}. In the present work, we follow the interior penalty approach and construct stable schemes by making use of the kinetic split fluxes to add additional stabilization terms similar to the NIPG and SIPG schemes. The kinetic model is used to obtain a flux splitting from which a numerical flux function is constructed even for the diffusion terms. The development of the schemes is motivated by considering the scalar convection-diffusion equation for which a kinetic formulation can be given. Through numerical experiments on this model problem, we find the NIPG scheme gives a convergence rate of $O(h^k)$ when using $P_k$ basis functions which is sub-optimal as compared to the usual Galerkin methods, while the SIPG scheme gives the optimal $O(h^{k+1})$ convergence rate. Similar convergence rates are observed through numerical experiments for the 1-D compressible Navier-Stokes equations. In the case of the Navier-Stokes equations, it is advantageous to use {\em entropy variables} since the equations assume a symmetric form in those variables~\cite{Amiram1983151,Shakib1991141,smd1986a}. It also facilitates the construction of entropy stable schemes which is an important thermodynamic constraint for compressible flows. Even the pure Galerkin scheme written in entropy variables satisfies the entropy condition~\cite{Shakib1991141} for any degree of the basis functions. For Euler equations, Barth~\cite{barth1998,Timothy20063311} has studied the entropy condition for discontinuous Galerkin methods and has derived a condition for a numerical flux function to be entropy consistent, which is a generalization of the E-flux condition~\cite{osher:947} to systems of equations. The numerical flux resulting from kinetic flux vector splitting  scheme~\cite{Mandal1994447} satisfies this entropy condition. The entropy variable formulation is also advantageous since the flux jacobians are well behaved in the limit of incompressible flows~\cite{Hauke:1998:0045-7825:1}; this is not the case with conserved variable formulation where some terms in the flux jacobian matrix become unbounded or assume indeterminate form.

The rest of the paper is organized as follows. In Section~(\ref{sec:con}) we first consider the scalar convection equation for which a kinetic formulation is given together with a DG scheme. The energy stability of the DG scheme is shown for the time continuous case. In Section~(\ref{sec:cd}), a kinetic formulation for the convection-diffusion equation is given together with DG discretization. The instability of the naive DG discretization is shown through numerical experiments. Stabilized versions of the schemes are constructed and their order of accuracy is studied through numerical experiments. In Section~(\ref{sec:ns}), the Navier-Stokes equations and their entropy variables are discussed. A novel DG scheme for the Navier-Stokes equations based on kinetic flux vector splitting is presented together with some numerical results which demonstrate the accuracy of the scheme.

\section{Kinetic model for convection equation}
\label{sec:con}

The scalar convection equation
\begin{equation}
u_{t}+cu_{x}=0\label{eq:lc}
\end{equation}
is the simplest hyperbolic partial differential equation. Throughout this work, we will assume periodic boundary conditions for the above partial differential equation. The exact solution of this equation is obtained by convecting the initial condition at a constant speed $c$ and the shape of the solution profile is unchanged. We can obtain this equation as moments of a "Maxwell-Boltzmann" distribution function
\begin{equation}
g(v,u)=u\sqrt{\frac{\beta}{\pi}}\exp\left[-\beta(v-c)^{2}\right]\label{eq:bolt}
\end{equation}
since $\langle g\rangle=u$ and $\langle vg\rangle=cu$, where the angle brackets denote integration over velocity space
\[
\langle\cdot\rangle=\int_{-\infty}^{+\infty}(\cdot)\textrm{d}v
\]
The quantity $\beta>0$ is a free parameter that can be chosen arbitrarily\footnote{In the case of compressible gases for which kinetic theory holds, there is no free parameter in the Maxwell-Boltzmann velocity distribution function. A parameter like $\beta$ is present but is related to the absolute temperature.}.
The "Boltzmann equation" for the distribution function is obtained by
differentiating equation \eqref{eq:bolt}
\[
g_{t}+vg_{x}=(u_{t}+vu_{x})\sqrt{\frac{\beta}{\pi}}\exp\left[-\beta(v-c)^{2}\right]
\]
and we use equation~\eqref{eq:lc} to replace $u_{t}$ to get
\begin{equation}
g_{t}+vg_{x}=(v-c)u_{x}\sqrt{\frac{\beta}{\pi}}\exp\left[-\beta(v-c)^{2}\right]=:Q\label{eq:bolt1}
\end{equation}
Note that the moment of the collision term $Q$ is zero, $\langle Q\rangle=0$,
so that taking moments of equation \eqref{eq:bolt1} gives us the
convection equation \eqref{eq:lc}. The vanishing of the moments of the collision term means that molecular collisions instantaneously drive the distribution function towards the local equilibrium distribution.

\subsection{Kinetic numerical flux: KFVS}

\begin{figure}
\begin{center}
\includegraphics[width=0.4\textwidth]{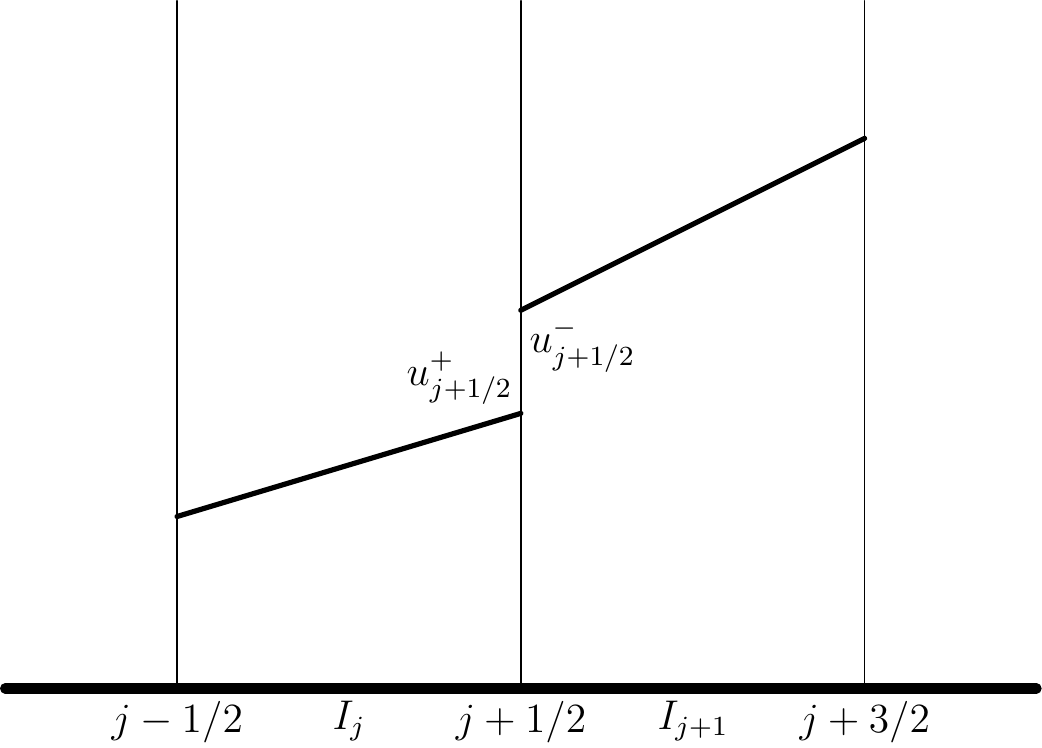}
\caption{Discontinuous finite element solution}
\label{fig:discsol}
\end{center}
\end{figure}
In a finite volume or DG method, we need to compute the flux across
the inter-element boundaries. The solution is possibly discontinuous
across the elements. This defines two states $u_{j+\half}^{+}$, $u_{j+\half}^{-}$
on either side of the cell face $x_{j+\half}$ as shown in figure~(\ref{fig:discsol}), where we have used the
notation
\[
u^{\pm}(x)=\lim_{\epsilon\searrow0}u(x\mp\epsilon)
\]
In the kinetic model, molecules carry information due to their streaming motion. This allows us to use the upwind principle to approximate the velocity distribution function at $x_{j+\half}$ based on the sign of the molecular velocity
as
\[
g_{j+\half}=\begin{cases}
g(v,u_{j+\half}^{+}) & v>0\\
g(v,u_{j+\half}^{-}) & v<0
\end{cases}
\]
The numerical flux function is obtained by integrating the above distribution function over the velocity space leading to
\[
F_{j+\half} = F(u_{j+\half}^{+},u_{j+\half}^{-}) = \langle vg_{j+\half}\rangle  = F^{+}(u_{j+\half}^{+})+F^{-}(u_{j+\half}^{-})
\]
The split fluxes are obtained by integrating over half velocity spaces,
and are given by
\[
F^{\pm}(u)=cuA^{\pm}+uB^{\pm}, \quad 
A^{\pm}=\frac{1}{2}[1\pm\textrm{erf}(s)],\quad B^{\pm}=\pm\frac{1}{2\sqrt{\pi\beta}}\exp(-s^{2}),\quad s=c\sqrt{\beta}
\]
The flux $F^+$ is due to all molecules moving to the right ($v>0$) while the flux $F^-$ is due to all molecules moving to the left ($v<0$); this is the upwind principle for kinetic schemes. The numerical flux function $F(u^{-},u^{+})$ is consistent, i.e.,
$F(u,u)=cu$ and is obviously a smooth function of $u$. The split
flux $F^{+}(u)$ is an increasing function of $u$ while $F^{-}(u)$
is a decreasing function. Hence the numerical flux $F(u^{-},u^{+})$
is a monotone flux. We can also write the numerical flux function
as a central flux with dissipation
\begin{equation}
F(u^{+},u^{-})=\frac{1}{2}(cu^{+}+cu^{-})+\frac{1}{2}D(u^{+}-u^{-}),\qquad D=c\,\textrm{erf}(s)+\frac{1}{\sqrt{\pi\beta}}\exp(-s^{2})
\label{eq:lcfc}
\end{equation}
and it is easy to check that $D>0$. In the limit of $\beta\to\infty$, the split fluxes become
\[
F^{+}(u)=\begin{cases}
cu, & c>0\\
0, & c<0
\end{cases},\qquad F^{-}(u)=\begin{cases}
0, & c>0\\
cu, & c<0
\end{cases}
\]
and the numerical flux function of equation~(\ref{eq:lcfc}) becomes the standard upwind flux
\[
F(u^{+},u^{-})=\begin{cases}
cu^{+}, & c>0\\
cu^{-}, & c<0
\end{cases}
\]
We will see that the quantity $D$ controls the rate of dissipation of energy and it is the smallest for the upwind scheme. In all the computations we use a value of $\beta=1$.

\subsection{DG scheme}

Consider a discretization of the one dimensional domain $\Omega$
by $N$ elements $I_{j}=[x_\jmh,x_\jph]$, $j=1,\ldots,N$. The DG method uses the
broken space of polynomials of degree $k$
\[
V_{h}^{k}=\left\{ \phi\in L^{2}(\Omega):\phi|_{I_{j}}\in P_{k}(I_{j})\right\} 
\]
The DG scheme is obtained by multiplying equation \eqref{eq:lc} by
a test function $\phi_h$ from the broken space $V_h^k$ and integrating over element
$I_{j}$
\[
\int_{I_{j}}\phi_{h}\partial_{t}u_{h}dx-\int_{I_{j}}cu_{h}\partial_{x}\phi_{h}dx+F(u_{h})_{j+\half}\phi_{h}(x_{j+\half}^{+})-F(u_{h})_{j-\half}\phi_{h}(x_{j-\half}^{-})=0
\]
The same equation would have been obtained if we had started from equation~(\ref{eq:bolt1}) by multiplying it with a test function and perform the integrals over velocity space and physical space; the collision term would have vanished due to the integration over velocity space. Note that we have used the numerical flux function to resolve the
discontinuity at cell boundaries, i.e.,
\[
F(u_{h})_{j+\half}=F(u_{j+\half}^{+},u_{j+\half}^{-})
\]
If we take $\phi_{h}=u_{h}$ and sum over all elements, we obtain
the energy equation
\begin{equation}
\frac{1}{2}\frac{d}{dt}\|u_{h}\|^{2}-\sum_{j=1}^{N}\int_{I_{j}}cu_{h}\partial_{x}u_{h}dx + \sum_{j=1}^{N}F(u_{h})_{j+\half}\jump{u_{h}}_{j+\half}=0
\label{eq:ee}
\end{equation}
where we have defined the inter-element jump in the solution by
\[
\jump{u_{h}}_{j+\half}:=u_{h}(x_{j+\half}^{+})-u_{h}(x_{j+\half}^{-})
\]
After some simple calculations, the energy equation reduces to
\[
\frac{1}{2}\frac{\ud}{\ud t}\|u_{h}\|^{2}+\frac{1}{2}D\sum_{j=1}^{N}\jump{u_{h}}_{j+\half}^{2}=0
\]
The jumps at the inter-element boundaries lead to dissipation of energy
and the scheme is stable. Note that even if $D=0$, i.e., if we use a central flux, the scheme preserves the energy for any degree of the basis functions, which is consistent with the exact solution of the partial differential equation.

\section{Kinetic model for convection-diffusion equation}
\label{sec:cd}

Consider the linear convection-diffusion equation
\begin{equation}
u_{t}+cu_{x}=\mu u_{xx},\qquad\mu>0
\label{eq:lcd}
\end{equation}
We can give a kinetic formulation for this equation using
the Boltzmann equation with BGK model~\cite{torr_xu}
\begin{equation}
f_{t}+vf_{x}=\frac{g-f}{\tr}\label{eq:bgk}
\end{equation}
where $g$ is the same Maxwell-Botzmann distribution as used in the
previous section. The basic idea of the BGK model is that the distribution
function $f$ is close to the local equilibrium distribution $g$,
and due to molecular collisions, it relaxes towards $g$ over a time interval $\tr$ which is
the {\em relaxation time}. We can look for solutions to the BGK model by
using the Chapman-Enskog expansion in the relaxation time $\tr$
\begin{equation}
f=f_{0}+\tr f_{1}+\tr^{2}f_{2}+\ldots
\label{eq:ceseries}
\end{equation}
Substituting this in equation \eqref{eq:bgk} and collecting terms
with the same coefficient of $\tr$, we get
\[
f_{0}=g,\qquad f_{1}=-(g_{t}+vg_{x})
\]
as the first two terms in the expansion. The zeroth order term in the Chapman-Enskog expansion leads to the convection equation as seen in the previous sections. The first order term becomes
\begin{equation*}
f_1 = -(v-c)u_x \sqrt{\frac{\beta}{\pi}}\exp\left[-\beta(v-c)^{2}\right]
\end{equation*}
where we have used the inviscid equation $u_t = -cu_x$ to replace the time derivative term since $f_1$ should depend only on the lower order terms in the Chapman-Enskog expansion. The first two terms constitute the Chapman-Enskog distribution function
\begin{equation}
\fce=f_0 + \tr f_1 = \left[u-\tr(v-c)u_{x}\right] \sqrt{\frac{\beta}{\pi}}\exp\left[-\beta(v-c)^{2}\right]
\label{eq:fce}
\end{equation}
Upon taking moments of this distribution function, we obtain
\[
\langle \fce\rangle=u,\qquad\langle v\fce\rangle=cu-\frac{\tr}{2\beta}u_{x}
\]
In order to recover the correct flux of equation \eqref{eq:lcd}, we have to choose $\tr=2\beta\mu$ so that $\langle v\fce\rangle=cu-\mu u_{x}$. Since $\fce$ is a truncated solution in the Chapman-Enskog expansion
it does not satisfy equation \eqref{eq:bgk}. We can derive an equation
for $\fce$ by differentiating equation \eqref{eq:fce} as follows.
\[
\fce_{t}+v\fce_{x}=\left[u_{t}+vu_{x}-\tr(v-c)(u_{xt}+vu_{xx})\right] \sqrt{\frac{\beta}{\pi}}\exp\left[-\beta(v-c)^{2}\right]
\]
We now use the equation \eqref{eq:lcd} to replace $u_{t}$
and $u_{xt}$ to obtain
\[
\fce_{t}+v\fce_{x}=\left[(v-c)u_{x}+(\mu-\tr(v-c)^{2})u_{xx}-\tr\mu(v-c)u_{xxx}\right] \sqrt{\frac{\beta}{\pi}}\exp\left[-\beta(v-c)^{2}\right] =: Q
\]
Upon taking moments of this equation, the right hand side vanishes, i.e., $\mom{Q}=0$,
and we obtain the convection-diffusion equation \eqref{eq:lcd}.

\subsection{Flux vector splitting}

The split fluxes for the viscous problem are obtained by integrating
the Chapman-Enskog distribution function $\fce$ over the half velocity
spaces, $-\infty < v < 0$ and $0 < v < +\infty$. The expressions for the kinetic split fluxes can be written as
\[
F^{\pm}=(cu-\mu u_{x})A^{\pm}+uB^{\pm}
\]
The numerical flux function $F(u^+,u^-)=F^+(u^+)+F^-(u^-)$ can be written as a sum of convective and dissipative fluxes
\[
F=F_{c}(u^{+},u^{-})+F_{d}(u_{x}^{+},u_{x}^{-})
\]
where the convective flux $F_{c}$ is identical to equation \eqref{eq:lcfc},
while the dissipative flux is given by
\[
F_{d}(u_{x}^{+},u_{x}^{-})=-\mu u_{x}^{+}A^{+}-\mu u_{x}^{-}A^{-}
\]
Note that the diffusive flux can also be written in split form as
\[
F_{d}(u_{x}^{+},u_{x}^{-})=F_{d}^{+}(u_{x}^{+})+F_{d}^{-}(u_{x}^{-}),\qquad F_{d}^{\pm}(u_{x})=-\mu u_{x}A^{\pm}
\]
In the case of the convection-diffusion equation the viscous fluxes depend only on the derivatives of the solution, but for the Navier-Stokes equations, the viscous fluxes will depend on the solution also.
\subsection{Test case}
We will consider the convection-diffusion equation~(\ref{eq:lcd}) in the domain $(-1,+1)$ together with the initial condition
\begin{equation}
u(x,0) = -\sin(\pi x)
\end{equation}
and periodic boundary conditions. The exact solution to this problem is given by
\begin{equation}
u(x,t) = -\exp(-\mu \pi^2 t) \sin(\pi (x-ct) )
\end{equation}
In all the tests, we set the convection speed to be $c=1$. We define two test cases based on different values of the viscosity coefficient $\mu$ as given in table~(\ref{tab:lcdtest}) with $t_f$ being the final time. Test case 1 corresponds to a situation in which convection dominates diffusion, while in Test case 2, diffusion dominates convection effects. The time integration is performed using the 3-stage, strong stability preserving Runge-Kutta scheme of Shu and Osher~\cite{Shu1988439}. The time-step is chosen based on a CFL condition of the form $\Delta t = \cfl \cdot \max(h/|c|,h^2/\mu)$. In the case of convection equation solved with a DG scheme, the CFL number depends on the degree of the polynomial basis functions~\cite{cockburn1989}.
\begin{table}[htbp]
   \centering
   \begin{tabular}{|c|c|c|}
   \hline
   Test Case & 1 & 2 \\
   \hline
   $\mu$ & 0.001 & 1 \\
   \hline
   $t_f$ & 30 & 0.5 \\
   \hline  
   \end{tabular}
   \caption{Parameters for the test case based on convection-diffusion equation~(\ref{eq:lcd})}
   \label{tab:lcdtest}
\end{table}

\subsection{DG scheme}

\begin{figure}
\begin{center}
\begin{tabular}{cc}
\includegraphics[width=0.35\textwidth]{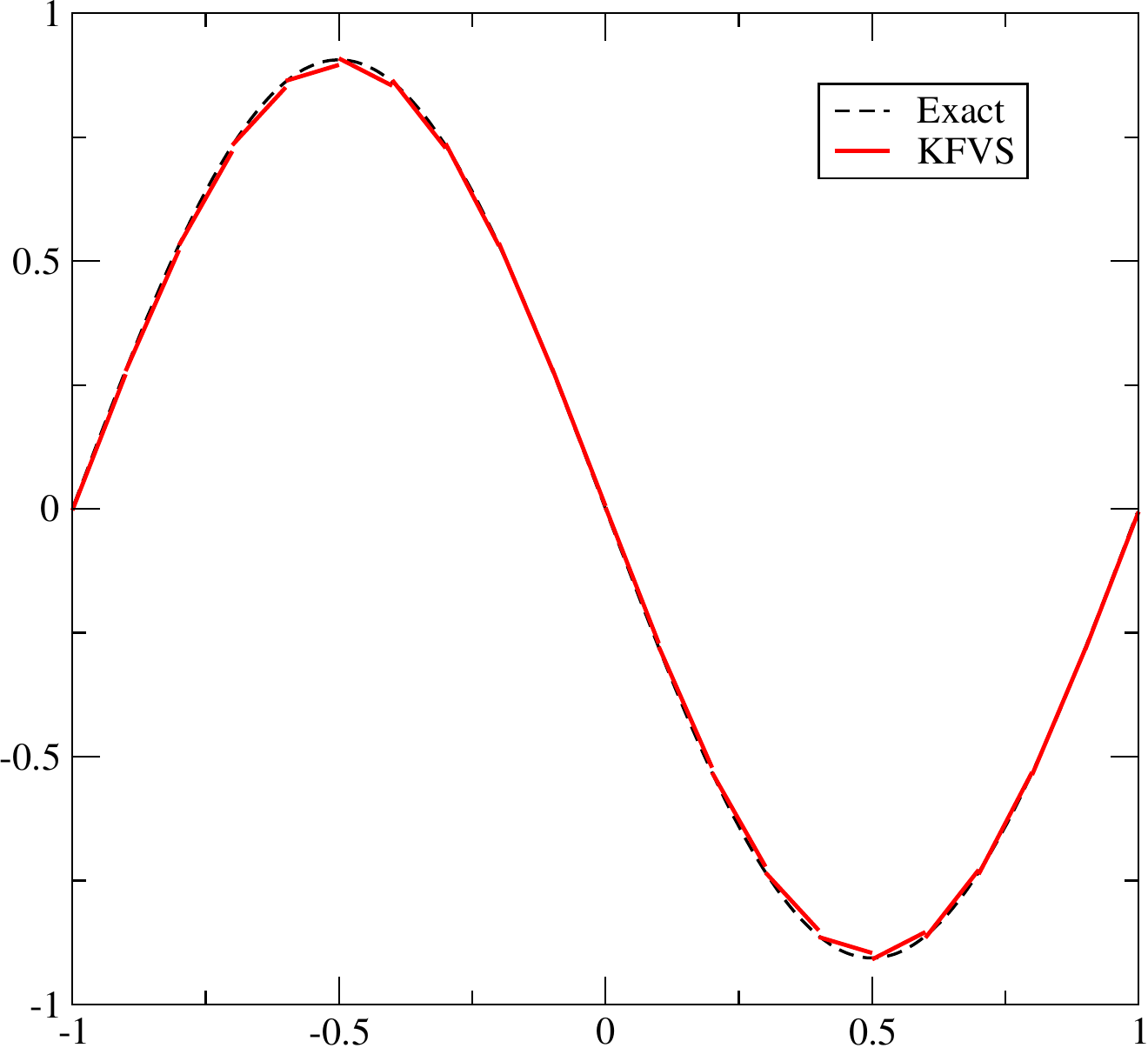} & 
\includegraphics[width=0.35\textwidth]{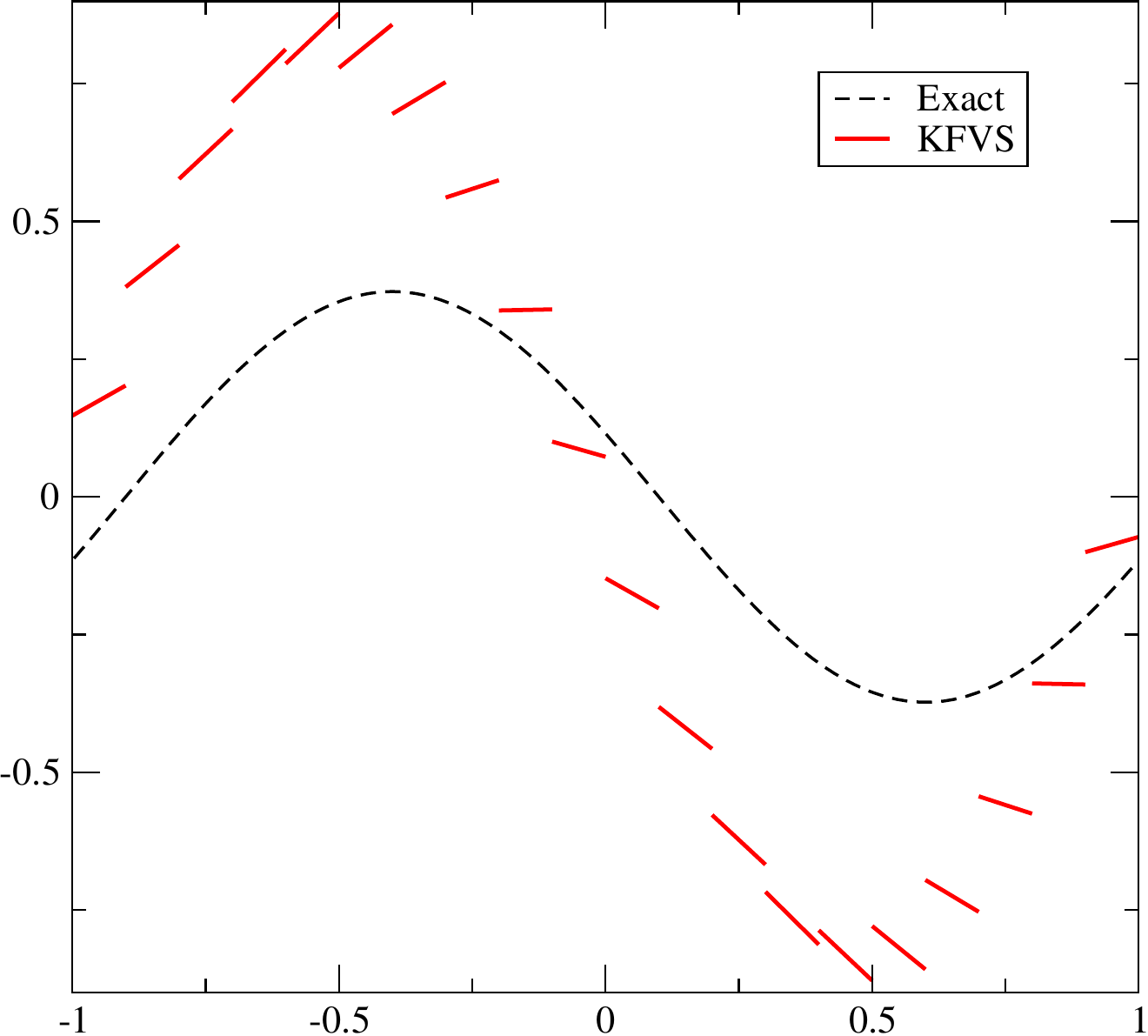} \\
(a) & (b) \\
\end{tabular}
\end{center}
\caption{Convection-diffusion problem using equation~(\ref{eq:dgelem}) and $P_{1}$ basis
functions. (a) Test case 1, (b) Test case 2}
\label{fig:lcd1}
\end{figure}

The DG discretization of equation \eqref{eq:lcd} is obtained by multiplying it by a test function $\phi_h$ and integrating over any element $I_j$
\begin{equation}
\int_{I_{j}}\phi_{h}\partial_{t}u_{h}dx  -\int_{I_{j}}cu_{h}\partial_{x}\phi_{h}dx+\mu\int_{I_{j}}\partial_{x}u_{h}\partial_{x}\phi_{h}dx
  +F(u_{h})_{j+\half}\phi_{h}(x_{j+\half}^{+})-F(u_{h})_{j-\half}\phi_{h}(x_{j-\half}^{-})=0
\label{eq:dgelem}
\end{equation}
In the above equation, the inter-element flux is approximated by the numerical flux which is defined as
\[
F(u_{h})=F_{c}(u_{h}^{+},u_{h}^{-})+F_{d}(\partial_{x}u_{h}^{+},\partial_{x}u_{h}^{-})
\]
and it includes both convective and diffusive fluxes. Substituting $\phi_{h}=u_{h}$ and summing over all elements we obtain
\begin{equation*}
\frac{1}{2}\frac{\ud}{\ud t}\|u_{h}\|^{2}-\sum_{j=1}^{N}\int_{I_{j}}cu_{h}\partial_{x}u_{h}dx  +\mu\sum_{j=1}^{N}\int_{I_{j}}(\partial_{x}u_{h})^{2}dx
  +\sum_{j=1}^{N}F(u_{h})_{j+\half} \jump{u_{h}}_{j+\half}=0
\end{equation*}
Following similar steps as in the case of linear convection equation,
we obtain the following energy equation
\begin{equation}
\frac{1}{2}\frac{\ud}{\ud t}\|u_{h}\|^{2}+\frac{1}{2}D\sum_{j=1}^{N} \jump{u_{h}}_{j+\half}^{2}+\mu\sum_{j=1}^{N}\int_{I_{j}}(\partial_{x}u_{h})^{2}dx + \sum_{j=1}^{N}F_{d}(\partial_{x}u_{h})_{j+\half} \jump{u_{h}}_{j+\half}=0
\label{eq:eprob}
\end{equation}
The second and third terms are dissipative since they lead to a decrease of the energy. But we do not know the sign of the last term, which can be either positive or negative. We apply this DG
scheme to the two test problems using $P_1$ basis functions and the results are shown in figure~(\ref{fig:lcd1}).
When $\mu=0.001$ the results look correct but this is not so when
$\mu=1$ which corresponds to a viscous dominated case. The scheme
is not able to dissipate energy at the correct rate. We see that a
naive discretization of convection-diffusion equation by the discontinuous Galerkin method leads to unstable schemes when viscosity is dominant.

\subsection{DG scheme with stabilization}

\begin{figure}
\begin{center}
\begin{tabular}{cc}
\includegraphics[width=0.35\textwidth]{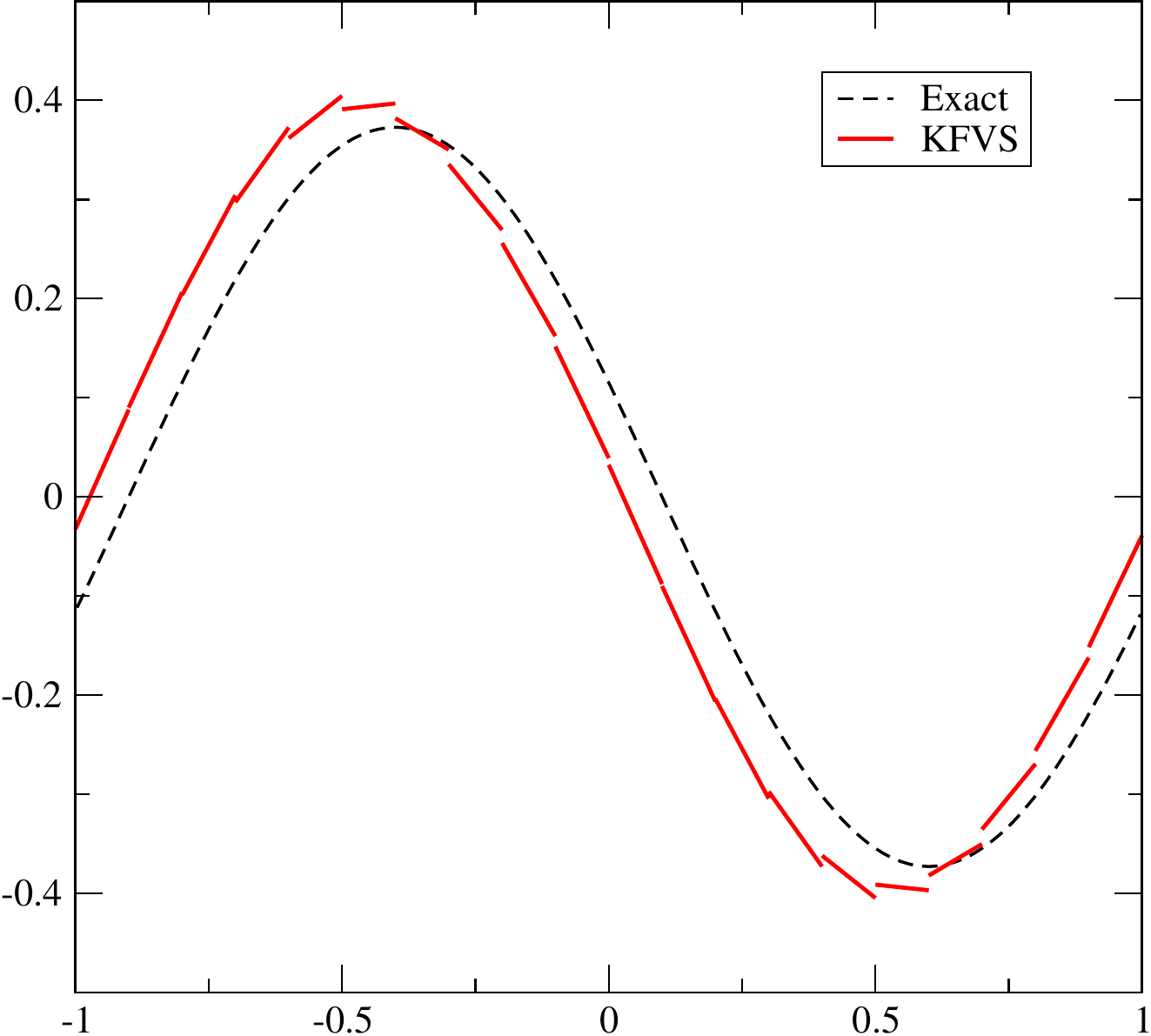} &
\includegraphics[width=0.35\textwidth]{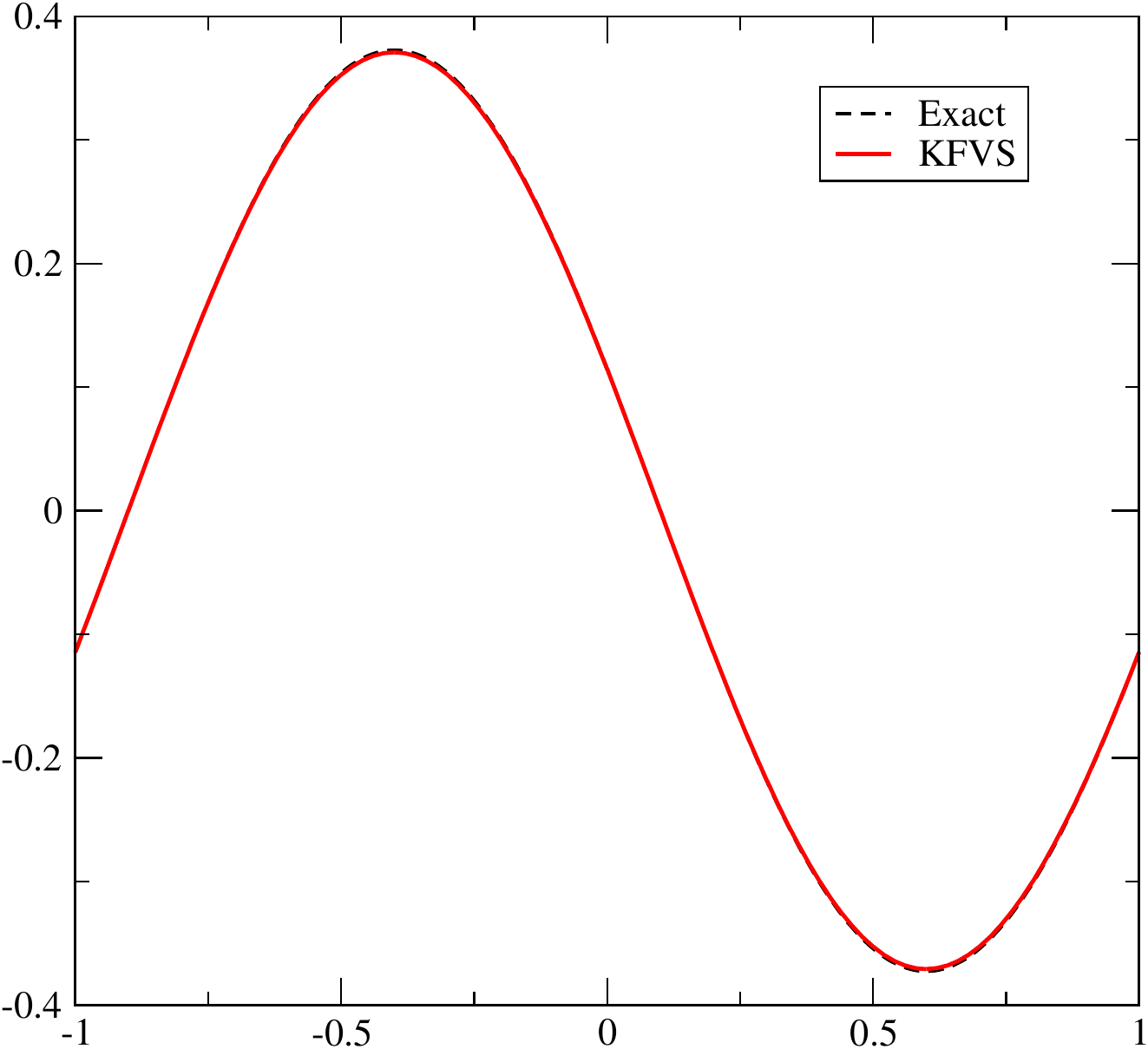} \\
(a) & (b) \\
\end{tabular}
\end{center}
\caption{Results for Test case 2 using equation
\eqref{eq:lcddg1}: (a) $P_{1}$  and (b) $P_{2}$ basis functions}
\label{fig:lcd2}
\end{figure}

Motivated by the energy analysis, we propose the following stabilized
DG scheme
\begin{equation}
\begin{aligned}\int_{I_{j}}\phi_{h}\partial_{t}u_{h}dx & -\int_{I_{j}}cu_{h}\partial_{x}\phi_{h}dx+\mu\int_{I_{j}}\partial_{x}u_{h}\partial_{x}\phi_{h}dx +F(u_{h})_{j+\half}\phi_{h}(x_{j+\half}^{+})-F(u_{h})_{j-\half}\phi_{h}(x_{j-\half}^{-})\\
 & - F_{d}^{+}(\partial_{x}\phi_{h}^{+})_{j+\half} \jump{u_{h}}_{j+\half} - F_{d}^{-}(\partial_{x}\phi_{h}^{-})_{j-\half} \jump{u_{h}}_{j-\half}=0
\end{aligned}
\label{eq:lcddg1}
\end{equation}
Note that the last two terms in the above equation are constructed from the kinetic
split viscous fluxes. Choosing $\phi_{h}=u_{h}$ and summing up over
all elements gives us
\begin{equation}
\frac{1}{2}\frac{\ud}{\ud t}\|u_{h}\|^{2}+\frac{1}{2}D\sum_{j=1}^{N} \jump{u_{h}}_{j+\half}^{2} + \mu\sum_{j=1}^{N}\int_{I_{j}}(\partial_{x}u_{h})^{2}dx=0
\label{eq:lcdge}
\end{equation}
which shows that the scheme is stable. The choice of the additional stabilization terms in equation \eqref{eq:lcddg1} was just enough to eliminate the problematic term in the energy equation~(\ref{eq:eprob}). The stabilized scheme is applied
to convection-diffusion problem with $\mu=1$ and the results are
shown in figure \ref{fig:lcd2}. The stabilized scheme is able to
dissipate the energy better than the unstabilized one, but it exhibits
a significant amount of phase error for $P_{1}$ basis functions; the energy dissipation is still seen to be incorrect with $P_1$ basis functions. However, the use of
$P_{2}$ basis functions is seen to give excellent results without
any phase error. In the next sections, we see that adding a penalty term gives correct results even with $P_1$ basis functions.

\subsection{Cell energy analysis}

In the previous sections, we performed a global energy analysis; here
we make a cell energy analysis. If we multiply the convection-diffusion
equation by $u$ and integrate over element $I_{j}$, we obtain
\begin{equation}
\frac{\ud}{\ud t}\|u\|_{I_{j}}^{2}+\left[\frac{1}{2}cu^{2}\right]_{j-\half}^{j+\half}+\left[-\mu uu_{x}\right]_{j-\half}^{j+\half}=-\mu\int_{I_{j}}(u_{x})^{2}dx \le 0
\label{eq:celle}
\end{equation}
We will now try to mimic this property for the DG scheme for which
we will have to define consistent numerical fluxes for the second
and third terms on the left hand side of the above equation. Substituting
$\phi_{h}=u_{h}$ in equation \eqref{eq:lcddg1}, we obtain
\[
\begin{aligned}\frac{\ud}{\ud t}\|u_{h}\|_{I_{j}}^{2} & -\frac{c}{2}\left[u_{h}^{2}(x_{j+\half}^{+})-u_{h}^{2}(x_{j-\half}^{-})\right] +F(u_{h})_{j+\half}u_{h}(x_{j+\half}^{+})-F(u_{h})_{j-\half}u_{h}(x_{j-\half}^{-})\\
 & - F_{d}^{+}(\partial_{x}u_{h}^{+})_{j+\half}\jump{u_{h}}_{j+\half}- F_{d}^{-}(\partial_{x}u_{h}^{-})_{j-\half}\jump{u_{h}}_{j-\half}=-\mu\int_{I_{j}}(\partial_{x}u_{h})^{2}dx
\end{aligned}
\]
The numerical flux $F(u_{h})_{j+\half}$ consists of the convective
and viscous fluxes. Define the average value
\[
\avg{u_{h}}_{j+\half}:=\frac{1}{2}(u_{h}(x_{j+\half}^{+})+u(x_{j+\half}^{-}))
\]
and convective and diffusive numerical fluxes
\begin{eqnarray*}
\bar{F}_{c}(u_{h})_{j+\half} &=& F_{c}(u_{h})_{j+\half}\avg{u_{h}}_{j+\half}-\frac{c}{2}\avg{u_{h}^{2}}_{j+\half} \\
\bar{F}_{d}(u_{h})_{j+\half} &=& F_{d}^{+}(\partial_{x}u_{h}^{+})_{j+\half}u_{h}(x_{j+\half}^{-})+F_{d}^{-}(\partial_{x}u_{h}^{-})_{j+\half}u_{h}(x_{j+\half}^{+})
\end{eqnarray*}
The last two equations provide consistent numerical fluxes for $\frac{1}{2}cu^{2}$ and $-\mu uu_{x}$, respectively. Then we obtain the cell energy equation for the DG scheme as
\[
\frac{\ud}{\ud t}\|u\|_{I_{j}}^{2}+\left[\bar{F}_{c}\right]_{j-\half}^{j+\half}+\left[\bar{F}_{d}\right]_{j-\half}^{j+\half}=-\mu\int_{I_{J}}(\partial_x u_{h})^{2}dx - \frac{1}{4}D\jump{u_{h}}_{\jmh}^{2} - \frac{1}{4}D\jump{u_{h}}_{\jph}^{2}\le 0
\]
which is consistent with the exact energy equation~(\ref{eq:celle}). The inter-element jumps in the solution add additional dissipation in the numerical solution.


\subsection{DG scheme with interior penalty}
The above examples have shown the importance of energy stability in constructing good DG scheme for the convection-diffusion equation. The schemes were written for each element but it is more instructive to write a weak formulation for the whole computational domain since it helps to understand some structure of the equations. By summing up the equation~(\ref{eq:lcddg1}) for each element, we write the DG scheme for convection-diffusion equation as a weak statement: Find $u_h(t) \in V_h^k$ such that for all $\phi_h \in V_h^k$, the following equation is satisfied
\begin{equation}
\begin{aligned}\int_{\Omega}\phi_{h}\partial_{t}u_{h}dx & -\int_{\Omega}cu_{h}\partial_{x}\phi_{h}dx+\mu\int_{\Omega}\partial_{x}u_{h}\partial_{x}\phi_{h}dx + \sum_{j}F_{c}(u_{h})_{j+\half}\jump{\phi_{h}}_{j+\half}\\
 & + \sum_{j}F_{d}(\partial_{x}u_{h})_{j+\half}\jump{\phi_{h}}_{j+\half}+\epsilon\sum_{j}F_{d}(\partial_{x}\phi_{h})_{j+\half}\jump{u_{h}}_{j+\half}
 +\sum_{j}\delta_h(u_{h})_{j+\half}\jump{\phi_h}_{j+\half}=0
\end{aligned}
\label{eq:lcddgip}
\end{equation}
where the last term is an additional term that has been added, called the {\em interior penalty term}, and is of the form~\cite{arnold1982,arnold2002}
\[
\delta_h(u)=\cip\frac{\mu}{h}\jump{u}
\]
where $\cip$ is a non-dimensional positive number. Two classes of schemes are obtained based on whether $\epsilon = -1$ or $+1$. Define the bilinear form consisting of all diffusive terms
\[
a_\epsilon(v_h, w_h) = \sum_{j}F_{d}(\partial_{x}v_{h})_{j+\half}\jump{w_{h}}_{j+\half}+\epsilon\sum_{j}F_{d}(\partial_{x}w_{h})_{j+\half}\jump{v_{h}}_{j+\half}
 +\sum_{j}\delta_h(v_{h})_{j+\half}\jump{w_h}_{j+\half}
 \]
\begin{enumerate}
\item If we choose $\epsilon=-1$ in (\ref{eq:lcddgip}) and omit the interior penalty term,
then we obtain the scheme of equation \eqref{eq:lcddg1}. The scheme is non-symmetric since $a_{-1}(v_h,w_h) \ne a_{-1}(w_h,v_h)$. As we have seen, this scheme satisfies a cell energy inequality. The interior penalty terms are necessary to make the scheme stable since
\[
a_{-1}(v_h,v_h) = \mu \int_\Omega (\partial_x v_h)^2 \ud x + \sum_{j=1}^N \cip\frac{\mu}{h}\jump{v_h}^2
\]
and the right handside is a norm on the space of discontinuous functions only in the presence of the penalty term for any $\cip>0$. Without the penalty term ($\cip=0$), the right hand side is only a semi-norm and  the scheme is only weakly stable~\cite{arnold2002}. In this case, the discretization of diffusive terms is similar to that of Baumann and Oden~\cite{Baumann1999311} for elliptic equation, who also show sub-optimal convergence in the case $k \ge 2$ only. There is no proof of convergence for $k=1$ case which combined with only weak stability explains the poor performance of the scheme shown in figure~(\ref{fig:lcd2}a). With the interior penalty terms, the global energy equation is
\begin{equation*}
\frac{1}{2}\frac{\ud}{\ud t}\|u_{h}\|^{2}+\frac{1}{2}D\sum_{j=1}^{N} \jump{u_{h}}_{j+\half}^{2} + \mu\sum_{j=1}^{N}\int_{I_{j}}(\partial_{x}u_{h})^{2}dx + \sum_{j=1}^N \cip\frac{\mu}{h}\jump{u_h}^2 =0
\end{equation*}
which shows that interior penalty adds additional stabilization to the scheme which is important when viscosity is dominant. The scheme with interior penalty is known as NIPG and is stable for any $k \ge 1$.
\item If
$\epsilon=+1$, the discretization of the diffusive terms becomes
symmetric, i.e., $a_1(v_h,w_h)=a_1(w_h,v_h)$, and the scheme is known as SIPG scheme. However, the SIPG
scheme does not satisy a cell energy inequality, though in the case of elliptic equation it can be
shown to be globally stable~\cite{arnold2002}. This scheme is stable only with the addition of a penalty term $\delta_h$, with a sufficiently large value of $\cip$. The determination of a good value of $\cip$ is an open problem; in practice a value of $\cip$ between 5 to 10 seems to yield satisfactory results. A larger value of $\cip$ reduces the allowable time step for explicit time integration schemes, requiring the necessity of implicit time integration schemes.
\end{enumerate}
The jump terms in the above DG scheme do not destroy the consistency of the scheme; if we substitute the exact solution of the convection-diffusion equation which is continuous in equation~(\ref{eq:lcddgip}), it exactly satisfies the equation since the jump terms vanish.
\subsection{Numerical results using interior penalty scheme}

\begin{figure}
\begin{center}
\begin{tabular}{cc}
\includegraphics[width=0.30\textwidth]{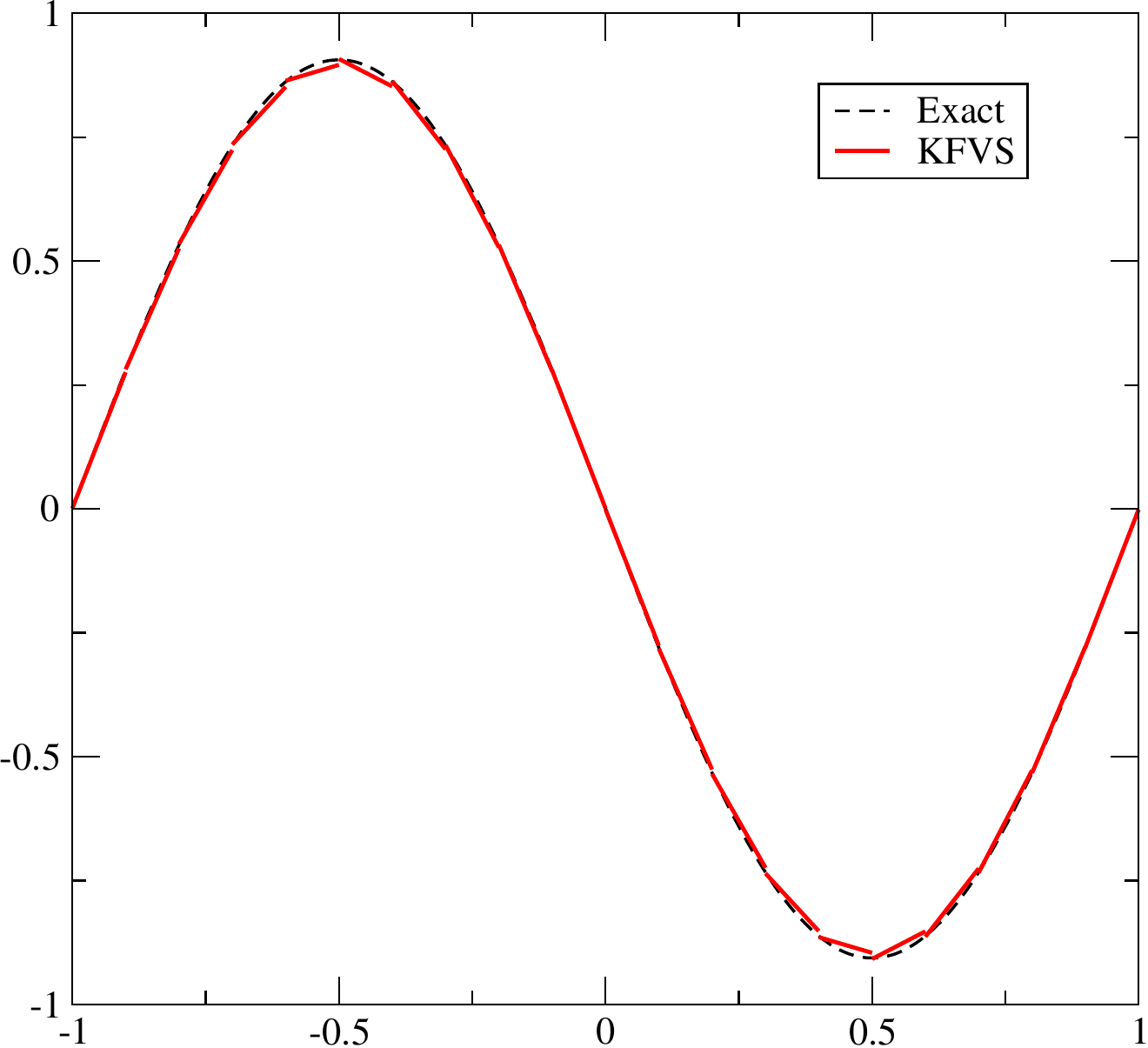} &
\includegraphics[width=0.30\textwidth]{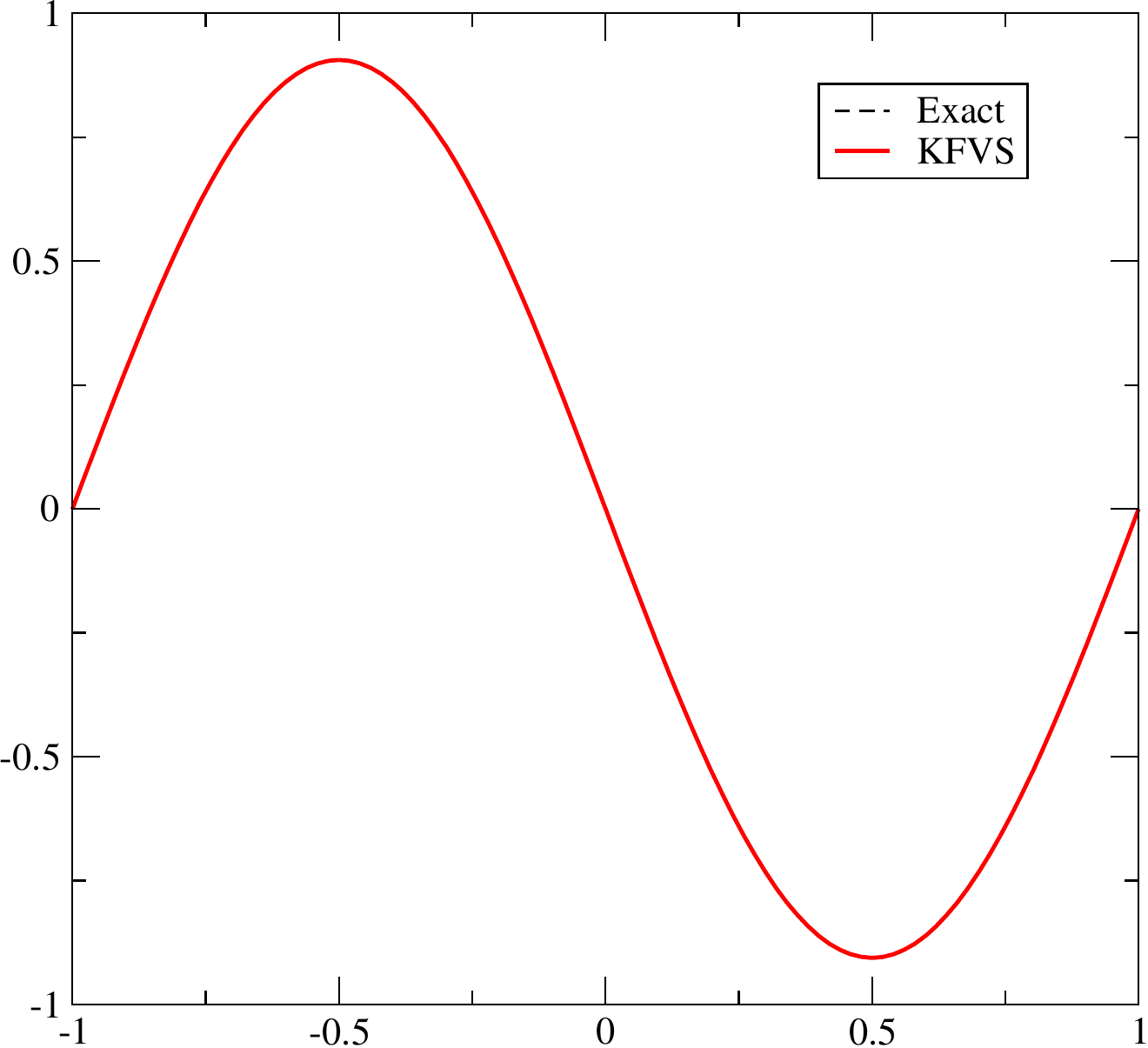} \\
NIPG, $P_1$  &  NIPG, $P_2$ \\
\includegraphics[width=0.30\textwidth]{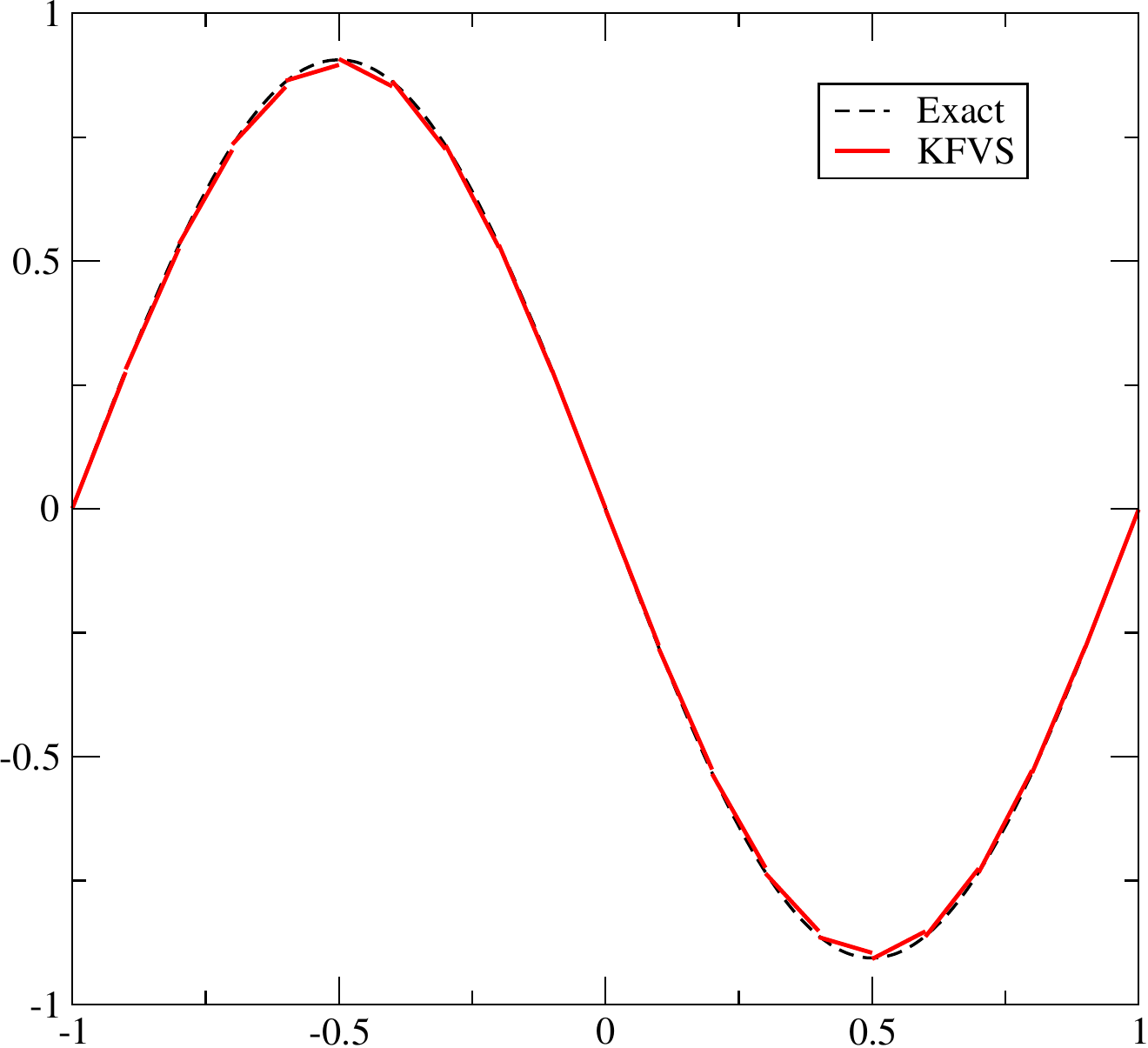} &
\includegraphics[width=0.30\textwidth]{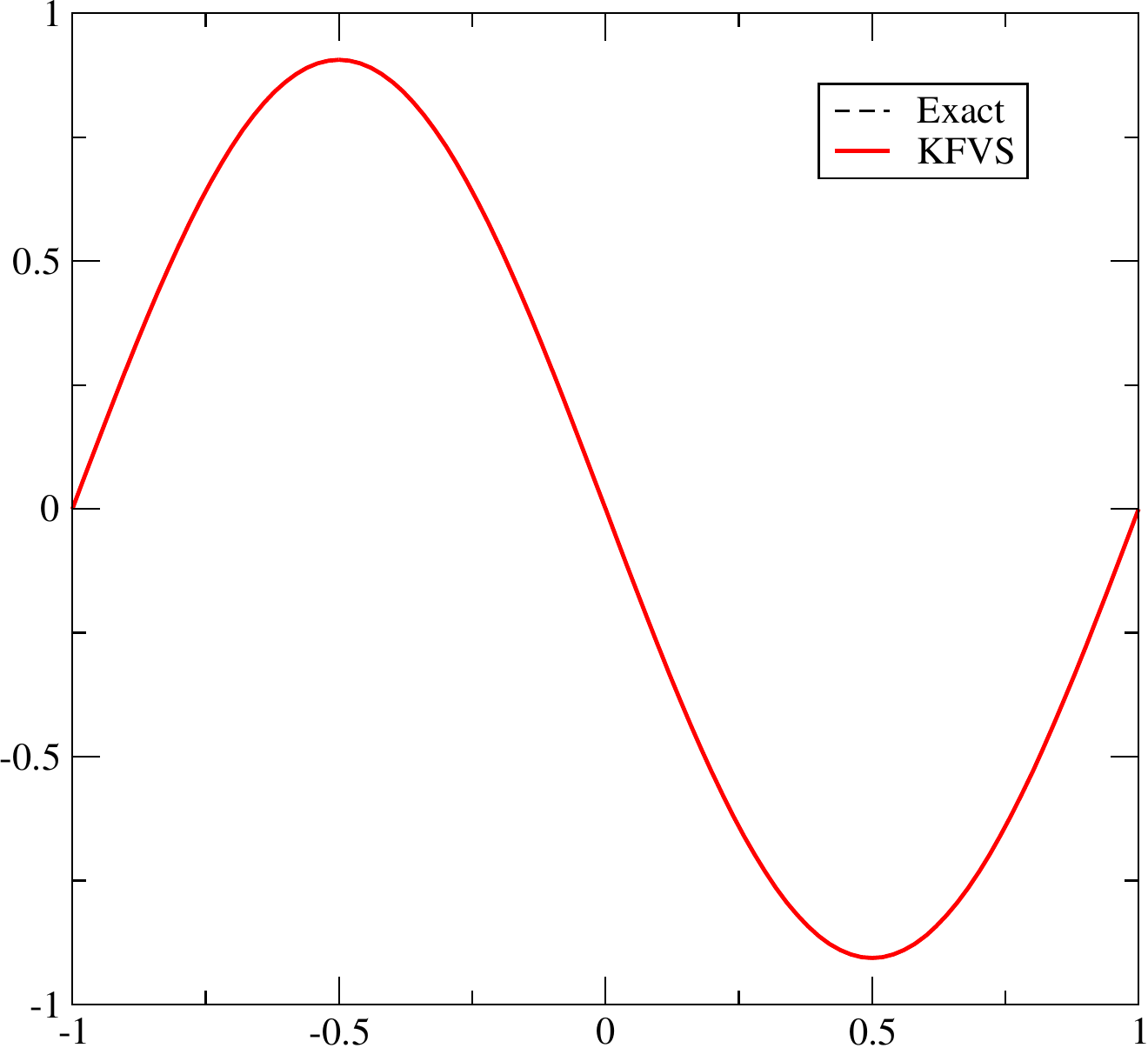}\\
SIPG, $P_1$  &  SIPG, $P_2$ \\
\end{tabular}
\end{center}
\caption{Results for Test case 1 using interior penalty scheme}
\label{fig:t1}
\end{figure}

\begin{figure}
\begin{center}
\begin{tabular}{cc}
\includegraphics[width=0.30\textwidth]{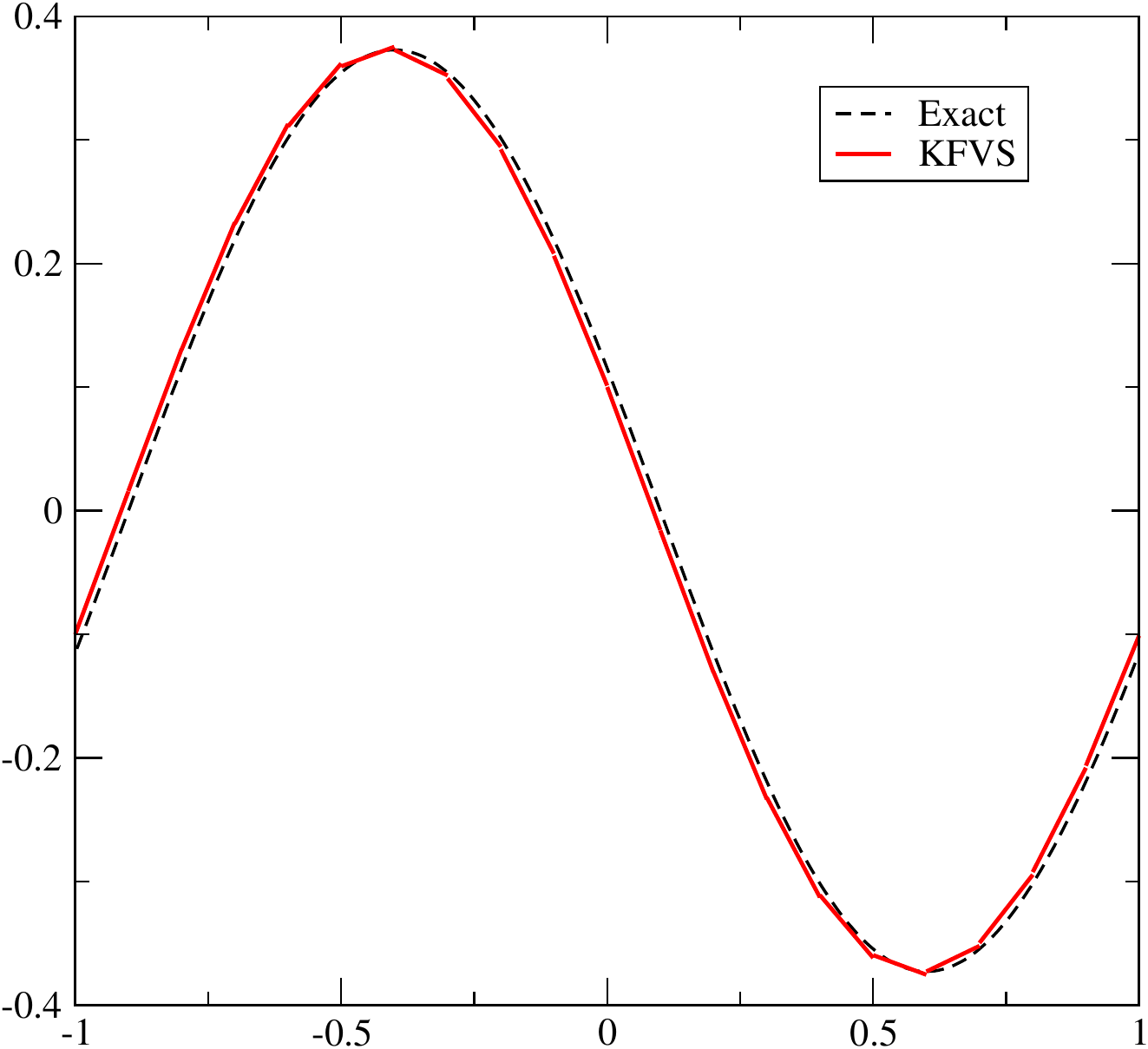} &
\includegraphics[width=0.30\textwidth]{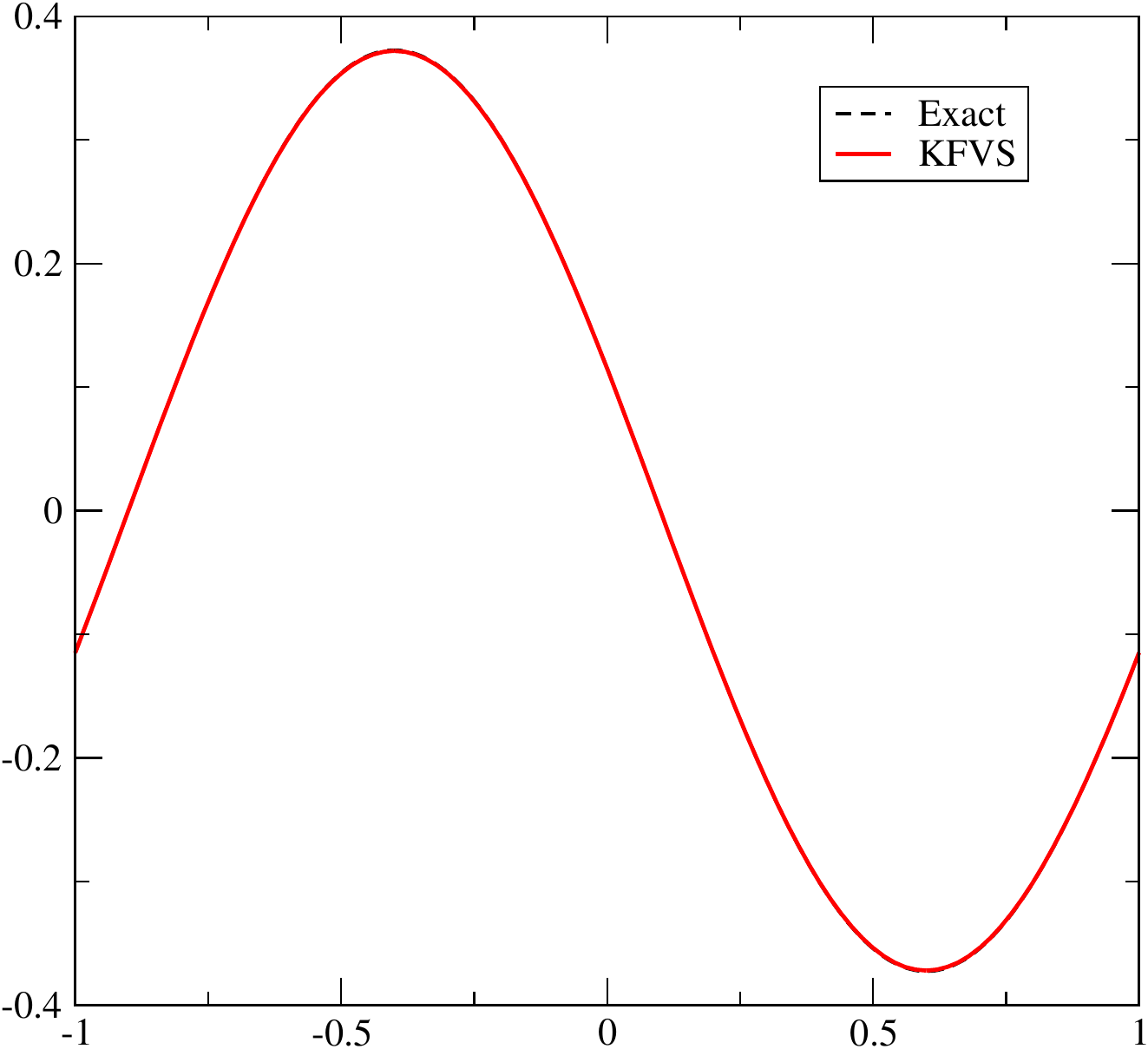} \\
NIPG, $P_1$  &  NIPG, $P_2$ \\
\includegraphics[width=0.30\textwidth]{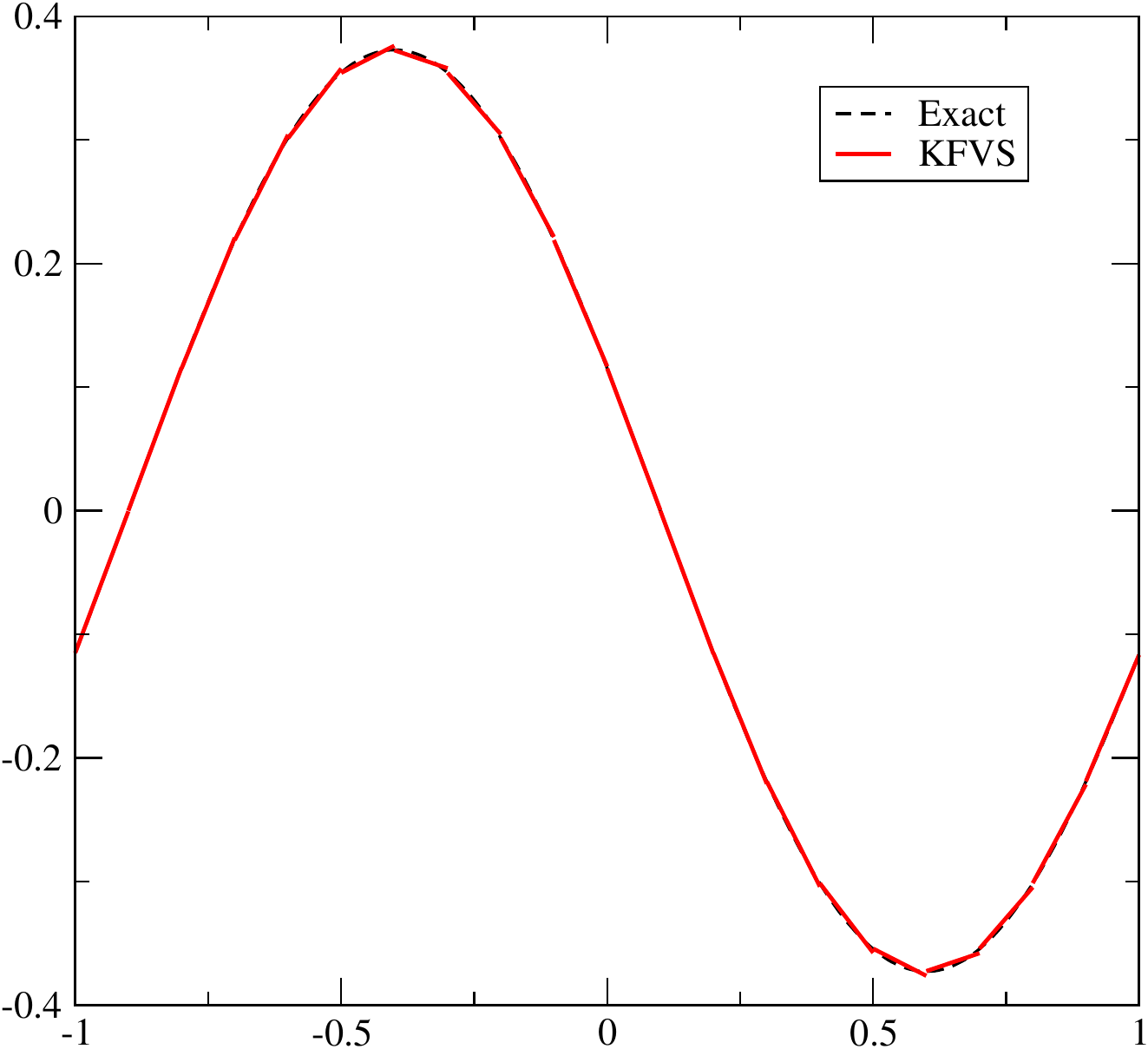} &
\includegraphics[width=0.30\textwidth]{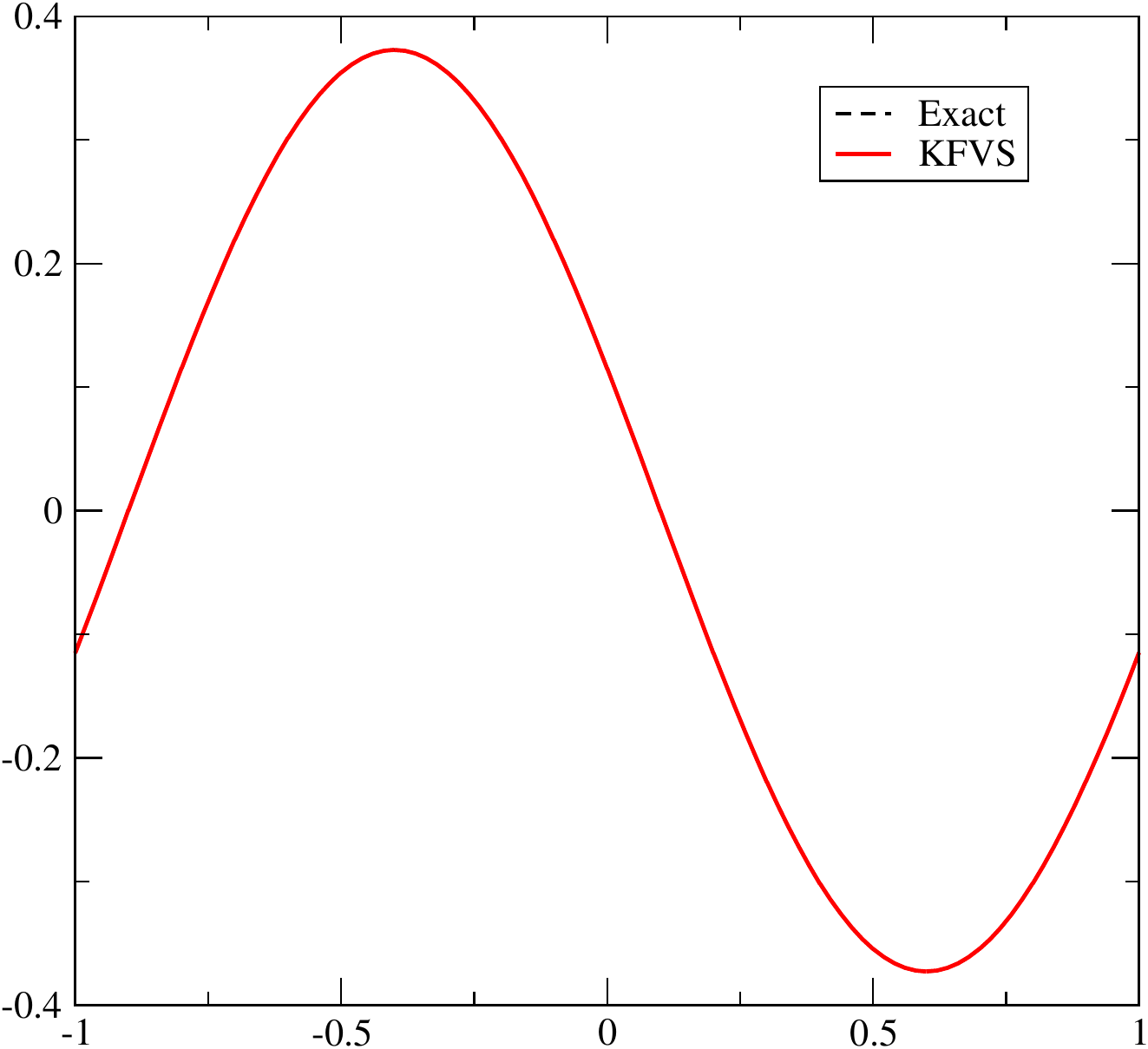}\\
SIPG, $P_1$  &  SIPG, $P_2$ \\
\end{tabular}
\end{center}
\caption{Results for Test case 2 using interior penalty scheme}
\label{fig:t2}
\end{figure}

In this section, we present results for the two test cases, using the interior penalty scheme given by equation~(\ref{eq:lcddgip}). We first present results using 20 cells. The schemes NIPG and SIPG are tested with a penalty parameter $\cip=10$. Figure~(\ref{fig:t1}) shows the results for Test case 1 which is convection dominated situation. Figure~(\ref{fig:t2}) shows results for Test case 2 which is a diffusion dominated problem. We see that the addition of the interior penalty leads to better solutions with $P_1$ basis functions, as compared to figure~(\ref{fig:lcd2}).

In order to compute the convergence rates of the various schemes, we perform a grid refinement study. Errors in the solution and its derivative are computed in $L^2$ norm and shown in tables~(\ref{tab:t1}-\ref{tab:t2}) for $P_1$ and $P_2$ basis functions. We observe the following trends:
\begin{itemize}
\item SIPG gives more accurate results than NIPG. This can be seen by comparing the error values for the two schemes given in the tables.
\item NIPG converges at sub-optimal rates; both the solution and its derivative converge at $O(h^k)$. 
\item SIPG converges at optimal rates; the solution converges at $O(h^{k+1})$ and the derivative converges at $O(h^k)$.
\end{itemize}
Actually both the schemes seem to exhibit slightly higher convergence rates on Test case 1 which is convection dominated. In Test case 2 which is viscosity dominated, the convergence rates are very close to the above mentioned values. These convergence rates are consistent with what is observed in the case of DG schemes applied to elliptic equation~\cite{arnold2002}. The NIPG scheme for elliptic problems is numerically found to show optimal convergence rate for odd degree polynomials but we have not observed this phenomenon for the convection-diffusion equation with the current NIPG scheme.

\begin{table}
\begin{center}
\begin{tabular}{|r|r|c|c|c|c|} \hline
 cells &  dofs & 
\multicolumn{2}{|c|}{L2} & 
\multicolumn{2}{|c|}{H1}\\ \hline
20 & 40 & 1.259e-02 & - & 2.757e-01 & -\\ \hline
40 & 80 & 3.535e-03 & 1.83 & 1.356e-01 & 1.02\\ \hline
80 & 160 & 1.245e-03 & 1.51 & 6.639e-02 & 1.03\\ \hline
160 & 320 & 3.511e-04 & 1.83 & 3.257e-02 & 1.03\\ \hline
320 & 640 & 3.986e-05 & 3.14 & 1.614e-02 & 1.01\\ \hline
\end{tabular}
\caption{Test case 1: NIPG, $P_1$ basis}
\label{tab:t1}
\end{center}
\end{table}

\begin{table}
\begin{center}
\begin{tabular}{|r|r|c|c|c|c|} \hline
 cells &  dofs & 
\multicolumn{2}{|c|}{L2} & 
\multicolumn{2}{|c|}{H1}\\ \hline
20 & 60 & 3.376e-04 & - & 1.416e-02 & -\\ \hline
40 & 120 & 7.361e-05 & 2.20 & 3.612e-03 & 1.97\\ \hline
80 & 240 & 1.660e-05 & 2.15 & 9.186e-04 & 1.98\\ \hline
160 & 480 & 3.433e-06 & 2.27 & 2.243e-04 & 2.03\\ \hline
320 & 960 & 6.892e-07 & 2.32 & 5.103e-05 & 2.14\\ \hline
\end{tabular}
\caption{Test case 1: NIPG, $P_2$ basis}
\end{center}
\end{table}

\begin{table}
\begin{center}
\begin{tabular}{|r|r|c|c|c|c|} \hline
cells & dofs & 
\multicolumn{2}{|c|}{L2} & 
\multicolumn{2}{|c|}{H1}\\ \hline
20 & 40 & 1.245e-02 & - & 2.866e-01 & -\\ \hline
40 & 80 & 1.843e-03 & 2.76 & 1.403e-01 & 1.03\\ \hline
80 & 160 & 3.170e-04 & 2.54 & 6.807e-02 & 1.04\\ \hline
160 & 320 & 6.181e-05 & 2.36 & 3.293e-02 & 1.05\\ \hline
320 & 640 & 1.351e-05 & 2.19 & 1.615e-02 & 1.03\\ \hline
\end{tabular}
\caption{Test case 1: SIPG, $P_1$ basis}
\end{center}
\end{table}

\begin{table}
\begin{center}
\begin{tabular}{|r|r|c|c|c|c|} \hline
 cells &  dofs & 
\multicolumn{2}{|c|}{L2} & 
\multicolumn{2}{|c|}{H1}\\ \hline
20 & 60 & 1.354e-04 & - & 1.367e-02 & -\\ \hline
40 & 120 & 1.622e-05 & 3.06 & 3.296e-03 & 2.05\\ \hline
80 & 240 & 1.900e-06 & 3.09 & 7.692e-04 & 2.10\\ \hline
160 & 480 & 2.153e-07 & 3.14 & 1.748e-04 & 2.14\\ \hline
320 & 960 & 2.487e-08 & 3.11 & 4.324e-05 & 2.02\\ \hline
\end{tabular}
\caption{Test case 1: SIPG, $P_2$ basis}
\end{center}
\end{table}

\begin{table}
\begin{center}
\begin{tabular}{|r|r|c|c|c|c|} \hline
 cells &  dofs & 
\multicolumn{2}{|c|}{L2} & 
\multicolumn{2}{|c|}{H1}\\ \hline
20 & 40 & 1.592e-02 & - & 1.088e-01 & -\\ \hline
40 & 80 & 8.018e-03 & 0.99 & 5.451e-02 & 1.00\\ \hline
80 & 160 & 4.031e-03 & 0.99 & 2.727e-02 & 1.00\\ \hline
160 & 320 & 2.022e-03 & 1.00 & 1.364e-02 & 1.00\\ \hline
320 & 640 & 1.013e-03 & 1.00 & 6.820e-03 & 1.00\\ \hline
\end{tabular}
\caption{Test case 2: NIPG, $P_1$ basis}
\end{center}
\end{table}

\begin{table}
\begin{center}
\begin{tabular}{|r|r|c|c|c|c|} \hline
cells & dofs & 
\multicolumn{2}{|c|}{L2} & 
\multicolumn{2}{|c|}{H1}\\ \hline
20 & 60 & 7.467e-04 & - & 4.778e-03 & -\\ \hline
40 & 120 & 1.860e-04 & 2.01 & 1.195e-03 & 2.00\\ \hline
80 & 240 & 4.648e-05 & 2.00 & 2.989e-04 & 2.00\\ \hline
160 & 480 & 1.162e-05 & 2.00 & 7.472e-05 & 2.00\\ \hline
320 & 960 & 2.906e-06 & 2.00 & 1.868e-05 & 2.00\\ \hline
\end{tabular}
\caption{Test case 2: NIPG, $P_2$ basis}
\end{center}
\end{table}

\begin{table}
\begin{center}
\begin{tabular}{|r|r|c|c|c|c|} \hline
 cells &  dofs & 
\multicolumn{2}{|c|}{L2} & 
\multicolumn{2}{|c|}{H1}\\ \hline
 20 & 40 & 2.255e-03 & - & 1.124e-01 & -\\ \hline
 40 & 80 & 5.653e-04 & 2.00 & 5.638e-02 & 1.00\\ \hline
 80 & 160 & 1.415e-04 & 2.00 & 2.823e-02 & 1.00\\ \hline
 160 & 320 & 3.538e-05 & 2.00 & 1.412e-02 & 1.00\\ \hline
 320 & 640 & 8.846e-06 & 2.00 & 7.063e-03 & 1.00\\ \hline
\end{tabular}
\caption{Test case 2: SIPG, $P_1$ basis}
\end{center}
\end{table}

\begin{table}
\begin{center}
\begin{tabular}{|r|r|c|c|c|c|} \hline
 cells &  dofs & 
\multicolumn{2}{|c|}{L2} & 
\multicolumn{2}{|c|}{H1}\\ \hline
 20 & 60 & 1.355e-04 & - & 9.936e-03 & -\\ \hline
 40 & 120 & 1.713e-05 & 2.98 & 2.497e-03 & 1.99\\ \hline
 80 & 240 & 2.158e-06 & 2.99 & 6.268e-04 & 1.99\\ \hline
 160 & 480 & 2.709e-07 & 2.99 & 1.571e-04 & 2.00\\ \hline
 320 & 960 & 3.394e-08 & 3.00 & 3.932e-05 & 2.00\\ \hline
\end{tabular}
\caption{Test case 2: SIPG, $P_2$ basis}
\label{tab:t2}
\end{center}
\end{table}

\subsection{Test case 3}

\begin{figure}
\begin{center}
\begin{tabular}{cc}
\includegraphics[width=0.48\textwidth]{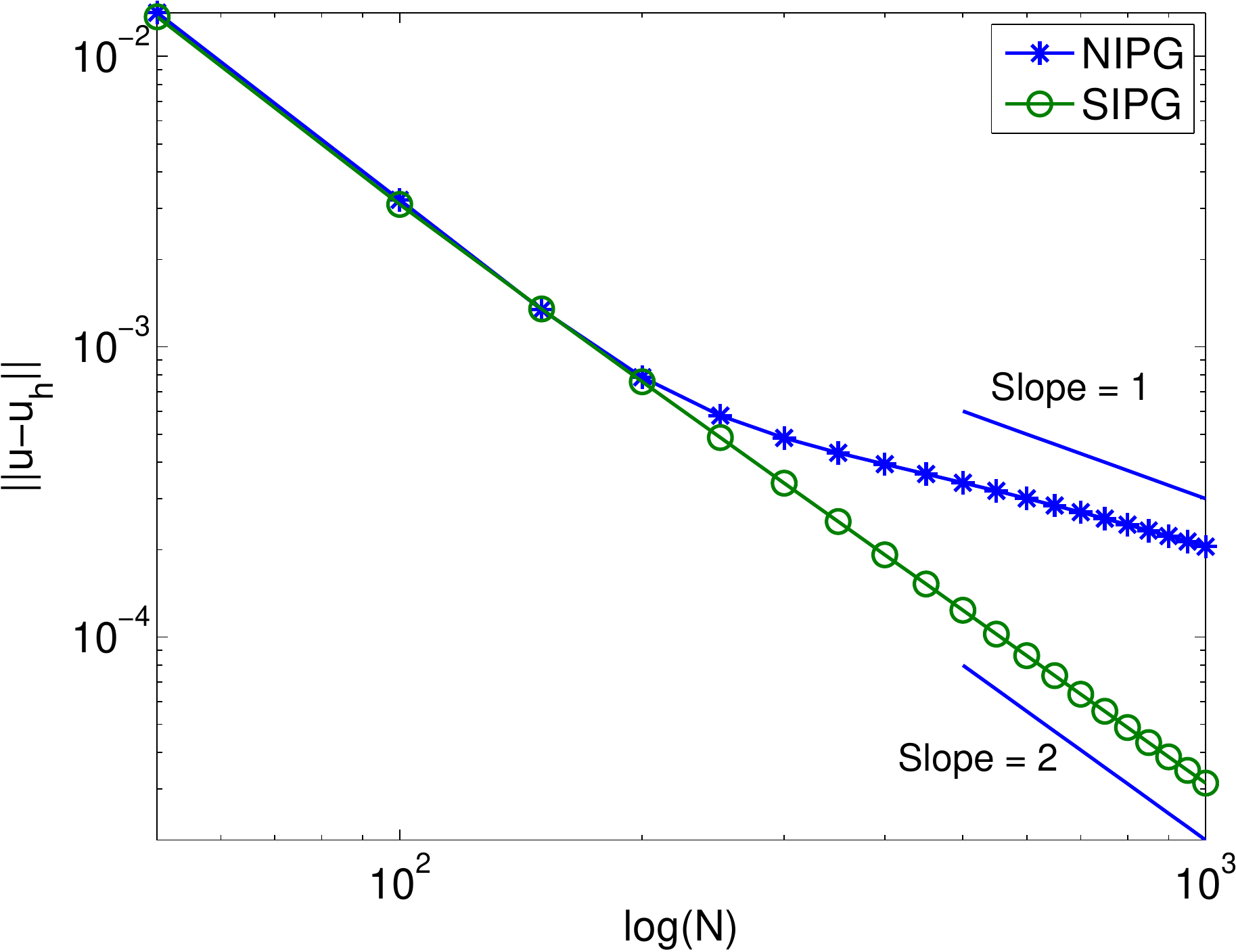} &
\includegraphics[width=0.48\textwidth]{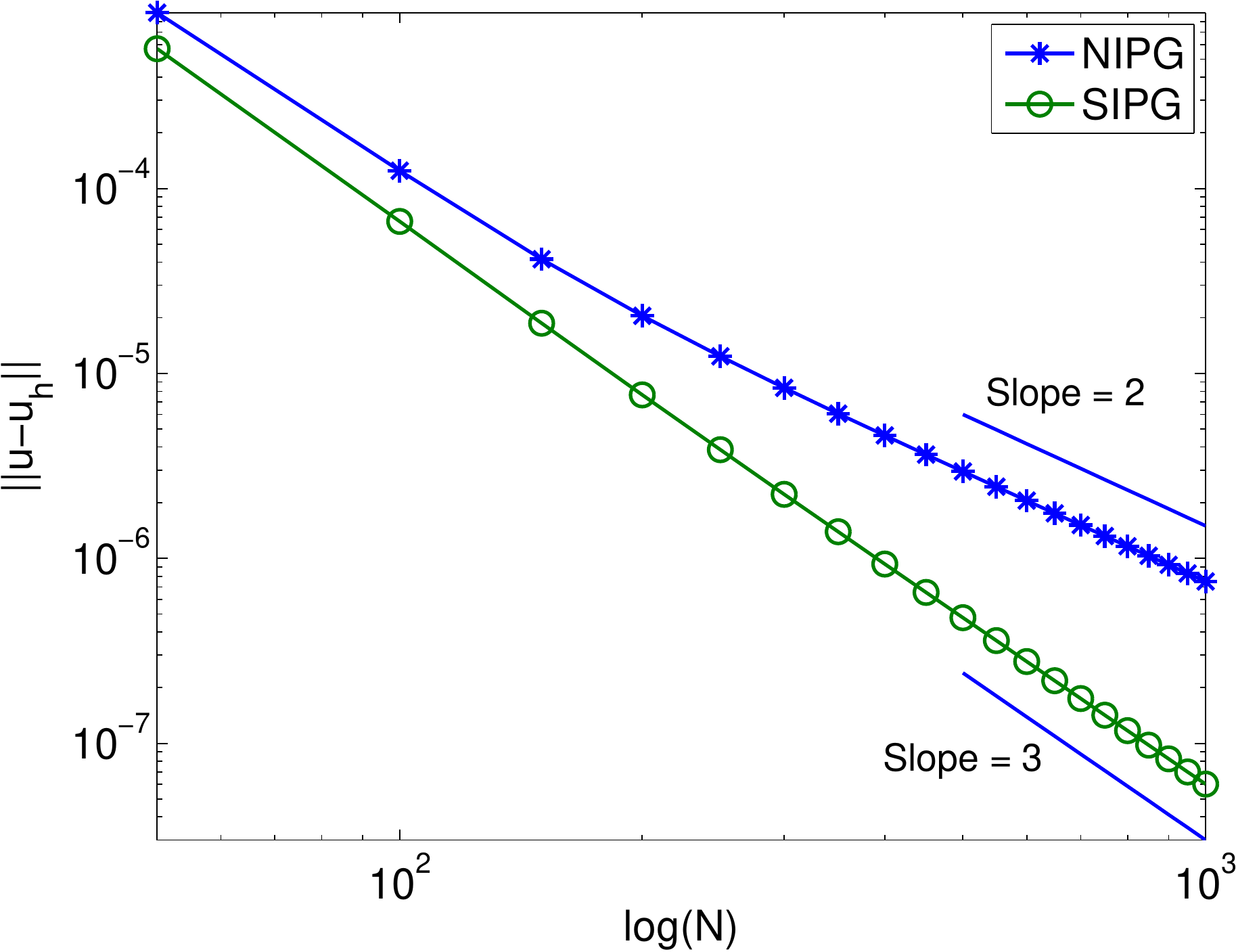} \\
(a) & (b)
\end{tabular}
\caption{Test case 3; convergence of error for (a) $P_1$ and (b) $P_2$ basis functions}
\label{fig:t3}
\end{center}
\end{figure}

Torrilhon and Xu~\cite{torr_xu} have analyzed the BGK finite volume scheme for scalar convection-diffusion equation and find that the BGK scheme exhibits only first order accuracy in the asymptotic limit, though on coarse grids it shows third order accuracy. They find a non-monotone convergence of the error with grid size. In order to study the behaviour of the error convergence, we apply the present DG scheme on a test case similar to the one given in \cite{torr_xu} with $\beta=1$, $\mu=0.005$, $c=1$, using $P_1$ basis functions. The grid is refined starting with 50 cells and increased by 50 cells in each refinement step. The initial condition is taken to be
\begin{equation*}
u(x,0) = 4 + \frac{8}{\pi} \sin(\pi x/2) + \frac{16}{3\pi} \sin(3\pi x/2)
\end{equation*}
for which an exact solution is available~\cite{torr_xu} for the case of periodic boundary conditions using Fourier approach. The convergence of the $L_2$ error in the solution as a function of the number of cells is shown in figure~(\ref{fig:t3}a) for $P_1$ basis functions; this plot must be compared with figure (7) in \cite{torr_xu}. It is clear that both the NIPG and SIPG schemes converge monotonically with grid refinement. Asymptotically, the NIPG scheme converges at first order while the SIPG scheme converges at second order. Surprisingly, the NIPG scheme exhibits second order convergence on coarse grids, a superconvergence phenomenon that is also observed in the BGK finite volume scheme~\cite{torr_xu}. The convergence of the error for $P_2$ basis functions is shown in figure~(\ref{fig:t3}b) which shows second and third order convergence for NIPG and SIPG schemes, respectively. On coarser grids, the NIPG scheme again shows optimal convergence rate. These convergence rates are consistent with what is observed in the previous two test cases, and also in the 1-D Navier-Stokes equations given in the later sections.
\section{Navier-Stokes equations}
\label{sec:ns}

The Navier-Stokes equations which represent conservation laws for mass, momentum and energy for a fluid can be written as a system of partial differential equations
\[
\partial_t U + \partial_{i}F_{i}(U)+\partial_{i}G_{i}(U,\nabla U)=0
\]
where $U$ is the vector of {\em conserved variables} corresponding to density
of mass, momentum, and energy, while $F_{i}$, $i=1,2,3$ are the
inviscid or convective fluxes and $G_{i}$ are the viscous fluxes.
The conserved variables are given by
\[
U=\left[\begin{array}{ccccc}
\rho & \rho u_{1} & \rho u_{2} & \rho u_{3} & \rho e\end{array}\right]^{\top}
\]
where $\rho$ is the density, $u=(u_{1},u_{2},u_{3})$ is the fluid
velocity and $e$ is the energy per unit mass, which can be written
in terms of the internal energy $\varepsilon$ and the kinetic energy
as $e=\varepsilon+\frac{1}{2}u^{2}$. An equation of state $\varepsilon=\varepsilon(\rho,p)$ is required
to close the system of equations, where $p$ is the pressure; for a polytropic ideal gas, the equation of state is given by $\varepsilon=\frac{p}{\rho(\gamma-1)}$ where $\gamma=\frac{C_p}{C_v}$ is the ratio of specific heats at constant pressure and constant volume~\cite{liepmann-roshko}. The flux vectors are given by
\begin{equation}
F_i = \begin{bmatrix}
\rho u_i \\
p \delta_{i1} + \rho u_i u_1 \\
p \delta_{i2} + \rho u_i u_2 \\
p \delta_{i3} + \rho u_i u_3 \\
(\rho e + p) u_i
\end{bmatrix}, \qquad
G_i = \begin{bmatrix}
0 \\
-\tau_{i1} \\
-\tau_{i2} \\
-\tau_{i3} \\
-\tau_{ij} u_j + q_i
\end{bmatrix}
\end{equation}
The shear stress tensor $\tau$ is related to the strain rate by Newton's constitutive law while the heat flux is related to the gradient of absolute temperature $T$ by Fourier's law
\begin{equation}
\tau_{ij} = \mu (\partial_i u_j + \partial_j u_i) - \frac{2}{3} \mu (\partial_k u_k) \delta_{ij}, \qquad q_i = -\kappa \partial_i T
\end{equation}
The coefficient of dynamic viscosity $\mu$ and the coefficient of heat conduction $\kappa$ are related through the Prandtl number $\Pr = \frac{\mu C_p}{\kappa}$. 
We will
use the notation $F=[F_{1},F_{2},F_{3}]$ and $G=[G_{1},G_{2,}G_{3}].$
Moreover, for any vector $n\in\mathbb{R}^{3}$, we denote $F\cdot n=F_{i}n_{i}$ and $G\cdot n=G_{i}n_{i}$ for the fluxes in the direction $n$. Since the numerical experiments in this paper are conducted for the 1-D Navier-Stokes equations, the corresponding equations are given in \ref{sec:ns1d}.


\subsection{Euler equations, entropy variables and symmetric hyperbolic form}
The system of equations containing only the convective terms
\begin{equation}
\partial_t U+\partial_{i}F_{i}=0\label{eq:euler}
\end{equation}
are called the Euler equations, and are hyperbolic in nature. This
means that for any vector $n=(n_{1},n_{2},n_{3})$, the matrices
\[
A(U,n)=A_{i}(U)n_{i},\qquad A_{i}(U)=F_{i}'(U)
\]
have all real eigenvalues and a complete set of eigenvectors. The
Euler equations admit an entropy function $\eta(U)$ with entropy
fluxes $\theta_{i}(U),\: i=1,2,3$ such that%
\footnote{For a scalar function of a vector like $\eta(U)$ we consider $\eta'(U)$
to be a row vector.%
}
\[
\theta_{i}'(U)=\eta'(U)F_{i}'(U)
\]
This implies that for smooth solutions, an additional conservation
law is satisfied, i.e., $\partial_t \eta+\partial_{i}\theta_{i}=0$, while for discontinuous solutions the inequality $\partial_t \eta+\partial_{i}\theta_{i}\le0$ is satisfied in the weak sense. The entropy and entropy fluxes can be 
related to the physical entropy $s=\ln(p/\rho^\gamma) + s_0$ by
\begin{equation}
\eta=-\frac{\rho s}{\gamma-1},\qquad\theta_{i}=-\frac{\rho su_{i}}{\gamma-1}
\label{eq:eulent}
\end{equation}
Actually any function of the form $\eta=\rho h(s)$ where $h$ is convex is an entropy for the Euler equations, but the above linear choice is the only one which is consistent with the second law of thermodynamics in the presence of heat conduction~\cite{Shakib1991141}. Due to the mass conservation equation, the constant term $s_0$ can be ignored and moreover it can be shown that $\eta(U)$ is a strictly convex function.
We can then define the entropy variables by
\begin{equation}
V^{\top}=\eta'(U)\label{eq:entvar}
\end{equation}
which can be inverted in a unique manner to obtain $U=U(V)$. Let us define the dual variables
by a Legendre transform
\[
\eta^{*}(V)=V\cdot U(V)-\eta(U(V)),\qquad\theta_{i}^{*}(V)=V\cdot F_{i}(U(V))-\theta_{i}(U(V))
\]
Differentiating these expressions it is easy to check that $\eta^{*'}(V)=U^{\top}$
and $\theta_{i}^{*'}(V)=F_{i}(U(V))^{\top}$.  It then follows that the matrices $U'(V)=\eta^{*''}(V)$ and $F_{i}'(U(V))U'(V)=\theta_{i}^{*''}(V)$
are symmetric. Moreover the matrix
\[
U'(V)=\eta''(U(V))^{-1}
\]
is positive-definite. Making a change of variables from $U$ to $V$ in equation \eqref{eq:euler},
we get
\[
\tilde{A}_{0}\partial_t V+\tilde{A}_{i}\partial_{i}V=0
\]
where we have defined
\[
\tilde{A}_{0}=U'(V)=\eta^{*''}(V),\qquad\tilde{A}_{i}=A_{i}\tilde{A}_{0}=\theta_{i}^{*''}(V)
\]
the matrix $\tilde{A}_{0}$ is symmetric positive definite, while
the matrices $\tilde{A}_{i}$ are symmetric.


\subsection{Kinetic description of Euler equations}

The Boltzmann equation~\cite{cercignani1988boltzmann} gives a molecular statistical description of
gas dynamics in the form of the velocity distribution function $f(v,x,t)$
which can be written as
\[
\partial_t f+v_{i}\partial_{i}f=Q(f,f)
\]
where the right hand side represents the effect of intermolecular
collisions. The meaning of the distribution function is that the quantity
$f(v,x,t)\textrm{d}x\textrm{d}v$ gives the average number of particles
found in the spatial volume $\textrm{d}x$ around the point $x$ at time $t$ and
having velocities in the range $\textrm{d}v$ around the velocity
$v$. Thus $f$ gives the number density of molecules in phase space
and by a suitable scaling, it can be considered as providing the mass
density. The macrosopic fluid motion is the average effect of all
molecular motions, so that the conserved variables are obtained by
integrating the distribution function over the velocity space
\[
U=\int_{\mathbb{R}^{3}}\psi(v)f(v,x,t)\textrm{d}v=\langle\psi f\rangle
\]
where
\[
\psi(v)=\left[\begin{array}{ccccc}
1 & v_{1} & v_{2} & v_{3} & \frac{1}{2}|v|^{2}\end{array}\right]^{\top}
\]
are the microscopic collisional invariants. Since collisions do not
destroy mass, momentum and energy, we have $\langle\psi Q\rangle=0$.
In fact, it can be shown~\cite{cercignani1988boltzmann} that any function $\phi(v)$ of the molecular velocities for which $\langle\phi Q\rangle=0$
must be of the form $\phi(v)=\sum a_{i}\psi_{i}$.

A gas in thermodynamic equilibrium is characterized by the fact that
there are no spatial or temporal variations in the distribution function,
i.e., $Q(f,f)=0$. It can be shown that the unique solution $g=f$
of this equation is of the form
\[
g=\exp(V\cdot\psi)
\]
for some constants $V\in\mathbb{R}^{5}$. This is the usual Maxwell-Boltzmann
distribution function and can be written in the more familiar form
\[
g(v)=\rho\left(\frac{\beta}{\pi}\right)^{3/2}\exp\left[-\beta|v-u|^{2}\right],\qquad\beta=\frac{1}{2RT}
\]
The Maxwell distribution can also be characterized as the solution
of the following minimization problem~\cite{perthame1990}
\[
\min\left\{ H[f]\::\:\langle\psi f\rangle=U\right\} ,\qquad H[f]=\langle f\ln f\rangle
\]
and the quantities $V$ appear as the Lagrange multipliers to enforce
the constraint. If we define
\[
\eta^{*}(V)=\langle f \rangle = \langle\exp(V\cdot\psi)\rangle
\]
then
\begin{equation}
\eta^{*'}(V)=\langle\psi^{\top}\exp(V\cdot\psi)\rangle=U^{\top}\label{eq:etaskin}
\end{equation}
and hence the Lagrange multipliers $V$ are precisely the entropy variables corresponding to the macroscopic state $U$. The modifications in the case of polyatomic gases are briefly discussed in \ref{sec:poly}.


\subsection{From Maxwell distribution to Euler equations}

The Maxwell distribution was obtained under the condition of thermodynamic
equilibrium corresponding to a gas which has uniform properties in
space and time. The appropriate moments of the distribution function
yield the conserved variables and the inviscid fluxes, i.e., $U=\mom{\psi(v)g}$ and $F_i=\mom{v_i \psi(v) g}$. Due to this
correspondence, it is common to use the collision-less Boltzmann equation
\[
\partial_t g+v_{i}\partial_{i}g=0
\]
to construct numerical schemes for the Euler equations. When the gas
is not uniform, then the assumption is that locally, it is still in
thermodynamic equilibrium, so that one has a local Maxwell distribution
\[
g(v,x,t)=\exp(V(x,t)\cdot\psi)
\]
defined at every point of space. This distribution however does not
satisfy the collision-less Boltzmann equation. We can derive an equation
for $g$ by differentiating the above equation
\[
\partial_t g+v_{i}\partial_{i}g=g\left[\partial_t V+v_{i}\partial_{i}V\right]\cdot\psi=g\left[-\tilde{A}_{0}^{-1}\tilde{A}_{i}+v_{i}I\right]\partial_{i}V\cdot\psi
\]
Taking moments on both sides, we obtain
\begin{equation}
\langle\psi(\partial_t g+v_{i}\partial_{i}g)\rangle=-\langle\psi\otimes\psi g\rangle\tilde{A}_{0}^{-1}\tilde{A}_{i}\partial_{i}V+\langle v_{i}\psi\otimes\psi g\rangle\partial_{i}V\label{eq:eul1}
\end{equation}
Now from equation \eqref{eq:etaskin}
\[
\eta^{*''}(V)=\langle\psi\otimes\psi g\rangle=U'(V)=\tilde{A}_{0}
\]
Moreover, since $F_{i}=\langle v_{i}\psi g\rangle$, we get $F_{i}'(V)=\langle v_{i}\psi\otimes\psi g\rangle=\tilde{A}_{i}$, which shows that the right hand side of equation \eqref{eq:eul1} vanishes, and we obtain the Euler equations. Due to this reason the right hand side is ignored and we use the collision-less Boltzmann equation with Maxwell-Boltzmann distribution function to derive numerical schemes for Euler equations.

\subsection{DG scheme for Euler equations}

\begin{figure}
\begin{center}
\includegraphics[width=0.35\textwidth]{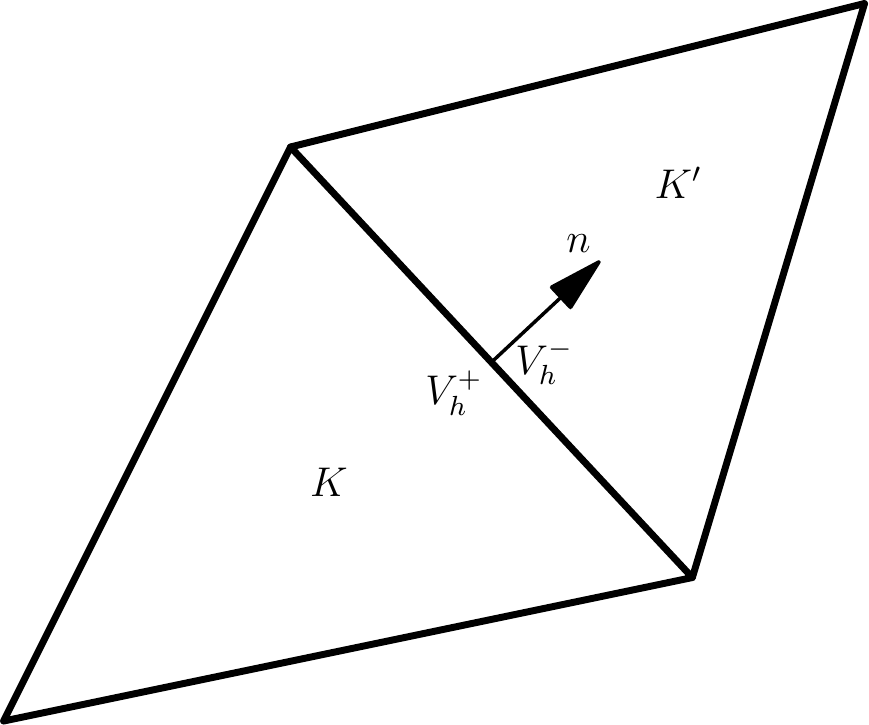}
\caption{Definition of edge normal and trace values}
\label{fig:edgenormal}
\end{center}
\end{figure}
Consider a triangulation $\Th=\{K\}$ of the domain $\Omega$ by polygonal, non-overlapping elements. We assume that the edges of the triangulation are oriented, so that each edge has a unit normal vector $n$ with trace values of the discontinuous functions on the edge defined as in figure~(\ref{fig:edgenormal}). Let us also define the jump and average operators as
\[
\jump{\cdot} = (\cdot)^+ - (\cdot)^-, \qquad \avg{\cdot} = \frac{1}{2}[ (\cdot)^+ + (\cdot)^-]
\]
respectively. Note that we do not involve the normal vector $n$ in the definition of these operators since it is used in the definition of the numerical fluxes. The DG scheme for an interior element $K$ is given by multiplying the Euler equation by a test function $W_h$ and doing some integration by parts to obtain
\begin{equation}
\int_{K}W_{h}\cdot\partial_{t}U(V_{h})\textrm{d}x-\int_{K}F_{i}(V_{h})\cdot\partial_{i}W_{h}\textrm{d}x+\int_{\partial K}\mathcal{H}_{c}(V_{h}^{+},V_{h}^{-},n)\cdot W_{h}^{+}\textrm{d}\sigma=0
\label{eq:dgeuler}
\end{equation}
The integral on the element boundary $\partial K$ is to be interpreted as the sum of the edge integrals, i.e., $\int_K (\cdot) \ud \sigma = \sum_{e \in \partial K} \int_e (\cdot) \ud\sigma$. Note that the flux across the element boundaries is approximated by a numerical flux function. The numerical flux is given by the concept of upwinding at the kinetic level; since molecules carry mass, momentum and energy, the flux across any edge is determined as a contribution due to molecules crossing the edge in both directions, which leads to the following numerical flux function
\begin{equation}
\mathcal{H}_{c}(V_{h}^{+},V_{h}^{-},n)=\left\langle (v\cdot n)^{+}\psi g(V_{h}^{+})\right\rangle +\left\langle (v\cdot n)^{-}\psi g(V_{h}^{-})\right\rangle \label{eq:eulkfvs}
\end{equation}
where we have used the notation
\[
(v\cdot n)^{\pm}=\frac{1}{2}\left[(v\cdot n)\pm|v\cdot n|\right]
\]
Note that each of the integrals in the two terms of the numerical flux are effectively taken over half the velocity space, i.e, $v\cdot n > 0$ and $v \cdot n < 0$, respectively. The kinetic flux can also be written as
\[
\mathcal{H}_{c}(V_{h}^{+},V_{h}^{-},n)=\frac{1}{2}\left\langle (v\cdot n)\psi g(V_{h}^{+})\right\rangle +\frac{1}{2}\left\langle (v\cdot n)\psi g(V_{h}^{-})\right\rangle +\frac{1}{2}\left\langle |v\cdot n|\psi\left[g(V_{h}^{+})-g(V_{h}^{-})\right]\right\rangle 
\]
and using
\[
g(V_{h}^{+})-g(V_{h}^{-})=\int_{0}^{1}\frac{\ud}{\ud s}g(\bar{V}(s))\ud s=\int_{0}^{1}g(\bar{V}(s))\psi\cdot\jump{V_{h}} \ud s, \qquad \bar{V}(s)=sV_{h}^{+}+(1-s)V_{h}^{-}
\]
we can write the numerical flux as
\[
\mathcal{H}_{c}(V_{h}^{+},V_{h}^{-},n)=\frac{1}{2}\left[F(V_{h}^{+}) + F(V_{h}^{-})\right]\cdot n + \frac{1}{2}\int_{0}^{1}\left\langle |v\cdot n|\psi\otimes\psi g(\bar{V}(s))\right\rangle \cdot\jump{V_{h}} \textrm{d}s
\]
This is precisely the kinetic mean-value (KMV) flux of Barth~\cite{barthcharrier2001,Timothy20063311}. Hence for the Euler equations, the KMV flux is identical to
the KFVS flux, which can be written in terms of explicit
formulae but involving error functions and exponentials~\cite{Mandal1994447}. The expressions for the one dimensional split fluxes are given in the~\ref{sec:kfvs1d}.


\subsubsection{Entropy stability of DG scheme}

To study entropy stability of the scheme, we follow Barth and substitute
$W_{h}=V_{h}$ in equation~(\ref{eq:dgeuler}) to get
\begin{equation*}
\int_{K}\partial_{t}\eta(V_{h})\textrm{d}x
 = -\int_{\partial K}\left[\Theta(V_{h}^{+},V_{h}^{-},n)+D(V_{h}^{+},V_{h}^{-},n)\right]\ud\sigma
\end{equation*}
where
\[
\Theta(V^{+},V^{-},n)=\avg{V} \cdot\mathcal{H}_{c}(V^{+},V^{-},n)-\avg{\theta^{*}\cdot n}
\]
is a consistent and conservative numerical flux for the entropy
flux $\theta\cdot n$, and
\[
D(V^{+},V^{-},n)=\frac{1}{2}\left(\jump{V} \cdot\mathcal{H}_{c}(V^{+},V^{-},n)-\jump{\theta^{*}\cdot n} \right) = \frac{1}{2}\int_0^1\jump{V} \cdot\left[\mathcal{H}_{c}(V^{+},V^{-},n)-F(\bar{V}(s))\cdot n\right]\ud s
\]
For entropy stability, we have to show that $D$ is positive; a sufficient condition  is given by~\cite{Cockburn:1998:RDG:287244.287254,Timothy20063311}
\begin{equation}
\jump{V} \cdot\left[\mathcal{H}_{c}(V^{+},V^{-},n)-F(\bar{V}(s))\cdot n\right]\ge 0,\qquad \forall s \in [0,1]\label{eq:eflux}
\end{equation}
which is a generalization to systems of equations, of Osher's E-flux
condition~\cite{osher:947} for scalar conservations laws. Barth~\cite{Timothy20063311} has shown that the KMV flux and hence also the KFVS flux, satisfies the E-flux condition; thus the KFVS-based DG scheme for the Euler equations satisfies the entropy condition for any degree of basis functions.

\subsection{Chapman-Enskog distribution function}
Consider the Boltzmann equation with the BGK model for the collision term
\begin{equation}
\partial_t f + v_i \partial_i f = \frac{g - f}{\tr}
\label{eq:bebgk}
\end{equation}
If we write the solution as a series in the relaxation time $\tr$ as in equation~(\ref{eq:ceseries}), we obtain the following first two terms of the solution
\begin{equation*}
f_0 = g, \qquad f_1 = -(\partial_t g + v_i \partial_i g)
\end{equation*}
The distribution function $g$ depends on the macroscopic variables $\rho, u, T$; we substitute the inviscid equations governing these quantities into the expression for $f_1$ to obtain the following Chapman-Enskog velocity distribution function~\cite{Chou:1997:KFS:254115.254120}
\begin{equation}
\fce = g \left[ 1 - \frac{\rho}{2p^2}\tau_{ij} c_i c_j - \frac{\rho q_i c_i}{p^2}\left(1 - \frac{2}{5} \beta |c|^2\right) \right]
\label{eq:ce}
\end{equation}
where the relaxation time is given by $\tr=\mu/p$ and $c_i = v_i - u_i$ is the peculiar velocity of the molecules. It can be checked that $U=\mom{\psi \fce}$ and $F_i + G_i = \mom{v_i \psi \fce}$, i.e., the moments of the Chapman-Enskog distribution yield the conserved variables and fluxes for the Navier-Stokes equations. Since $\fce$ is a truncated solution in the Chapman-Enskog expansion, it does not satisfy the Boltzmann equation~(\ref{eq:bebgk}); we can derive an equation for $\fce$ as $\partial_t \fce + v_i \partial_i \fce = Q$ and it can be shown through explicit but lengthy computations that $\mom{\psi Q}=0$. The Chapman-Enskog distribution is derived for a monatic gas for which Prandtl number is one; however written in the form of equation~(\ref{eq:ce}), the effect of Prandtl number is accounted in the computation of the heat flux $q$ and we take the above distribution function to be valid for polyatomic gases also~\cite{Chou:1997:KFS:254115.254120,May:2007:IGB:1223679.1223717}.

\subsection{DG scheme for Navier-Stokes equations}

In terms of the entropy variables, the Navier-Stokes equations can
be written as
\[
\tilde{A}_{0}\partial_{t}V+\tilde{A}_{i}\partial_{i}V-\partial_{i}(K_{ij}\partial_{j}V)=0,\qquad G_{i}=-K_{ij}\partial_{j}V
\]
and the matrices $K_{ij}$ are symmetric positive semi-definite~\cite{dutt1988,Shakib1991141}.
Taking dot product with entropy variables, we obtain
\begin{equation}
\partial_{t}\eta+\partial_{i}\theta_{i}+\partial_{i}(V\cdot G_{i})=-(\partial_{i}V)^{\top}K_{ij}\partial_{j}V\le0\label{eq:nsent}
\end{equation}
which is essentially the entropy condition. In order to mimic the above property, it is useful to write the DG scheme for Navier-Stokes equations using entropy variables. For an interior element $K$, the DG scheme which is the analoque of the scheme given by equation~(\ref{eq:lcddg1}) is
\begin{equation}
\begin{aligned}
&\int_{K}W_{h}\cdot\partial_{t}U(V_{h})\textrm{d}x-\int_{K}[F_{i}(V_{h}) + G_i(V_h, \nabla V_h)]\cdot\partial_{i}W_{h}\textrm{d}x \\
&+ \int_{\partial K} \mathcal{H}(V_{h}^{+},\nabla V_{h}^{+},V_h^-,\nabla V_h^-, n)  \cdot W_{h}^{+}\textrm{d}\sigma - \int_{\partial K} \mathcal{H}_d^+(V_{h}^{+},\nabla W_{h}^{+}, n)  \cdot \jump{V_{h}}\textrm{d}\sigma = 0
\end{aligned}
\label{eq:dgnselem}
\end{equation}
where $n$ is taken to be the outward normal to the element $K$. The numerical flux function $\mathcal{H}$ which consists of both inviscid and viscous fluxes is given by kinetic splitting as
\[
\mathcal{H}(V^{+},\nabla V^{+},V^{-},\nabla V^{-},n)=\mathcal{H}_{c}(V^{+},V^{-},n)+\mathcal{H}_{d}(V^{+},\nabla V^{+},V^{-},\nabla V^{-},n)
\]
Here $\mathcal{H}_{c}$ is the KFVS numerical flux for Euler equations
as given by equation \eqref{eq:eulkfvs} while $\mathcal{H}_{d}$
is the KFVS numerical flux for the viscous fluxes obtained from Chapman-Enskog
distribution. This flux can be written as
\[
\mathcal{H}_{d}(V^{+},\nabla V^{+},V^{-},\nabla V^{-},n)=\mathcal{H}_{d}^{+}(V^{+},\nabla V^{+},n) + \mathcal{H}_{d}^{-}(V^{-},\nabla V^{-},n)
\]
Note that equation~(\ref{eq:dgnselem}) is a non-symmetric DG scheme for the NS equation at the element level.
\subsubsection{Entropy stability}

In order to study entropy stability of the scheme given by equation~(\ref{eq:dgnselem}) we substitute $W_{h}=V_{h}$ to obtain
\[
\begin{aligned}\int_{K}\partial_{t}\eta(V_{h})\textrm{d}x & +\int_{\partial K}\left[\Theta(V_{h}^{+},V_{h}^{-},n)+D(V_{h}^{+},V_{h}^{-},n)\right]\textrm{d}\sigma\\
 & +\int_{\partial K}\left[\mathcal{H}_{d}^{+}(V_{h}^{+},\nabla V_{h}^{+},n)\cdot V_{h}^{-}+\mathcal{H}_{d}^{-}(V_{h}^{-},\nabla V_{h}^{-},n)\cdot V_{h}^{+}\right]\textrm{d}\sigma -\int_{K}G_{i}(V_{h},\nabla V_{h})\cdot\partial_{i}V\textrm{d}x=0
\end{aligned}
\]
The first two integrals are common with the Euler equations, while the last two terms are the contributions from NS equations. We observe that
\[
\mathcal{H}_{d}^{+}(V_{h}^{+},\nabla V_{h}^{+},n)\cdot V_{h}^{-}+\mathcal{H}_{d}^{-}(V_{h}^{-},\nabla V_{h}^{-},n)\cdot V_{h}^{+}
\]
is a consistent and conservative numerical flux for the viscous entropy
flux $V\cdot(G_{i}n_{i})$, while from the Navier-Stokes equations,
we have
\[
G_{i}(V_h,\nabla V_h)\cdot\partial_{i}V_h = - (\partial_{i}V_h)^\top K_{ij} \partial_j V_h \le 0
\]
which leads to a cell entropy inequality for the DG scheme which mimics
equation \eqref{eq:nsent}. The non-symmetric nature of the scheme~(\ref{eq:dgnselem}) is required to obtain the above cell entropy inequality. This type of local or cell entropy condition does not seem to be valid in the case of the symmetric version of the scheme.
\subsubsection{DG scheme for NS equation}
We now state the DG scheme for NS equation in a global formulation including penalty terms which is analogous to the scheme given in equation~(\ref{eq:lcddgip}). Let $\Gi$ denote the set of all interior edges/faces of the finite element grid and let $\Gamma$ denote the boundary faces. The time continuous DG scheme for Navier-Stokes equations is given by
\begin{equation}
\begin{aligned}
&\int_{\Omega}W_{h} \cdot \partial_{t}U(V_{h})\textrm{d}x  - \int_{\Omega} \left[F_{i}(V_{h})+G_{i}(V_{h},\nabla V_{h})\right]\cdot\partial_{i}W_{h}\textrm{d}x
 + \int_{\Gi}\mathbb{\mathcal{H}}(V_{h}^{+},\nabla V_{h}^{+},V_{h}^{-},\nabla V_{h}^{-},n)\cdot \jump{W_{h}} \textrm{d}\sigma \\
 +& \epsilon\int_{\Gi}\mathcal{H}_{d}(V_{h}^{+},\nabla W_{h}^{+},V_h^-,\nabla W_h^-, n)\cdot \jump{V_{h}} \textrm{d}\sigma 
+ \int_{\Gi} \delta_h(V_h) \cdot \jump{W_h} \ud \sigma + N_\Gamma(V_h,W_h) = 0
\end{aligned}
\label{eq:dgns}
\end{equation}
Choosing $\epsilon=-1$ yields the NIPG scheme while $\epsilon=+1$ yields the SIPG scheme. 
The interior penalty term is of a similar form as in the scalar problem and given by
\begin{equation*}
\delta_h(V_h) = \cip \frac{k^2 \nu}{h} \jump{U(V_h)}
\end{equation*}
where we use the coefficient of kinematic viscosity in order to have dimensional consistency. The boundary conditions are imposed through the boundary terms $N_\Gamma$ in a weak form~\cite{Hartmann20089670} and given by
\begin{equation*}
\begin{aligned}
N_\Gamma(V_h, W_h) =& \int_{\Gamma}\mathcal{H}(\Gamma V_{h}^{+}, \Gamma\nabla V_{h}^{+}, \Gamma V_{h}^{-}, \Gamma\nabla V_{h}^{-},n)\cdot W_{h}^+ \textrm{d}\sigma \\
&+ \int_{\Gamma}\mathcal{H}_d(\Gamma V_{h}^{+},\Gamma\nabla W_{h}^{+}, \Gamma V_{h}^{-}, \Gamma\nabla W_{h}^{-},n)\cdot [V_{h}^+ - \Gamma V_h^+]  \textrm{d}\sigma + \int_\Gamma \delta_\Gamma(V_h) \cdot W_h^+ \ud\sigma
\end{aligned}
\end{equation*}
The operator $\Gamma$ modifies the boundary trace of the solution to apply the boundary conditions; for example, in the case of a no-slip boundary
\begin{equation*}
\Gamma V_h^+ = \Gamma V_h^- = \left[
V_1^+, 
0, 
0, 
0, 
V_5^+
\right]^\top, \qquad \Gamma\nabla V_h^+ = \Gamma\nabla V_h^- = \nabla V_h^+
\end{equation*}
For an isothermal wall, $V_5^+$ is modified by using the specified wall temperature. For an adiabatic wall, $\nabla V_h^+$ is modified so that the heat flux is zero while at a farfield boundary with conditions $V_\infty$, we take $\Gamma V_h^+ = V_h^+$, $\Gamma V_h^- = V_\infty$ and $\Gamma \nabla V_h^+ = \Gamma \nabla V_h^- = \nabla V_h^+$. The boundary penalty term $\delta_\Gamma$ is similar to the interior penalty and is of the form
\[
\delta_\Gamma(V_h) = \cip \frac{k^2 \nu}{h} [U( V_h^+) - U(\Gamma V_h^+)]
\]

\subsection{Test case: Order of accuracy}
\begin{figure}
\begin{center}
\begin{tabular}{cc}
\includegraphics[width=0.48\textwidth]{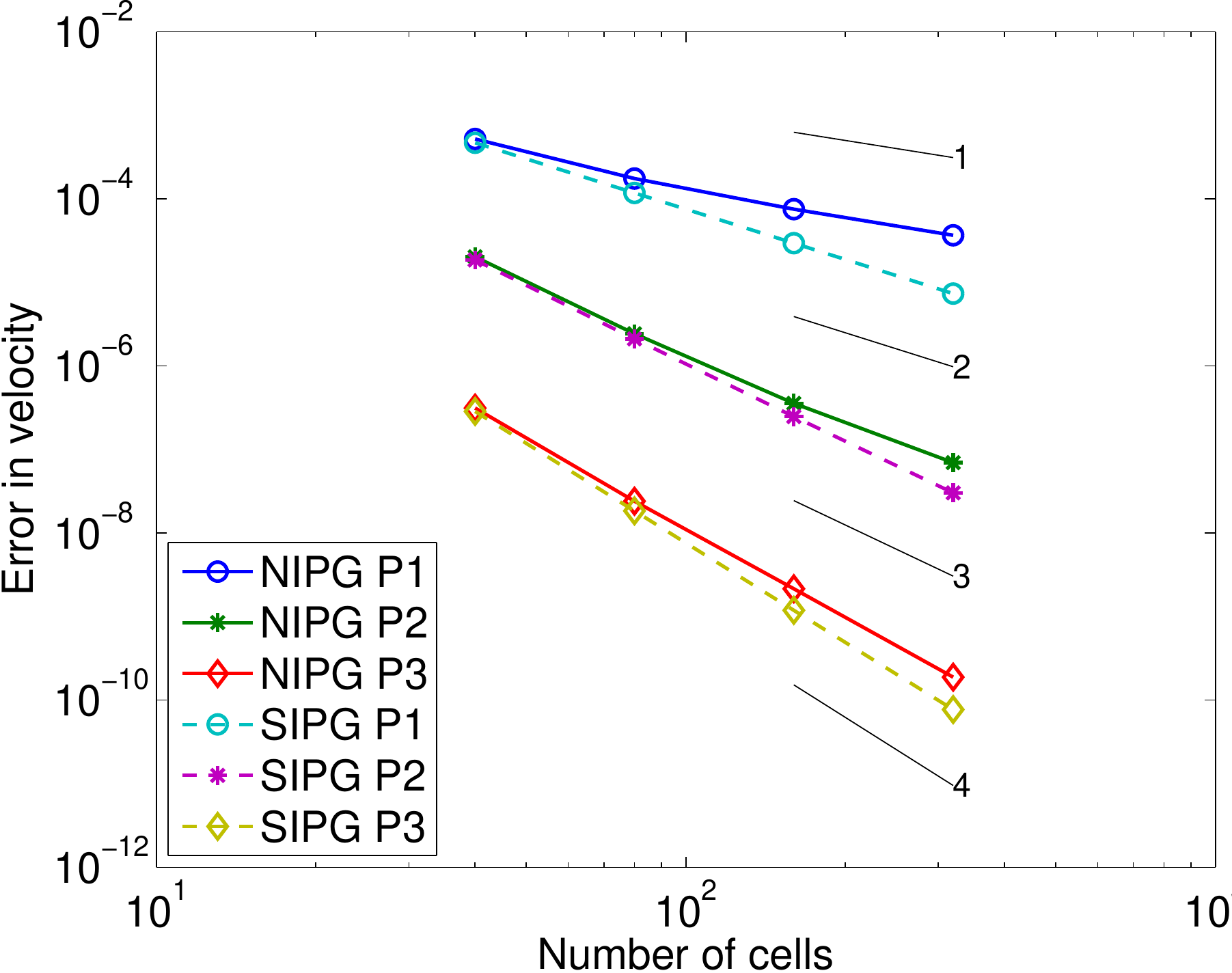} &
\includegraphics[width=0.48\textwidth]{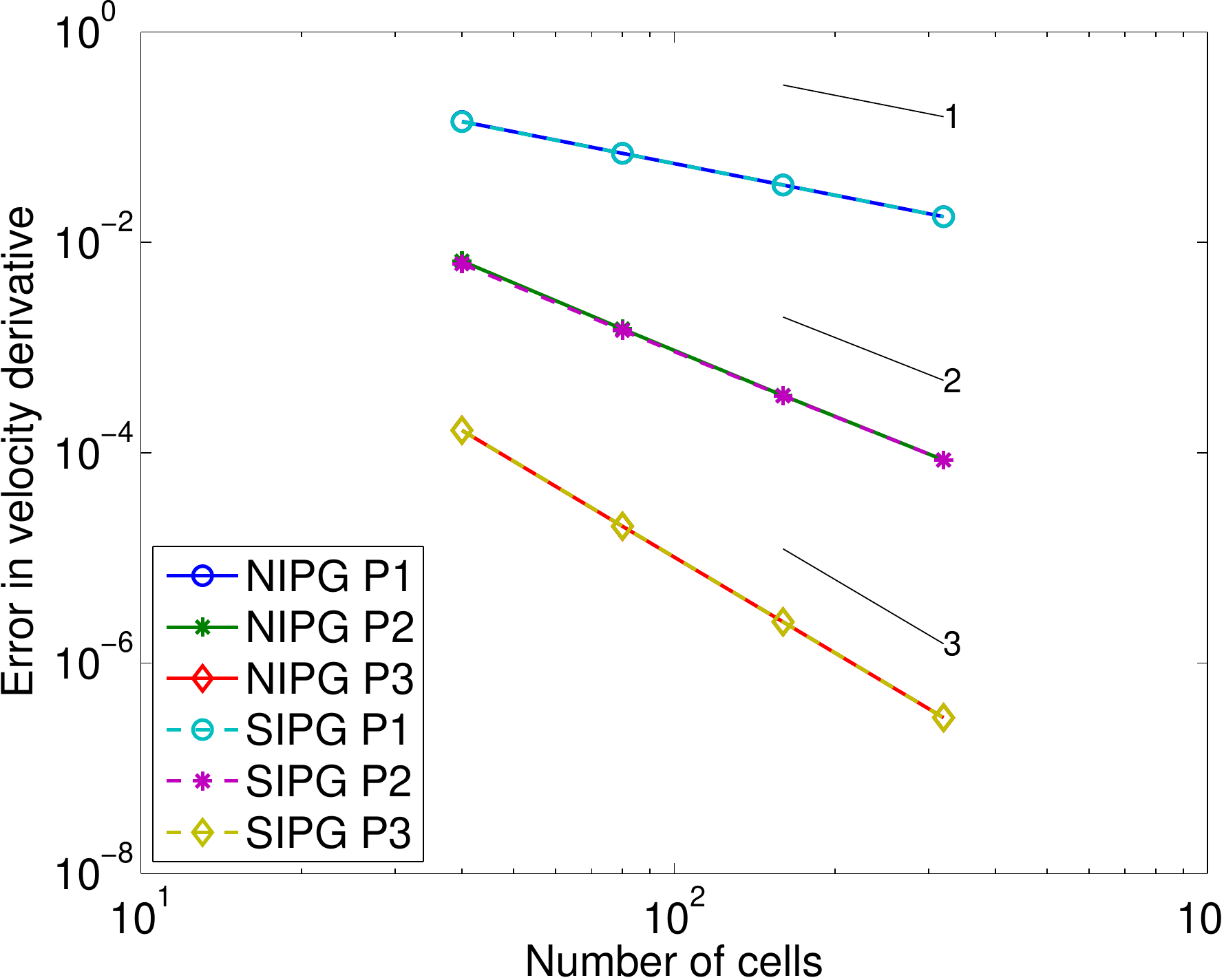} \\
(a) & (b) \\
\includegraphics[width=0.48\textwidth]{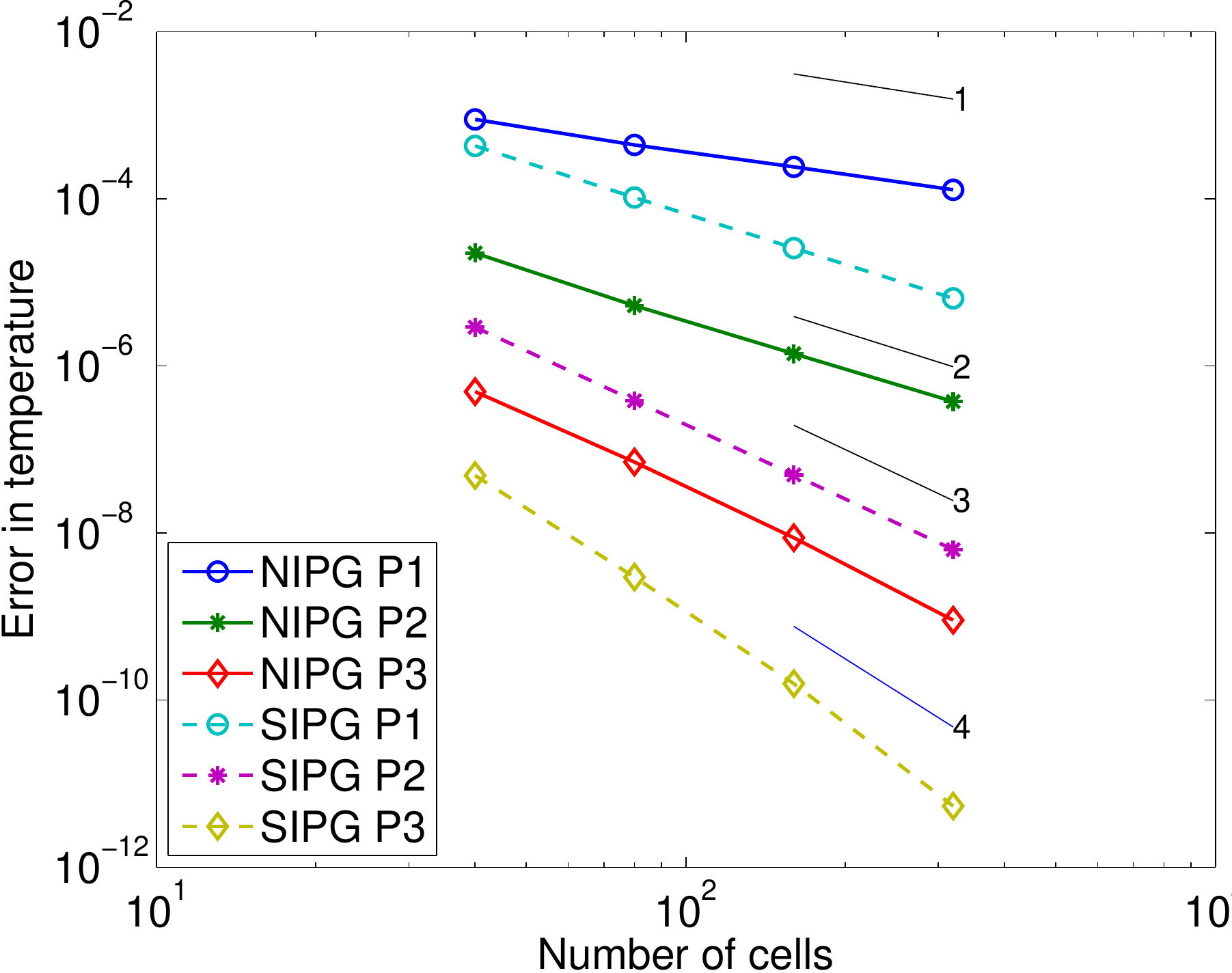} &
\includegraphics[width=0.48\textwidth]{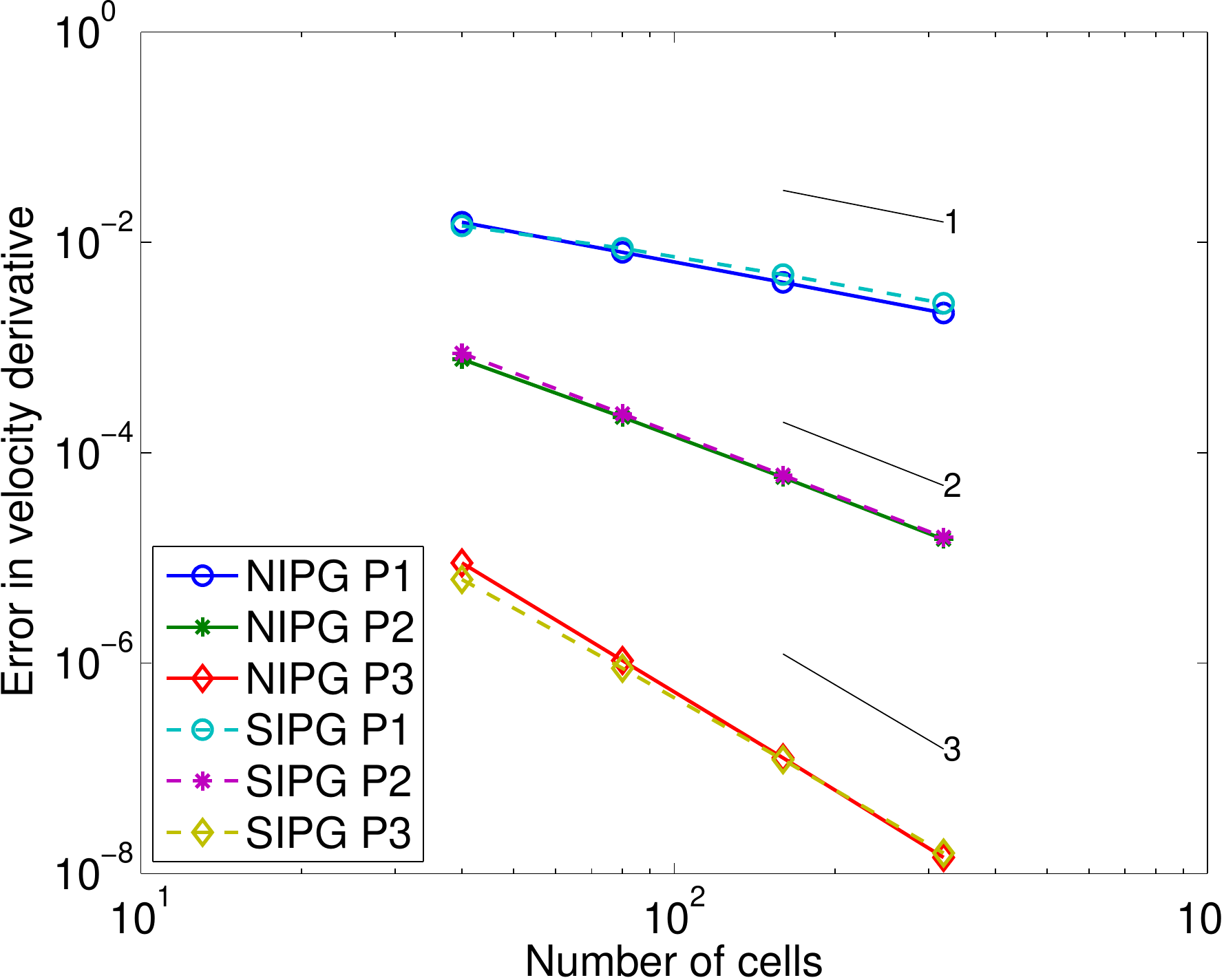} \\
(c) & (d)
\end{tabular}
\caption{Order of accuracy study for 1-D NS, $\mu=0.01$: (a) Velocity (b) Velocity derivative (c) Temperature (d) Temperature derivative}
\label{fig:nsordmu0p01}
\end{center}
\end{figure}

\begin{figure}
\begin{center}
\begin{tabular}{cc}
\includegraphics[width=0.48\textwidth]{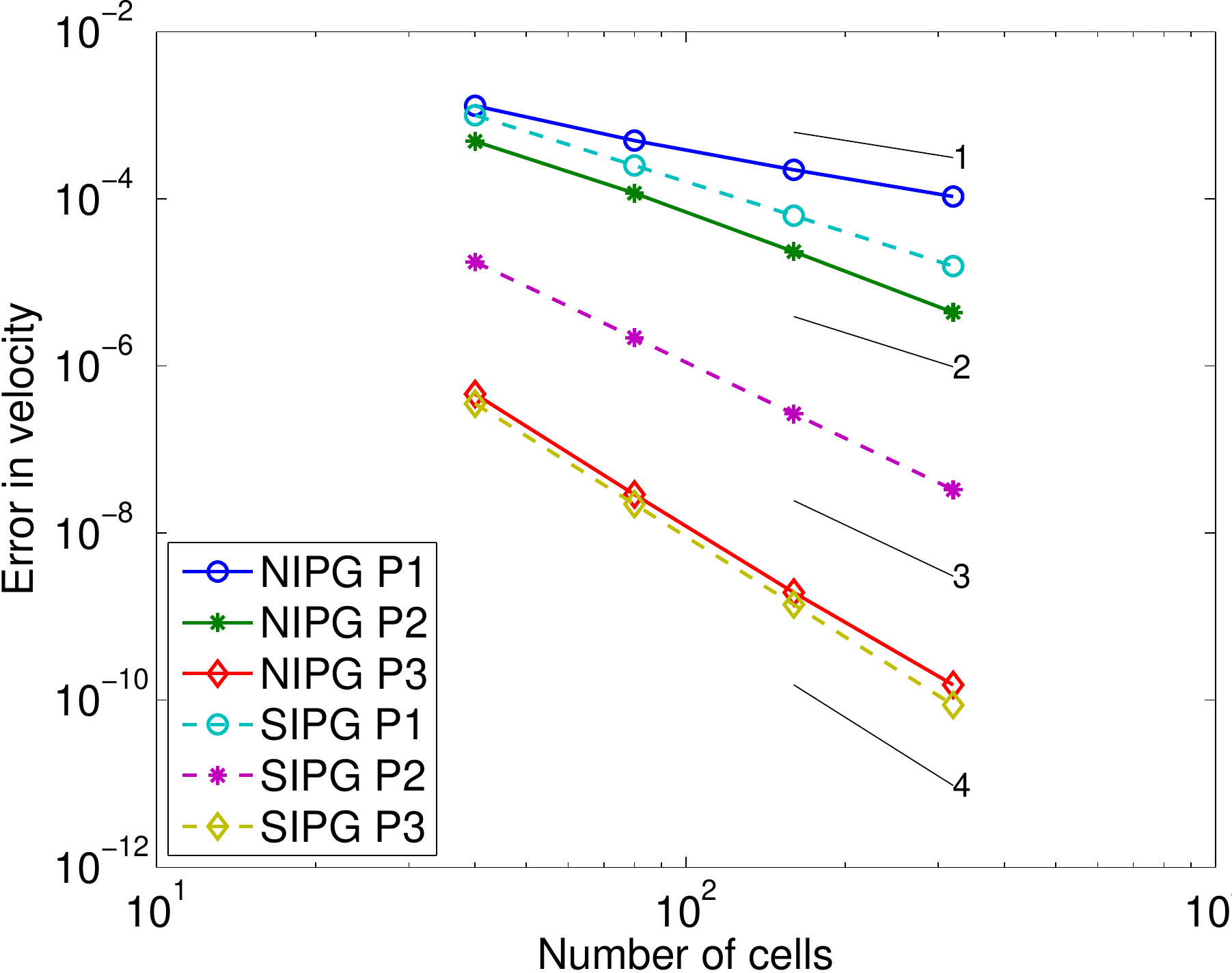} &
\includegraphics[width=0.48\textwidth]{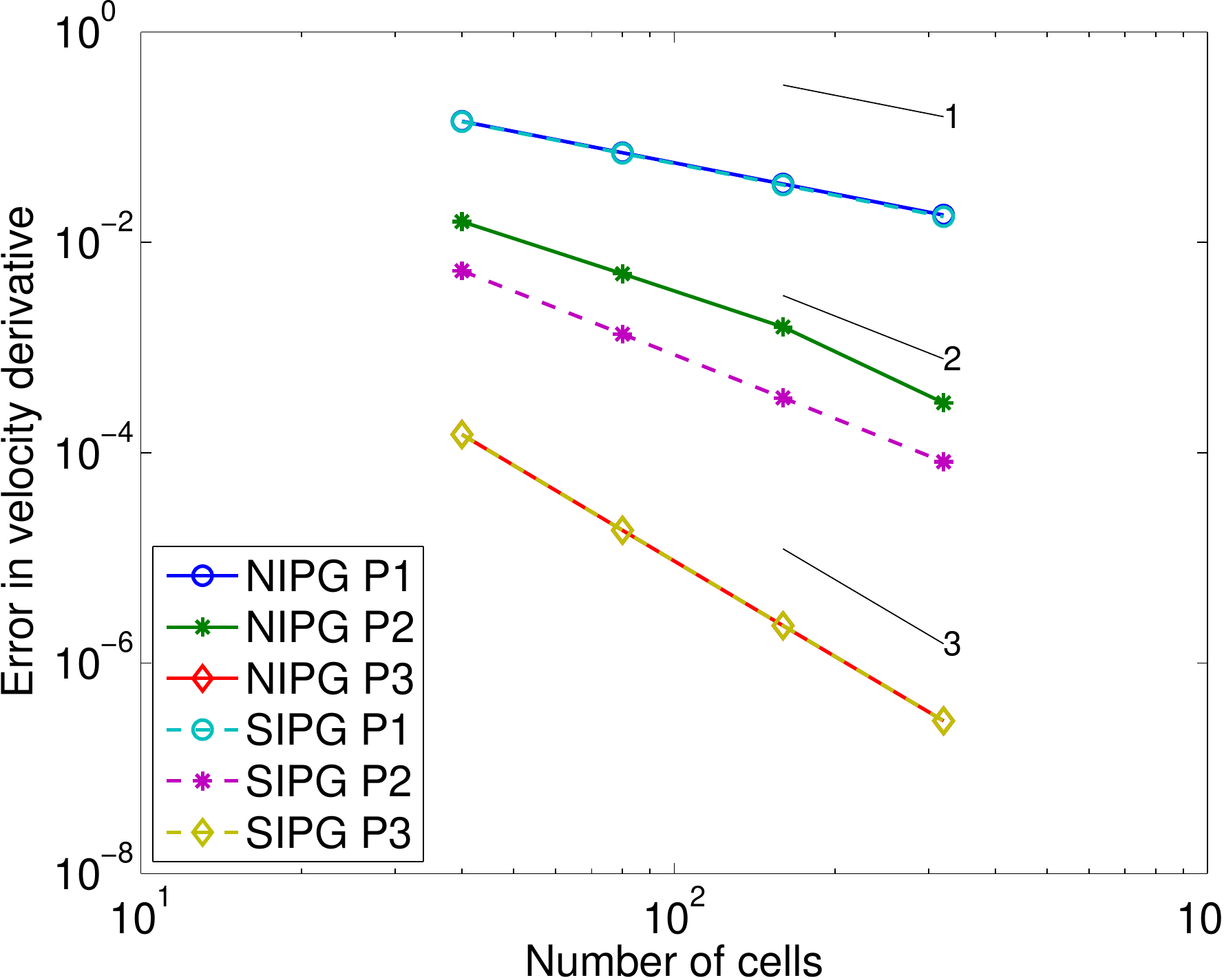} \\
(a) & (b) \\
\includegraphics[width=0.48\textwidth]{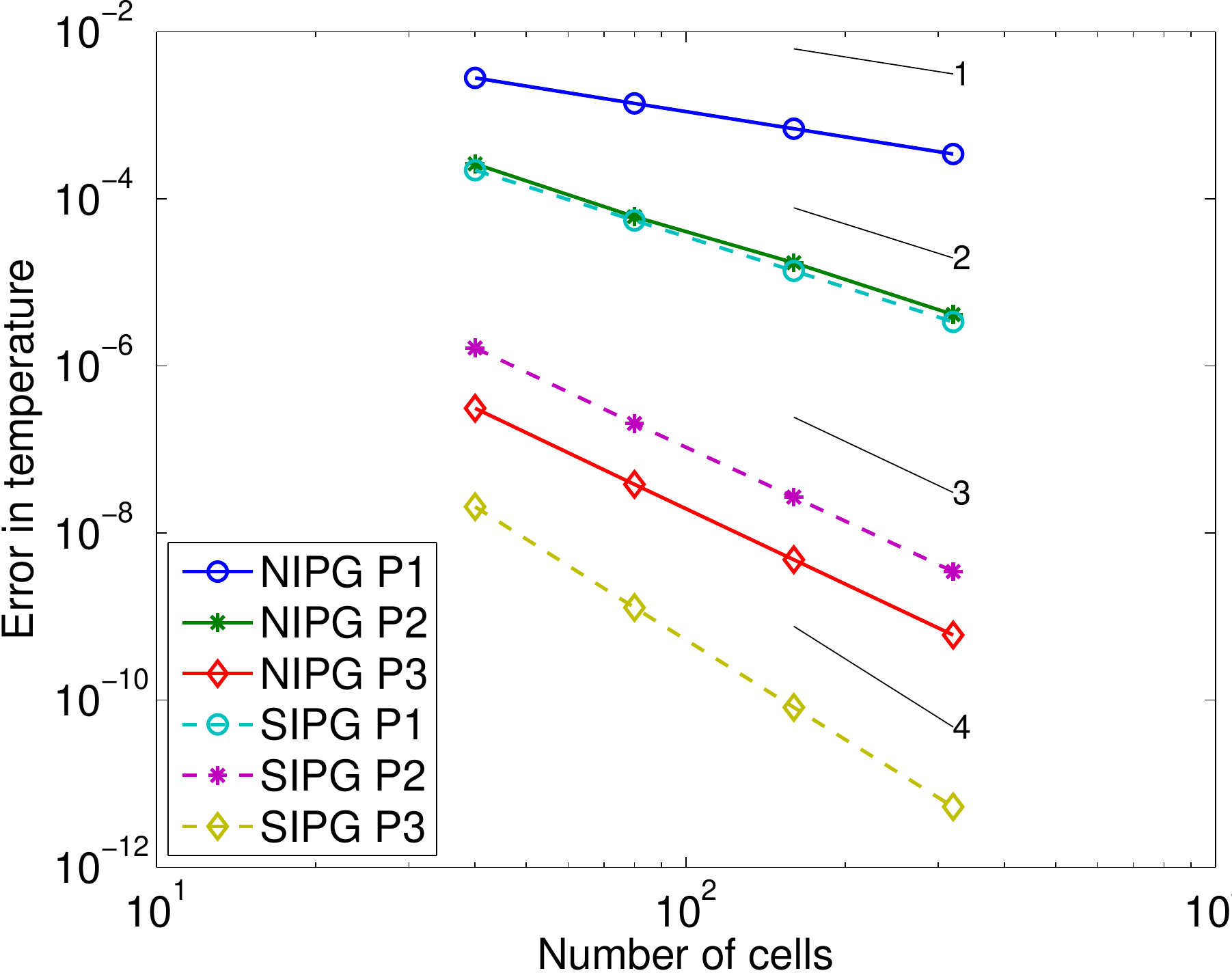} &
\includegraphics[width=0.48\textwidth]{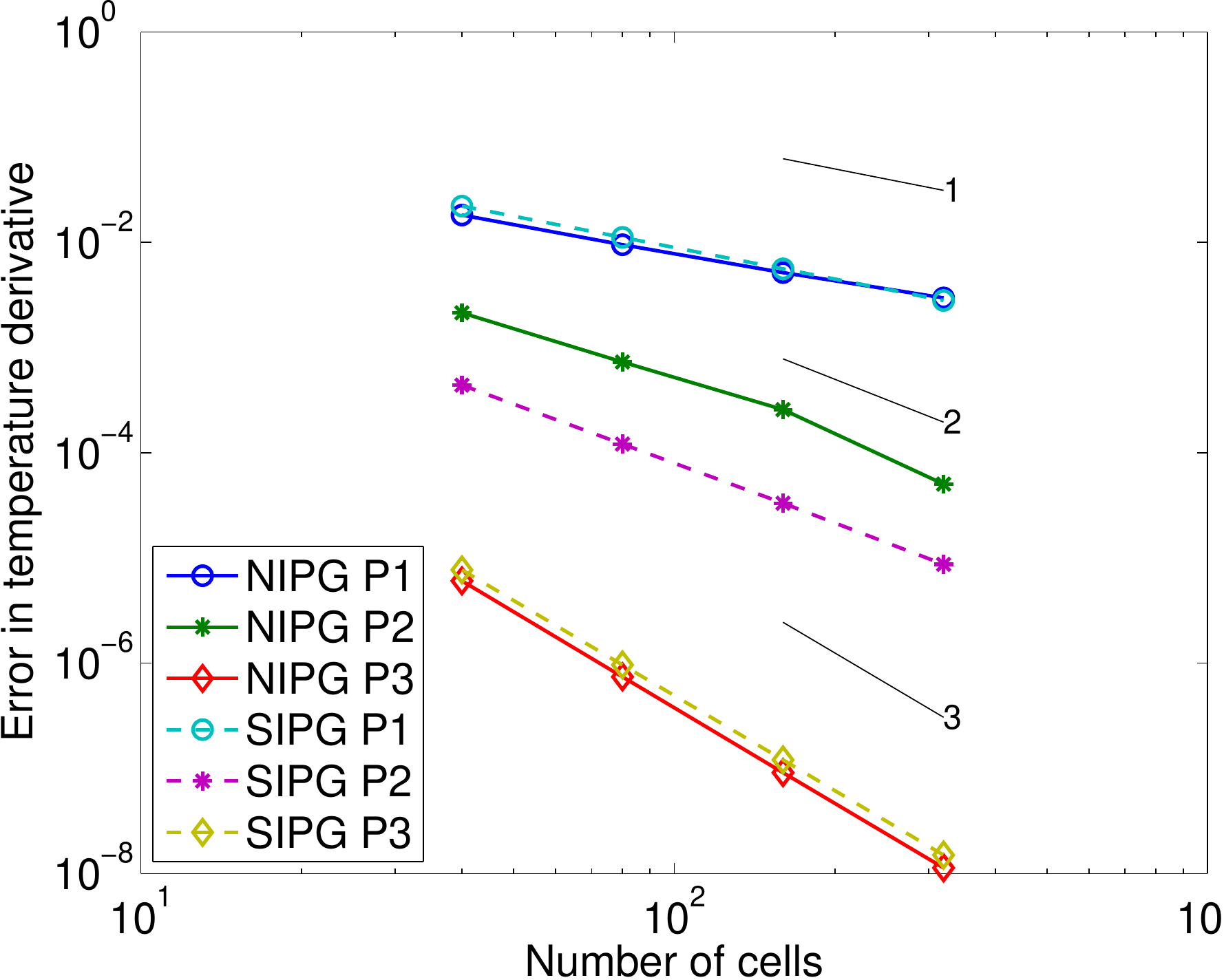} \\
(c) & (d)
\end{tabular}
\caption{Order of accuracy study for 1-D NS, $\mu=1.0$: (a) Velocity (b) Velocity derivative (c) Temperature (d) Temperature derivative}
\label{fig:nsordmu1}
\end{center}
\end{figure}

The NS equations do not have 1-D analytical solutions which could be used to compute the convergence rate of the numerical schemes. Hence we construct exact solutions by the method of manufactured solutions. We assume the following form for the exact solution
\begin{equation}
\rho(x) = 1 + \frac{1}{2}\cos(2\pi x), \quad u(x)=10 x^2(1-x)^2\sin(2\pi x), \quad T(x)=1+2x^2(1-x)^2
\end{equation}
which solves the following stationary 1-D NS equations with source terms
\begin{equation}
(\rho u)_x = f_1, \quad (p+\rho u^2)_x - \tau_x = f_2, \quad ((\rho e+p)u)_x - (\tau u)_x + q_x = f_3
\end{equation}
together with the boundary conditions $u(0)=u(1)=0$ and $q(0)=q(1)=0$, which are the adiabatic boundary conditions corresponding to zero heat flux. The source terms $f_1, f_2, f_3$ are computed by substituting the exact solutions into the above equations. The viscosity coefficient is $\mu=0.01$ or $\mu=1$ which corresponds to convection dominated and viscous dominated cases, while the ratio of specific heats $\gamma=7/5$ corresponding to polyatomic gas like air.  Computations are performed on a sequence of four uniform grids with $N = 40, 80, 160$ and 320 elements using $P_1$, $P_2$ and $P_3$ basis functions. The $L_2$ norm of the error in velocity and temperature and their derivatives are computed from the exact solution by performing a high order quadrature. The convergence of the error in velocity and temperature and their derivatives are plotted in figure~(\ref{fig:nsordmu0p01}) and (\ref{fig:nsordmu1}) for the case of $\mu=0.01$ and $\mu=1$ respectively. The short lines are drawn with slopes indicated beside them. From these figures we can observe that for the NIPG scheme, the solution and its derivative converge approximately at the rate of $O(h^k)$. For the SIPG scheme, the solution converges at the rate $O(h^{k+1})$ while the derivative converges at the rate $O(h^k)$. The SIPG scheme is slightly better in terms of the error in the solution while both schemes yield nearly similar solution for the derivatives.

\subsection{Test case: 1-D shock structure} 
\begin{figure}
\begin{center}
\begin{tabular}{cc}
N=100, $P_1$ & N=100, $P_2$ \\
\includegraphics[width=0.40\textwidth]{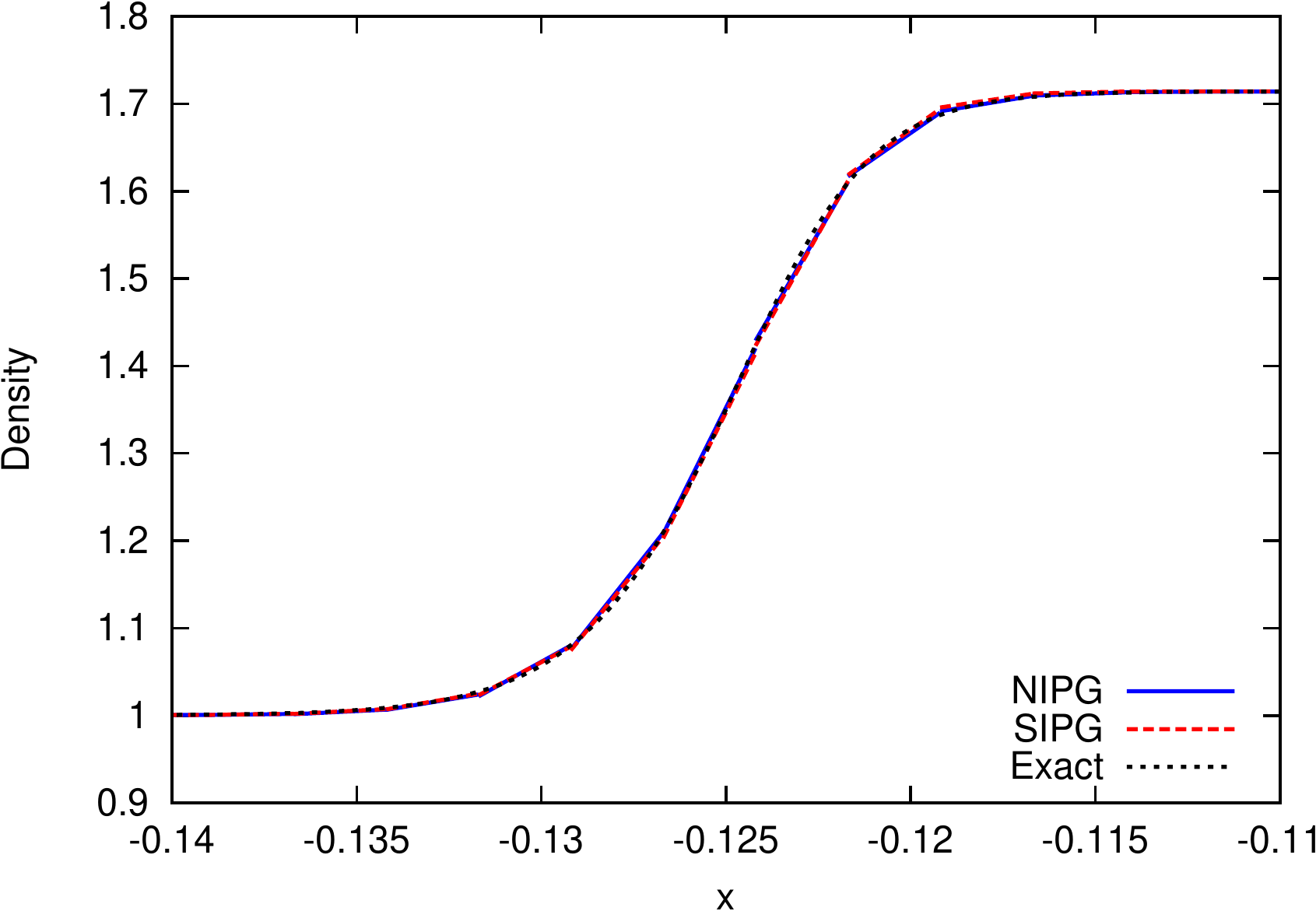} &
\includegraphics[width=0.40\textwidth]{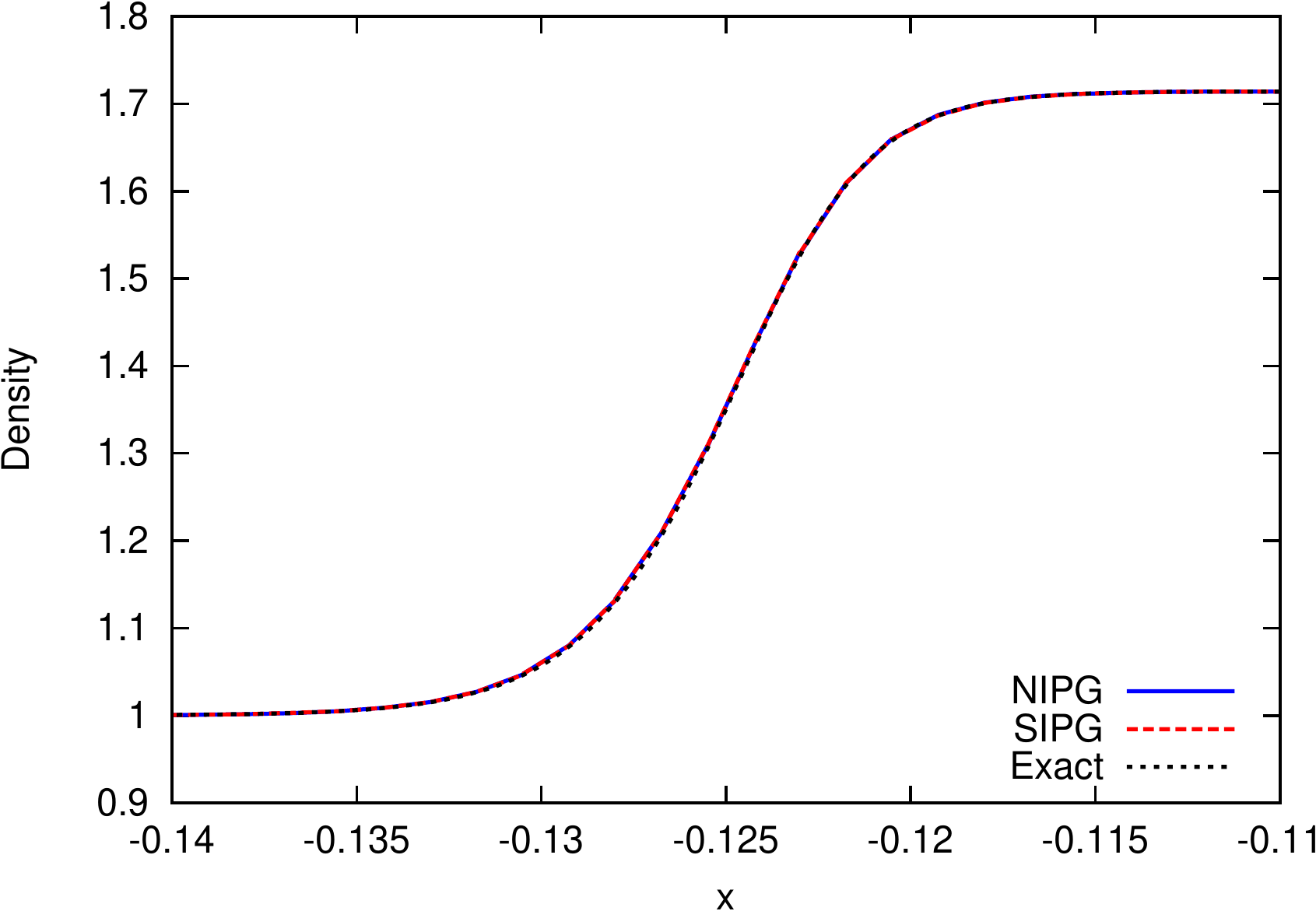} \\
N=200, $P_1$ & N=200, $P_2$ \\
\includegraphics[width=0.40\textwidth]{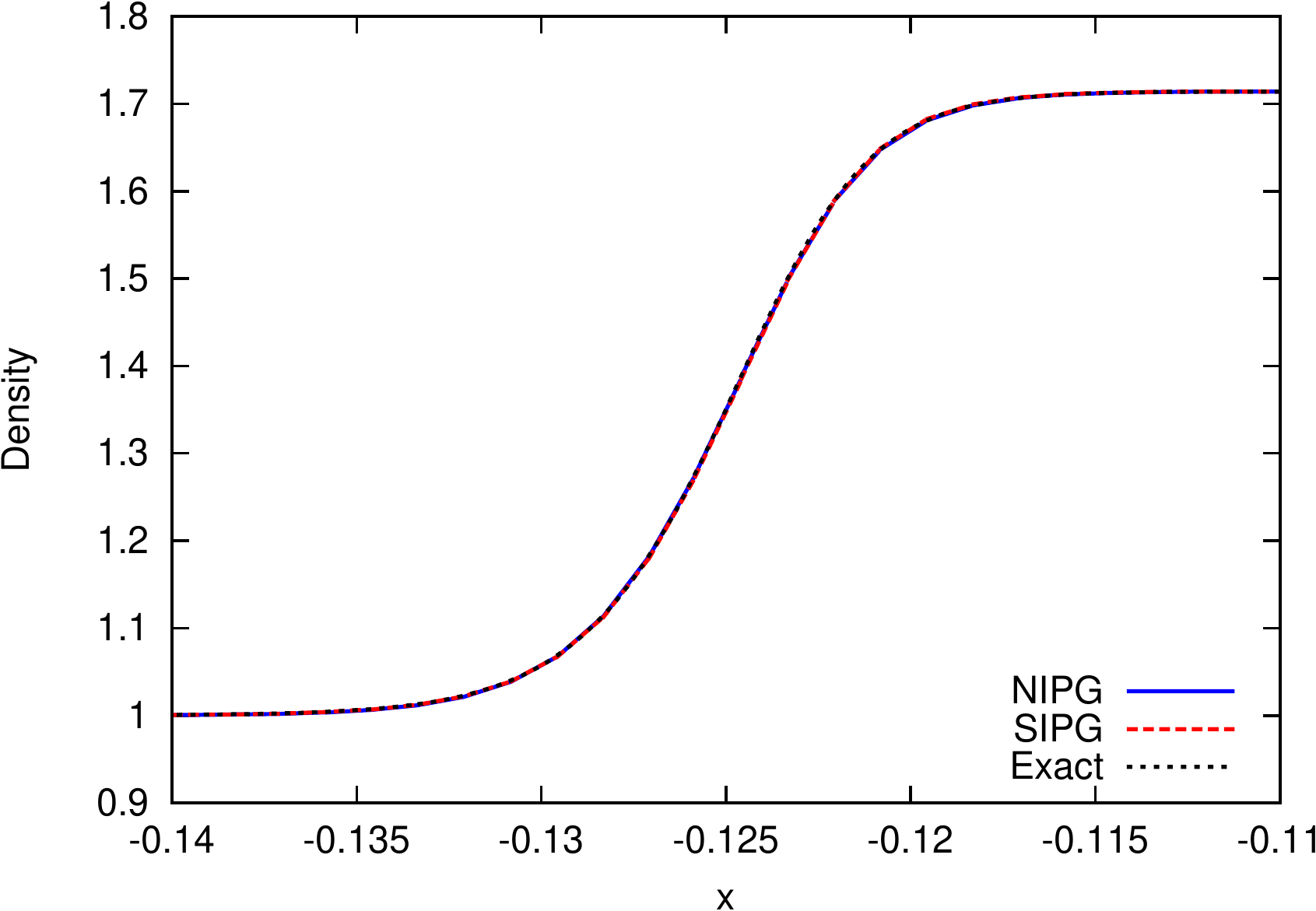} &
\includegraphics[width=0.40\textwidth]{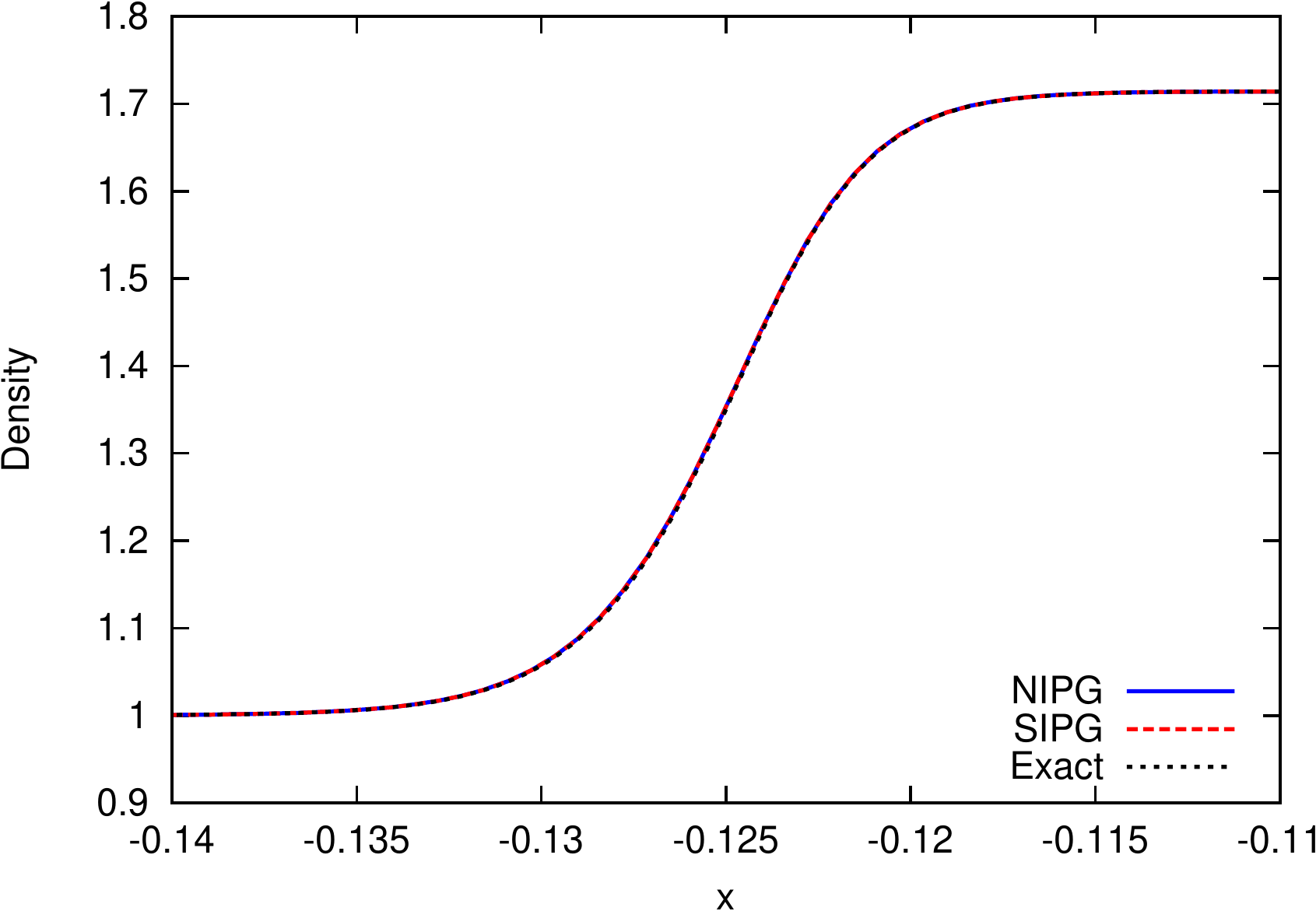}
\end{tabular}
\end{center}
\caption{Density for shock structure problem using entropy variables}
\label{fig:qrho}
\end{figure}

\begin{figure}
\begin{center}
\begin{tabular}{cc}
N=100, $P_1$ & N=100, $P_2$ \\
\includegraphics[width=0.40\textwidth]{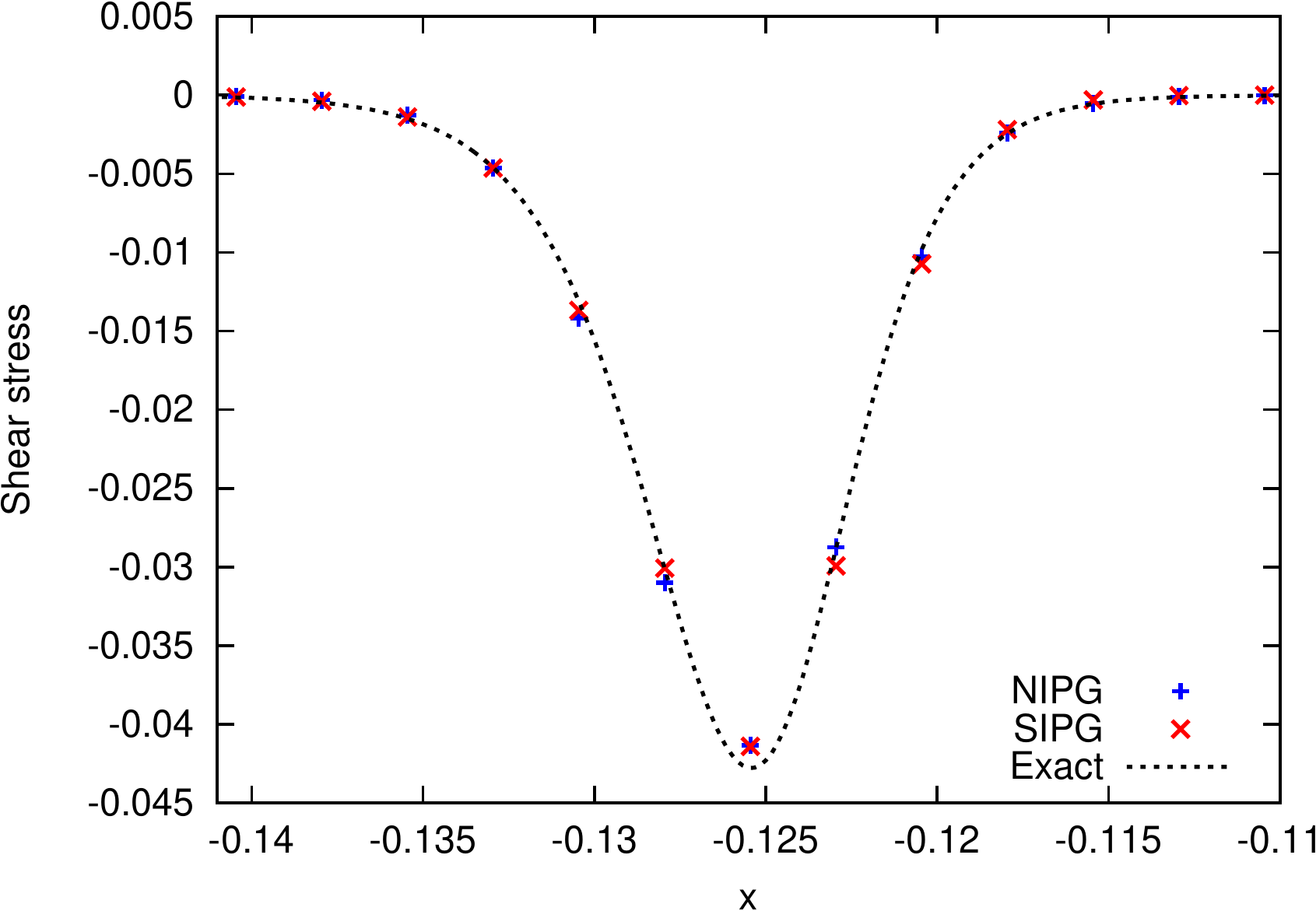} &
\includegraphics[width=0.40\textwidth]{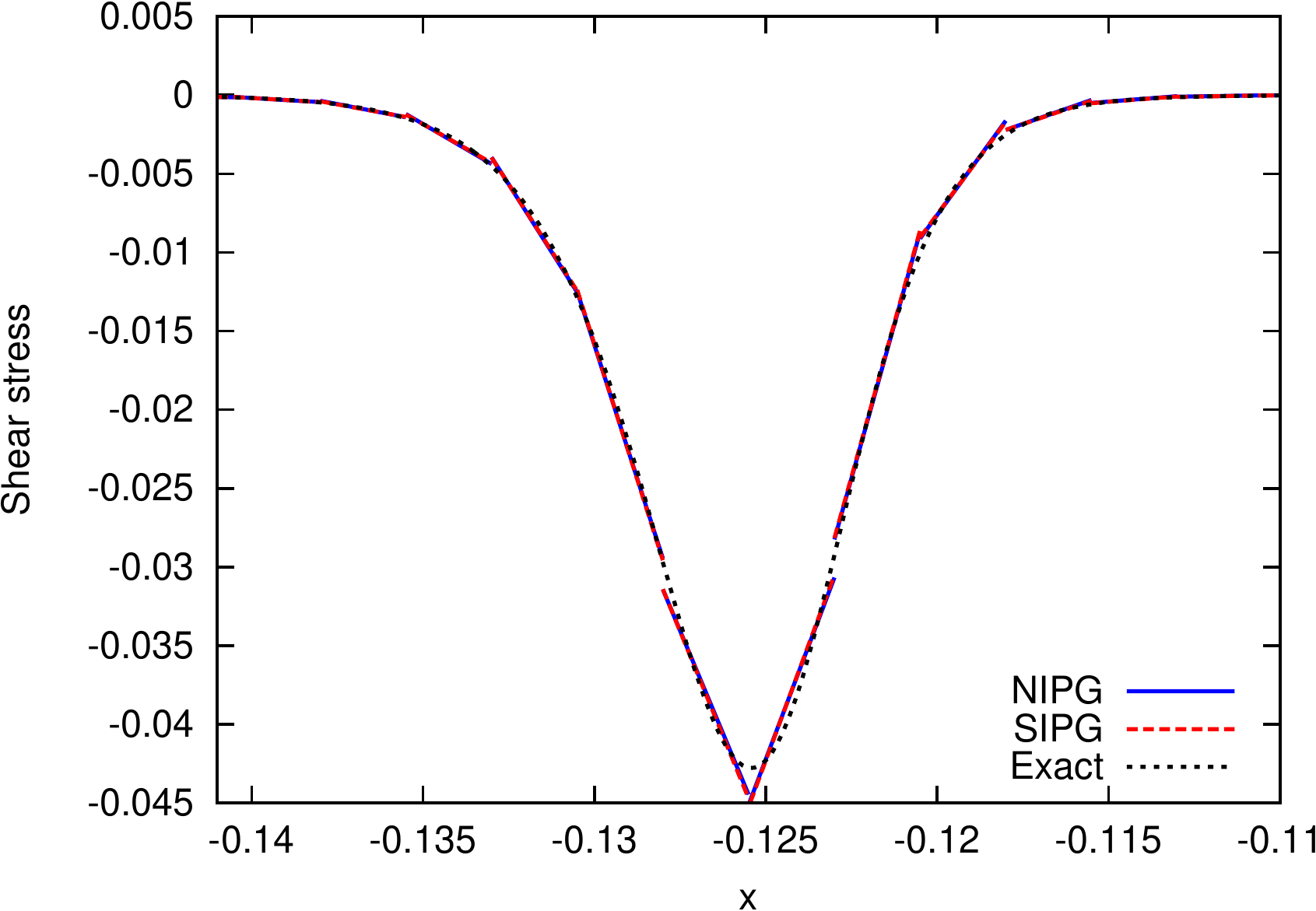} \\
N=200, $P_1$ & N=200, $P_2$ \\
\includegraphics[width=0.40\textwidth]{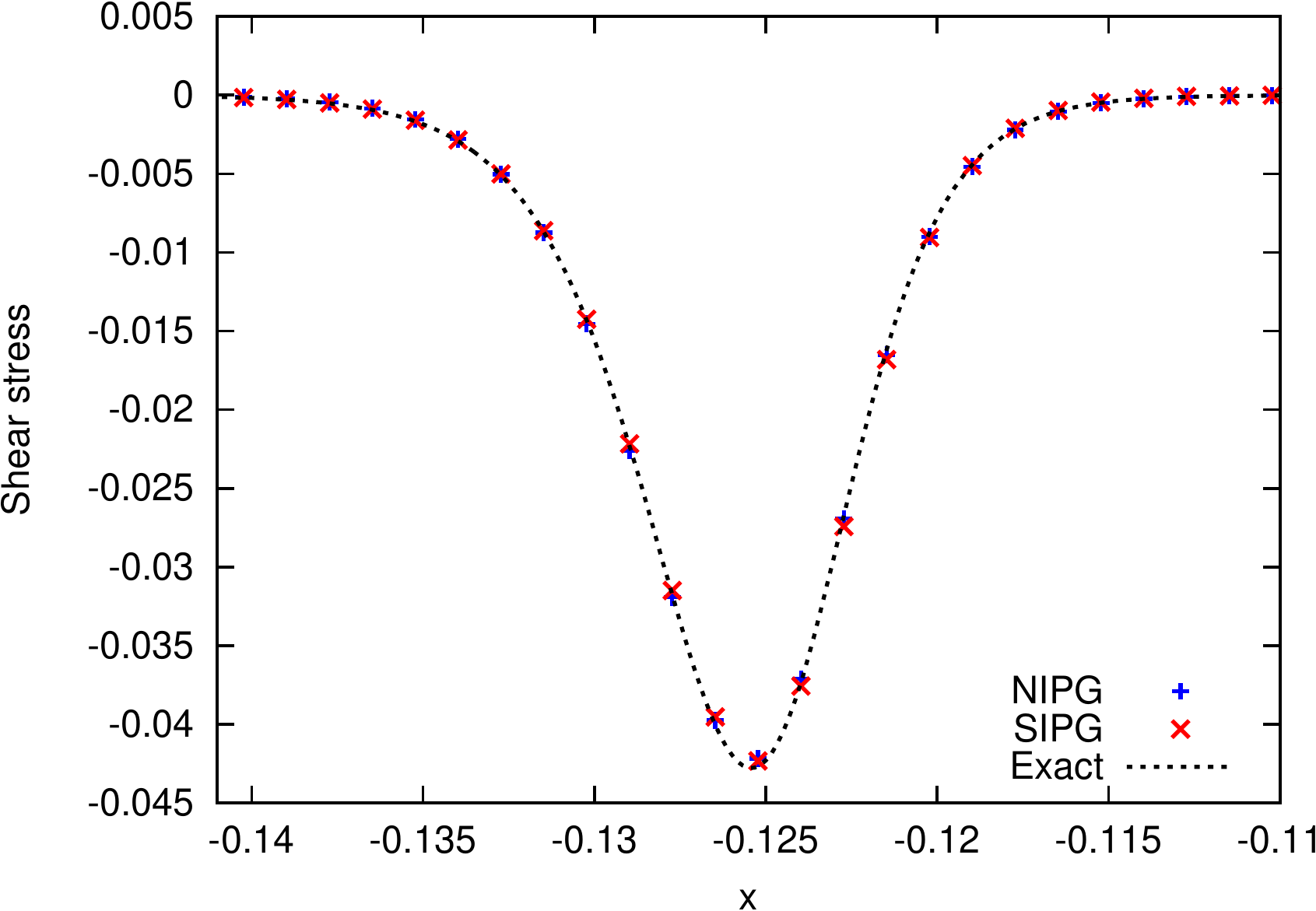} &
\includegraphics[width=0.40\textwidth]{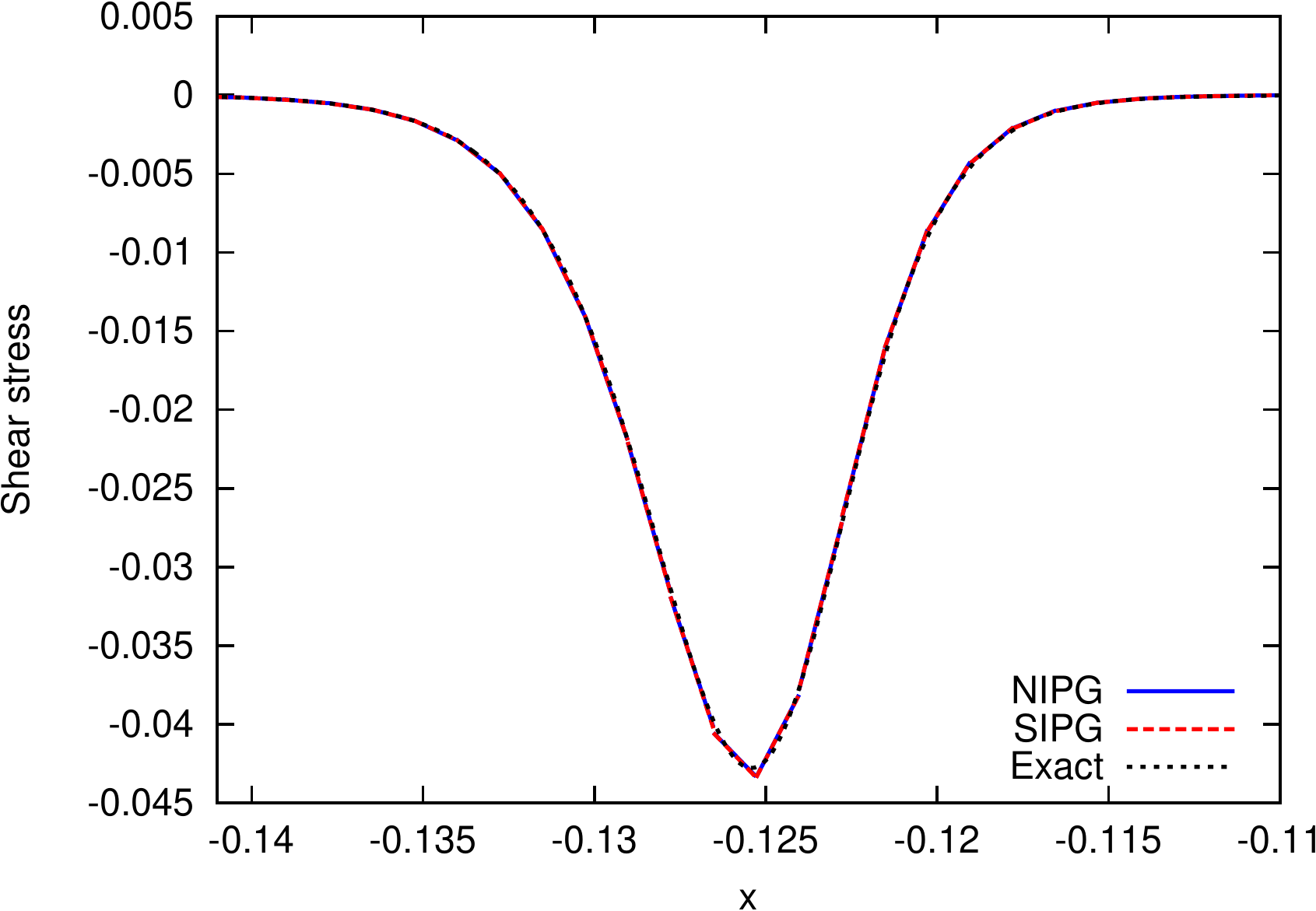}
\end{tabular}
\end{center}
\caption{Shear stress for shock structure problem using entropy variables}
\label{fig:qtau}
\end{figure}

\begin{figure}
\begin{center}
\begin{tabular}{cc}
N=100, $P_1$ & N=100, $P_2$ \\
\includegraphics[width=0.40\textwidth]{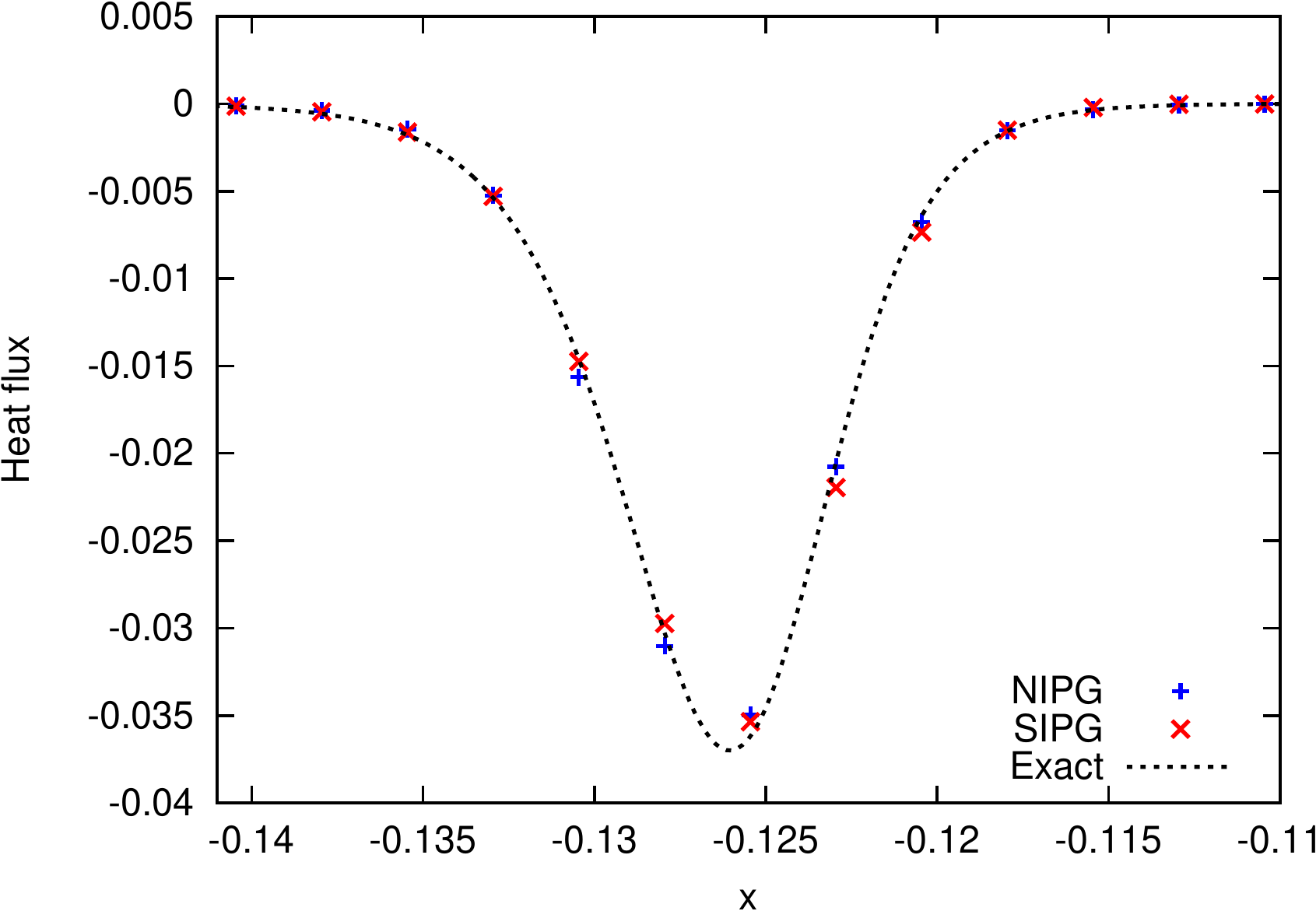} &
\includegraphics[width=0.40\textwidth]{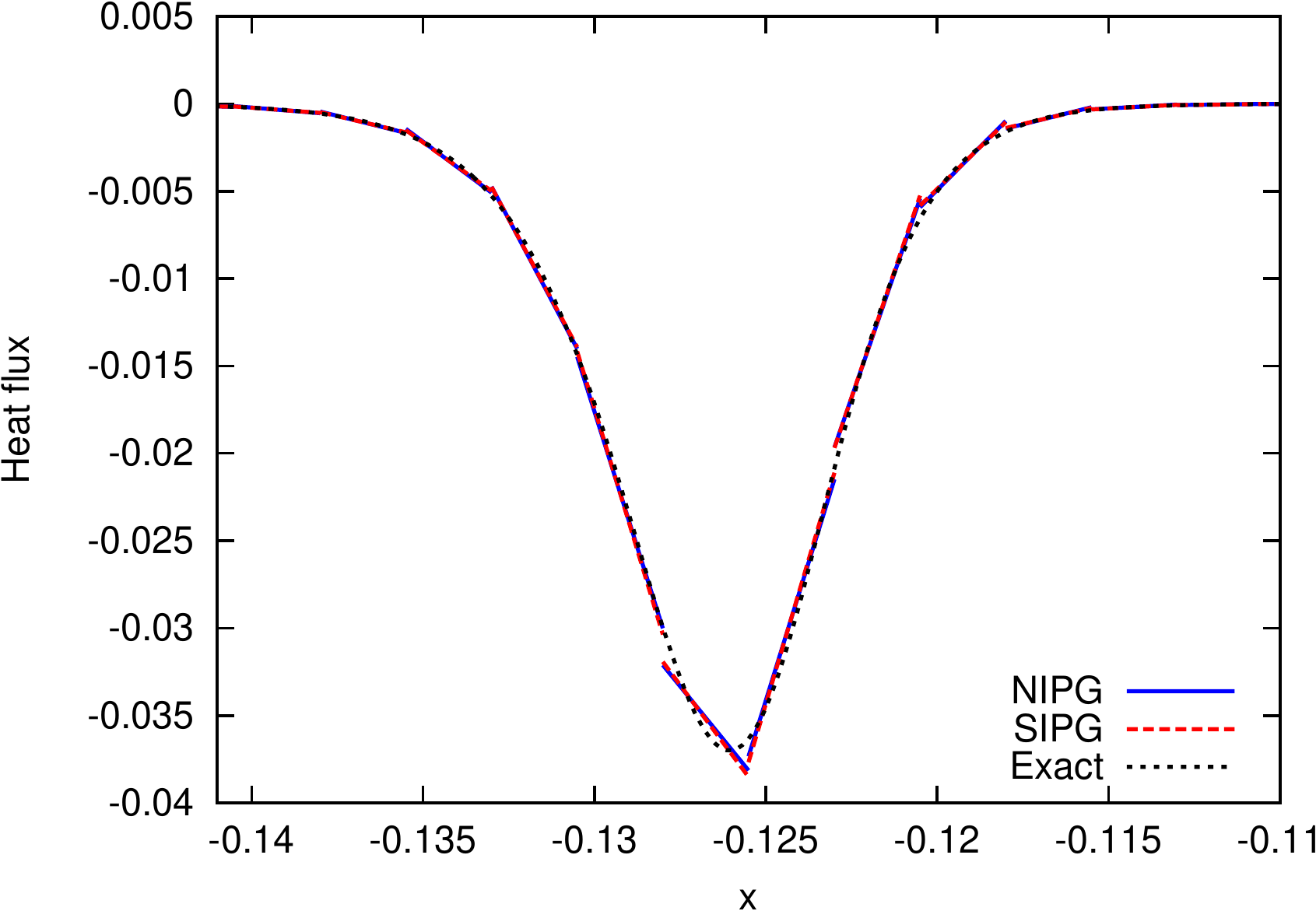} \\
N=200, $P_1$ & N=200, $P_2$ \\
\includegraphics[width=0.40\textwidth]{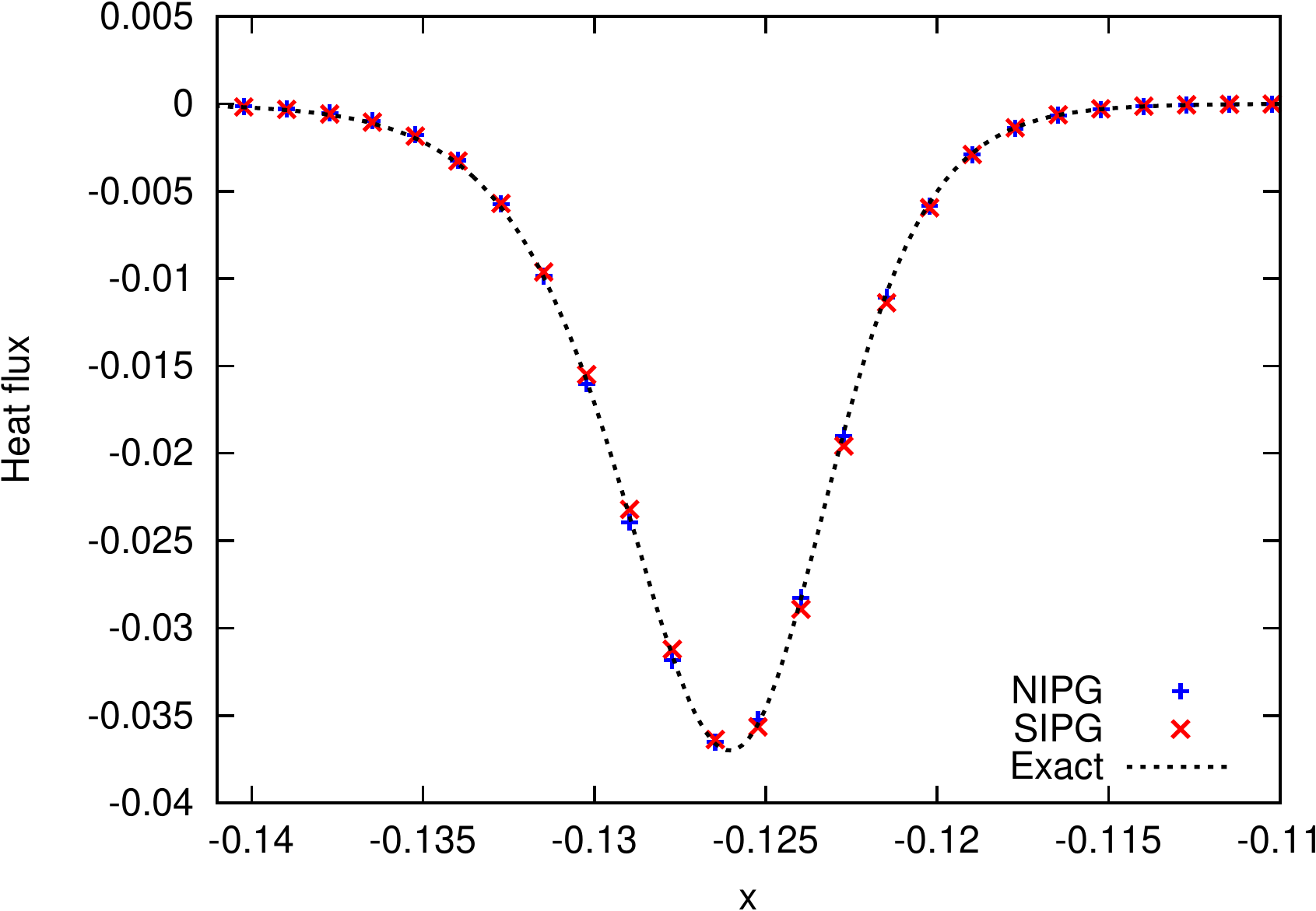} &
\includegraphics[width=0.40\textwidth]{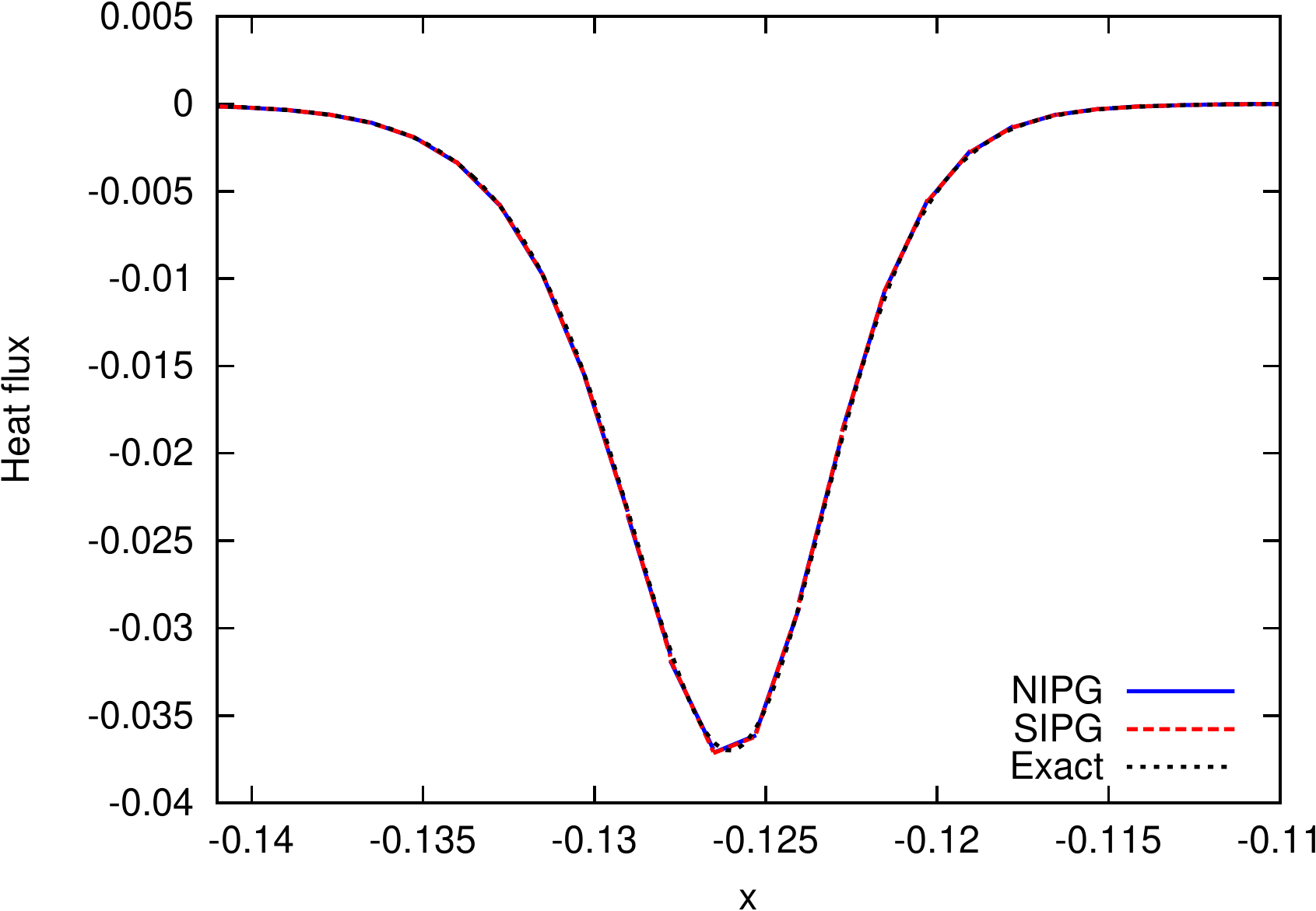}
\end{tabular}
\end{center}
\caption{Heat flux for shock structure problem using using entropy variables}
\label{fig:qq}
\end{figure}

We compute the shock structure using the Navier-Stokes equations. While it is known that the NS equations are not valid inside the shock, we use this as a verification test case. The exact solution of the NS equations for a stationary shock are obtained using the method of~\cite{gilbarg1953} which is implemented  as a matlab code in~\cite{Xu2001289}. This solution is not truly exact since it uses an ODE solver, and hence we are not able to compute errors in the DG solution. The parameters defining the problem are: Mach number ahead of the shock is $M_1=1.5$, $\gamma=5/3$, $\Pr = 2/3$ while the viscosity law is given by $\mu = \mu_1 (T/T_1)^{0.8}$ where the subscript ``1" denotes pre-shock conditions and $\mu_1=0.0005$. The time discretization is made with backward Euler implicit scheme and the time step is chosen based on the convective speeds as $\Delta t = \cfl \cdot \frac{h}{\max(|u|+c)}$ with $\cfl=5$. Figure~(\ref{fig:qrho}) shows the results for the density using NIPG/SIPG schemes and $P_1/P_2$ basis functions on grids with $N=100, h=1/400$ and $N=200,h=1/800$ elements, while figures~(\ref{fig:qtau}) and (\ref{fig:qq}) show the corresponding results for shear stress and heat flux. It is clear that the numerical solutions are able to compute the shock structure quite accurately. It is more important to resolve the variation of shear stress and heat flux since these exhibit an extremum inside the shock. In the case of $P_1$ basis functions, the solution derivative is constant inside each cell and so we plot the constant value of shear stress and heat flux at the center of each cell with a symbol. It is clear from these plots that both the schemes give accurate solutions and are able to resolve the shock structure well, even for the shear stress and heat flux. The DG solution compares very well with the analytical solution of the NS equation. With $P_2$ basis functions, we see that the DG solution is able to capture the variation of shear stress $\tau$ and heat flux $q$ inside each cell very accurately; the accuracy is seen to improve with grid refinement.

\section{Summary and conclusions}
We have proposed a DG scheme for Navier-Stokes equations which makes use of kinetic flux vector splitting to construct a numerical flux function for both the convective and diffusive fluxes. These schemes are motivated by applying them to scalar convection-diffusion problem and studying their accuracy. Energy/entropy stable schemes are obtained by adding additional stabilization terms which are constructed from the kinetic numerical flux and interior penalty terms. The resulting schemes which can be either symmetric or non-symmetric in the discretization of diffusive terms, exhibit good accuracy properties. The symmetric schemes show optimal convergence rates in numerical tests on the scalar convection-diffusion equation. These good properties are observed also for Navier-Stokes equations based on entropy variables.  Future work with these schemes will be focused on applying them for two dimensional viscous flow problems to study if the observed convergence rates carry over to higher dimensions, and trying to couple them with particle methods for rarefied flow situations.

\appendix

\section{1-D Navier-Stokes equations}
\label{sec:ns1d}

In conservation form the 1-D Navier-Stokes equations can be written as
\[
\frac{\partial U}{\partial t}+\frac{\partial F}{\partial x} + \frac{\partial G}{\partial x} = 0
\]
where
\[
U=\left[\begin{array}{c}
\rho\\
\rho u\\
\rho e
\end{array}\right],\qquad F=\left[\begin{array}{c}
\rho u\\
p+\rho u^{2}\\
(\rho e+p)u
\end{array}\right],
\qquad G=\left[\begin{array}{c}
0\\
-\tau\\
-\tau u+q
\end{array}\right]
\]
and the shear stress and heat flux are given by
\[
\tau=\frac{4}{3}\mu\frac{\partial u}{\partial x},\qquad q=-\kappa\frac{\partial T}{\partial x}
\]
The coefficient of heat conduction is given by$
\kappa=\frac{\mu C_{p}}{\Pr}=\frac{\mu\gamma R}{(\gamma-1)\Pr}$ where $R$ is the gas constant and $\Pr$ is the Prandtl number.
\section{Polyatomic gases}
\label{sec:poly}
The Maxwell-Boltzmann distribution function discussed above corresponds to a monatic gas for which the ratio of specific heats is $\gamma=\frac{5}{3}$. Real gases like air which are mostly composed of diatomic molecules $N_2$ and $O_2$ have a value of $\gamma=1.4$ and effectively behave like a diatomic gas. The internal energy of such gases is different since they have contributions from rotational degrees of freedom. In order to get the correct value of $\gamma$ and hence the internal energy, additional internal energy variables can be introduced. In the approach of Deshpande~\cite{smd1986a}, the collisional invariant corresponding to energy is modified to $\frac{1}{2}|v|^2+I$ while in the BGK schemes~\cite{Xu2001289}, it is modified as $\frac{1}{2}|v|^2+\xi^2$, with $\xi^2=\xi_1^2+\ldots+\xi_K^2$, the value of $K$ depending on $\gamma$. A third approach is to use $\frac{1}{2}|v|^2+I^\delta$ as introduced by Perthame~\cite{perthame1990}; this has the advantage of preserving the form of the Boltzmann entropy function $H=f\ln f$~\cite{perthame1990,Timothy20063311} and will be adopted in the present work. Note that the kinetic split fluxes resulting from any of the above approaches are identical since there is no splitting with respect to the additional degrees of freedom. With the above choice of the collisional invariants, the Maxwell-Boltzmann distribution function is given by~\cite{Timothy20063311}
\[
g(\rho,u,T; v, I) = \frac{\rho}{\alpha (RT)^{1/(\gamma-1)}} e^{-\frac{1}{RT}(\frac{1}{2}|v-u|^2+I^\delta)}, \qquad \alpha = \int_{\re^d} e^{-\frac{1}{2}|v|^2} \ud v \cdot \int_{\re^+} e^{-I^\delta} \ud I, \qquad \delta = \frac{1}{\frac{1}{\gamma-1} - \frac{d}{2}}
\]
The moments of the above distribution function yield the correct expressions for the energy and energy flux for a polyatomic gas while the mass and momentum equations are unaffected. The entropy variables for a polytropic gas corresponding to the convex entropy (\ref{eq:eulent}) are given by~\cite{Timothy20063311}
\begin{equation*}
V = \left[
-\frac{s}{\gamma-1} - \frac{|u|^2}{2RT}, \ \ 
\frac{u}{RT}, \ \ 
-\frac{1}{RT}
\right]^\top
\end{equation*}
Any constant terms in the entropy variables, e.g., in the definition of the physical entropy $s$, can be ignored since they do not change the entropy stability properties of the NS equations or of the DG schemes since only the jump in the entropy variables or their derivatives appear in the equations.

\section{KFVS fluxes for NS equation in 1-D}
\label{sec:kfvs1d}

The inviscid kinetic split fluxes are given by
\[
F^{\pm}=\left[\begin{array}{c}
\rho uA^{\pm}+\rho B^{\pm}\\
(p+\rho u^{2})A^{\pm}+\rho uB^{\pm}\\
(\rho e+p)uA^{\pm}+(\rho e+p/2)B^{\pm}
\end{array}\right]
\]
where
\[
s=u\sqrt{\beta},\qquad A^{\pm}=\frac{1}{2}(1\pm\text{erf}(s)),\qquad B^{\pm}=\pm\frac{1}{2\sqrt{\pi\beta}}\exp(-s^{2})
\]
The viscous kinetic split fluxes are given by
\[
G^{\pm}=\left[\begin{array}{c}
-\left[\tau + \frac{4}{5} \beta u q \right] \beta B^{\pm}\\
-\tau A^{\pm}+\frac{4}{5}\beta q B^{\pm}\\
(-u\tau+q)A^{\pm} - \left[\frac{3}{2}\tau + \frac{2}{5}\beta u q\right]B^{\pm}
\end{array}\right]
\]
where $\tau$ and $q$ are the shear stress and heat flux.
\section*{References}
\bibliographystyle{elsarticle-num}
\bibliography{bibdesk}

\begin{thebibliography}{10}
\expandafter\ifx\csname url\endcsname\relax
  \def\url#1{\texttt{#1}}\fi
\expandafter\ifx\csname urlprefix\endcsname\relax\def\urlprefix{URL }\fi
\expandafter\ifx\csname href\endcsname\relax
  \def\href#1#2{#2} \def\path#1{#1}\fi

\bibitem{cercignani1988boltzmann}
C.~Cercignani, The {B}oltzmann equation and its applications, Applied
  mathematical sciences, Springer-Verlag, 1988.

\bibitem{DI1980231}
D.~I. Pullin, Direct simulation methods for compressible inviscid ideal-gas
  flow, Journal of Computational Physics 34~(2) (1980) 231 -- 244.
\newblock \href {http://dx.doi.org/10.1016/0021-9991(80)90107-2}
  {\path{doi:10.1016/0021-9991(80)90107-2}}.

\bibitem{smd-nasatp-2613}
S.~M. Deshpande, A second-order accurate kinetic theory based method for
  inviscid compressible flows, Tech. Rep. NASA TP-2613, NASA (1986).

\bibitem{Mandal1994447}
J.~C. Mandal, S.~M. Deshpande, Kinetic flux vector splitting for {E}uler
  equations, Computers and Fluids 23~(2) (1994) 447 -- 478.
\newblock \href {http://dx.doi.org/DOI: 10.1016/0045-7930(94)90050-7}
  {\path{doi:DOI: 10.1016/0045-7930(94)90050-7}}.

\bibitem{PhysRev.94.511}
P.~L. Bhatnagar, E.~P. Gross, M.~Krook, A model for collision processes in
  gases. {I}. {S}mall amplitude processes in charged and neutral one-component
  systems, Phys. Rev. 94 (1954) 511--525.
\newblock \href {http://dx.doi.org/10.1103/PhysRev.94.511}
  {\path{doi:10.1103/PhysRev.94.511}}.

\bibitem{Chou:1997:KFS:254115.254120}
S.~Y. Chou, D.~Baganoff, Kinetic flux-vector splitting for the
  {N}avier-{S}tokes equations, J. Comput. Phys. 130 (1997) 217--230.
\newblock \href {http://dx.doi.org/10.1006/jcph.1996.5579}
  {\path{doi:10.1006/jcph.1996.5579}}.

\bibitem{Prendergast:1993:NHG:167236.167246}
K.~H. Prendergast, K.~Xu, Numerical hydrodynamics from gas-kinetic theory, J.
  Comput. Phys. 109 (1993) 53--66.
\newblock \href {http://dx.doi.org/10.1006/jcph.1993.1198}
  {\path{doi:10.1006/jcph.1993.1198}}.

\bibitem{xu_martinelli_jameson_1995}
K.~Xu, L.~Martinelli, A.~Jameson, Gas-kinetic finite volume methods, in: S.~M.
  Deshpande, S.~Desai, R.~Narasimha (Eds.), Fourteenth International Conference
  on Numerical Methods in Fluid Dynamics, Vol. 453 of Lecture Notes in Physics,
  Springer, Berlin, 1995, pp. 106--111.

\bibitem{Xu2001289}
K.~Xu, A gas-kinetic {BGK} scheme for the {N}avier-{S}tokes equations and its
  connection with artificial dissipation and {G}odunov method, Journal of
  Computational Physics 171~(1) (2001) 289 -- 335.
\newblock \href {http://dx.doi.org/DOI: 10.1006/jcph.2001.6790} {\path{doi:DOI:
  10.1006/jcph.2001.6790}}.

\bibitem{Ohwada2002156}
T.~Ohwada, On the construction of kinetic schemes, Journal of Computational
  Physics 177~(1) (2002) 156 -- 175.

\bibitem{torr_xu}
M.~Torrilhon, K.~Xu, Stability and consistency of kinetic upwinding for
  advection--diffusion equations, IMA Journal of Numerical Analysis 26~(4)
  (2006) 686--722.

\bibitem{reedhill}
W.~H. Reed, T.~R. Hill, Triangular mesh methods for the neutron transport
  equation, Tech. Rep. LA-UR-73-476, Los Alamos Scientific Laboratory (1973).

\bibitem{Cockburn:1989:TRL:69978.69982}
B.~Cockburn, S.-Y. Lin, C.-W. Shu, {TVB} {R}unge-{K}utta local projection
  discontinuous {G}alerkin finite element method for conservation laws {III}:
  {O}ne-dimensional systems, J. Comput. Phys. 84 (1989) 90--113.

\bibitem{Cockburn:1998:RDG:287244.287254}
B.~Cockburn, C.-W. Shu, The {R}unge-{K}utta discontinuous {G}alerkin method for
  conservation laws {V}: {M}ultidimensional systems, J. Comput. Phys. 141
  (1998) 199--224.

\bibitem{arnold2002}
D.~N. Arnold, F.~Brezzi, B.~Cockburn, L.~D. Marini, Unified analysis of
  discontinuous {G}alerkin methods for elliptic problems, SIAM Journal on
  Numerical Analysis 39~(5) (2002) pp. 1749--1779.

\bibitem{Bassi1997267}
F.~Bassi, S.~Rebay, A high-order accurate discontinuous finite element method
  for the numerical solution of the compressible {N}avier-{S}tokes equations,
  Journal of Computational Physics 131~(2) (1997) 267 -- 279.
\newblock \href {http://dx.doi.org/10.1006/jcph.1996.5572}
  {\path{doi:10.1006/jcph.1996.5572}}.

\bibitem{Cockburn:1998:LDG:305653.305671}
B.~Cockburn, C.-W. Shu, The local discontinuous {G}alerkin method for
  time-dependent convection-diffusion systems, SIAM J. Numer. Anal. 35 (1998)
  2440--2463.

\bibitem{arnold1982}
D.~N. Arnold, An interior penalty finite element method with discontinuous
  elements, SIAM Journal on Numerical Analysis 19~(4) (1982) pp. 742--760.

\bibitem{Oden1998491}
J.~Oden, I.~Babuska, C.~E. Baumann, A discontinuous hp finite element method
  for diffusion problems, Journal of Computational Physics 146~(2) (1998) 491
  -- 519.

\bibitem{Baumann1999311}
C.~E. Baumann, J.~T. Oden, A discontinuous hp finite element method for
  convection-diffusion problems, Computer Methods in Applied Mechanics and
  Engineering 175~(3-4) (1999) 311 -- 341.

\bibitem{FLD:FLD338}
F.~Bassi, S.~Rebay, Numerical evaluation of two discontinuous {G}alerkin
  methods for the compressible {N}avier--{S}tokes equations, International
  Journal for Numerical Methods in Fluids 40~(1-2) (2002) 197--207.
\newblock \href {http://dx.doi.org/10.1002/fld.338}
  {\path{doi:10.1002/fld.338}}.

\bibitem{Hartmann20089670}
R.~Hartmann, P.~Houston, An optimal order interior penalty discontinuous
  {G}alerkin discretization of the compressible {N}avier-{S}tokes equations,
  Journal of Computational Physics 227~(22) (2008) 9670 -- 9685.
\newblock \href {http://dx.doi.org/10.1016/j.jcp.2008.07.015}
  {\path{doi:10.1016/j.jcp.2008.07.015}}.

\bibitem{Gassner:2008:DGS:1342061.1342078}
G.~Gassner, F.~L\"{o}rcher, C.~D. Munz,
  \href{http://dx.doi.org/10.1007/s10915-007-9169-1}{A discontinuous {G}alerkin
  scheme based on a space-time expansion {II}. {V}iscous flow equations in
  multi dimensions}, J. Sci. Comput. 34~(3) (2008) 260--286.
\newblock \href {http://dx.doi.org/10.1007/s10915-007-9169-1}
  {\path{doi:10.1007/s10915-007-9169-1}}.
\newline\urlprefix\url{http://dx.doi.org/10.1007/s10915-007-9169-1}

\bibitem{Xu:2004:DGB:996978.997050}
K.~Xu, Discontinuous {G}alerkin {BGK} method for viscous flow equations:
  {O}ne-dimensional systems, SIAM J. Sci. Comput. 25 (2004) 1941--1963.
\newblock \href {http://dx.doi.org/10.1137/S1064827502416113}
  {\path{doi:10.1137/S1064827502416113}}.

\bibitem{springerlink:10.1007/BF03177419}
H.~Liu, K.~Xu, A gas-kinetic discontinuous {G}alerkin method for viscous flow
  equations, Journal of Mechanical Science and Technology 21 (2007) 1344--1351,
  10.1007/BF03177419.

\bibitem{Amiram1983151}
A.~Harten, On the symmetric form of systems of conservation laws with entropy,
  Journal of Computational Physics 49~(1) (1983) 151 -- 164.
\newblock \href {http://dx.doi.org/10.1016/0021-9991(83)90118-3}
  {\path{doi:10.1016/0021-9991(83)90118-3}}.

\bibitem{Shakib1991141}
F.~Shakib, T.~J. Hughes, Z.~Johan, A new finite element formulation for
  computational fluid dynamics: X. {T}he compressible {E}uler and
  {N}avier-{S}tokes equations, Computer Methods in Applied Mechanics and
  Engineering 89~(1-3) (1991) 141 -- 219.
\newblock \href {http://dx.doi.org/10.1016/0045-7825(91)90041-4}
  {\path{doi:10.1016/0045-7825(91)90041-4}}.

\bibitem{smd1986a}
S.~M. Deshpande, On the {M}axwellian distribution, symmetric form and entropy
  conservation for the {E}uler equations, Tech. Rep. TP-2583, NASA Langley,
  Hampton, VA (1986).

\bibitem{barth1998}
T.~Barth, Numerical methods for gasdynamic systems on unstructured grids, in:
  Kroner, Ohlberger, Rohde (Eds.), An introduction to recent developments in
  theory and numerics for conservation laws, Vol.~5 of Lecture notes in
  computational science and engineering, Springer-Verlag, 1998, pp. 198--285.

\bibitem{Timothy20063311}
T.~Barth, On discontinuous {G}alerkin approximations of {B}oltzmann moment
  systems with {L}evermore closure, Computer Methods in Applied Mechanics and
  Engineering 195~(25-28) (2006) 3311 -- 3330.
\newblock \href {http://dx.doi.org/10.1016/j.cma.2005.06.016}
  {\path{doi:10.1016/j.cma.2005.06.016}}.

\bibitem{osher:947}
S.~Osher, Convergence of generalized {MUSCL} schemes, SIAM Journal on Numerical
  Analysis 22~(5) (1985) 947--961.
\newblock \href {http://dx.doi.org/10.1137/0722057}
  {\path{doi:10.1137/0722057}}.

\bibitem{Hauke:1998:0045-7825:1}
G.~Hauke, T.~J.~R. Hughes, A comparative study of different sets of variables
  for solving compressible and incompressible flows, Computer Methods in
  Applied Mechanics and Engineering 153~(1) (1998) 1--44.
\newblock \href {http://dx.doi.org/doi:10.1016/S0045-7825(97)00043-1}
  {\path{doi:doi:10.1016/S0045-7825(97)00043-1}}.

\bibitem{Shu1988439}
C.-W. Shu, S.~Osher, Efficient implementation of essentially non-oscillatory
  shock-capturing schemes, Journal of Computational Physics 77~(2) (1988) 439
  -- 471.
\newblock \href {http://dx.doi.org/10.1016/0021-9991(88)90177-5}
  {\path{doi:10.1016/0021-9991(88)90177-5}}.

\bibitem{cockburn1989}
B.~Cockburn, C.-W. Shu, {TVB} {R}unge-{K}utta local projection discontinuous
  {G}alerkin finite element method for conservation laws {II}: {G}eneral
  framework, Mathematics of Computation 52~(186) (1989) pp. 411--435.

\bibitem{liepmann-roshko}
H.~W. Liepmann, A.~Roshko, Elements of {G}asdynamics, Dover, 2002.

\bibitem{perthame1990}
B.~Perthame, Boltzmann type schemes for gas dynamics and the entropy property,
  SIAM Journal on Numerical Analysis 27~(6) (1990) pp. 1405--1421.

\bibitem{barthcharrier2001}
T.~J. Barth, P.~Charrier, Energy stable flux formulas for the discontinuous
  {G}alerkin discretization of first-order nonlinear conservation laws, Tech.
  Rep. NAS-01-001, NASA (2001).

\bibitem{May:2007:IGB:1223679.1223717}
G.~May, B.~Srinivasan, A.~Jameson, An improved gas-kinetic {BGK} finite-volume
  method for three-dimensional transonic flow, J. Comput. Phys. 220 (2007)
  856--878.
\newblock \href {http://dx.doi.org/10.1016/j.jcp.2006.05.027}
  {\path{doi:10.1016/j.jcp.2006.05.027}}.

\bibitem{dutt1988}
P.~Dutt, Stable boundary conditions and difference schemes for
  {N}avier-{S}tokes equations, SIAM Journal on Numerical Analysis 25~(2) (1988)
  pp. 245--267.

\bibitem{gilbarg1953}
D.~Gilbarg, D.~Paolucci, The structure of shock waves in the continuum theory
  of fluids, Journal of Rational Mechanics and Analysis 2 (1953) 617.

\end{thebibliography}

\end{document}